\documentclass[a4paper,10pt]{article}
 
 \usepackage{amsmath, amssymb, latexsym, amscd, amsthm,amsfonts,amstext}
 \usepackage[mathscr]{eucal}
 \usepackage{graphicx}
 \usepackage{subfig}

\title{Reconstruction of shapes
and refractive indices from backscattering experimental data using the
adaptivity} 

\author{Larisa Beilina$^{\circ}$, Nguyen Trung Th\`anh$^{\diamond}$, Michael V. Klibanov$^{\diamond}$\footnote{Corresponding author} \\
and John Bondestam Malmberg$^{\circ}$\\
$^{\circ}$Department of Mathematical Sciences, \\
Chalmers University of Technology and
Gothenburg University, Gothenburg, Sweden \\
Emails: \texttt{larisa@chalmers.se, john.bondestam.malmberg@chalmers.se}.\\
$^\diamond$Department of Mathematics \& Statistics, University of North Carolina at Charlotte,\\
Charlotte, NC, USA. Emails: \texttt{tnguy152@uncc.edu, mklibanv@uncc.edu}}
\date{}

\begin{document}
\maketitle

\begin{abstract}
We consider the inverse problem of the reconstruction of
the spatially distributed dielectric constant $\varepsilon _{r}\left( 
\mathbf{x}\right), \ \mathbf{x}\in \mathbb{R}^{3}$, which is an unknown
coefficient in the Maxwell's equations, from time-dependent backscattering experimental radar data associated with a single source of electric pulses. The refractive index is $n\left( 
\mathbf{x}\right) =\sqrt{\varepsilon _{r}\left( \mathbf{x}\right) }.$ The
coefficient $\varepsilon _{r}\left( \mathbf{x}\right) $ is reconstructed
using a two-stage reconstruction procedure. In the first
stage an approximately globally convergent method proposed is applied
to get a good first approximation of the exact solution. In the second stage a locally
convergent adaptive finite element method is applied, taking the solution of the first stage as the starting point of the minimization of the Tikhonov functional. This functional is minimized on a
sequence of locally refined meshes. It is shown here that all three components of interest
of targets can be simultaneously accurately imaged: refractive indices,
shapes and locations.
\end{abstract}


\textbf{AMS classification codes:} 65N15, 65N30, 35J25.

\def\bR{\mathbb{R}}
\def\bx{\mathbf{x}}

\graphicspath{{FIGURES/}
{Figures/}
{FiguresJ/newfigures/}
{pics/}} 

\section{Introduction}
\label{sec:1}

In this paper we investigate the problem of imaging objects placed in air from time-dependent backscattering radar measurements, using a two-stage reconstruction method. In the first stage, initial images are calculated using the globally
convergent method for Coefficient Inverse Problems (CIPs),
which was originated in \cite{BK1} with a number of follow up publications;
results were summarized in the book \cite{BK}. In the second stage, those images are refined
using an adaptive finite element method (adaptivity) of \cite{BMaxwell2}. Results of the first stage for the data sets considered in this paper were presented in \cite{BTKF,NBKF}. Here, we present the results of the second stage. Only the maximum value of the dielectric constant and the location of a target were
accurately reconstructed in \cite{BTKF,NBKF} using the globally convergent method. The accuracy of the reconstruction of the shape of the target was limited. Using the two-stage reconstruction procedure, it is shown here that we can simultaneously and accurately reconstruct all three components  of interest of objects: refractive indices, shapes, and locations.

\begin{figure}[tbp]
\centering
\subfloat[]{\includegraphics[width = 0.4\textwidth,height=0.3\textwidth]{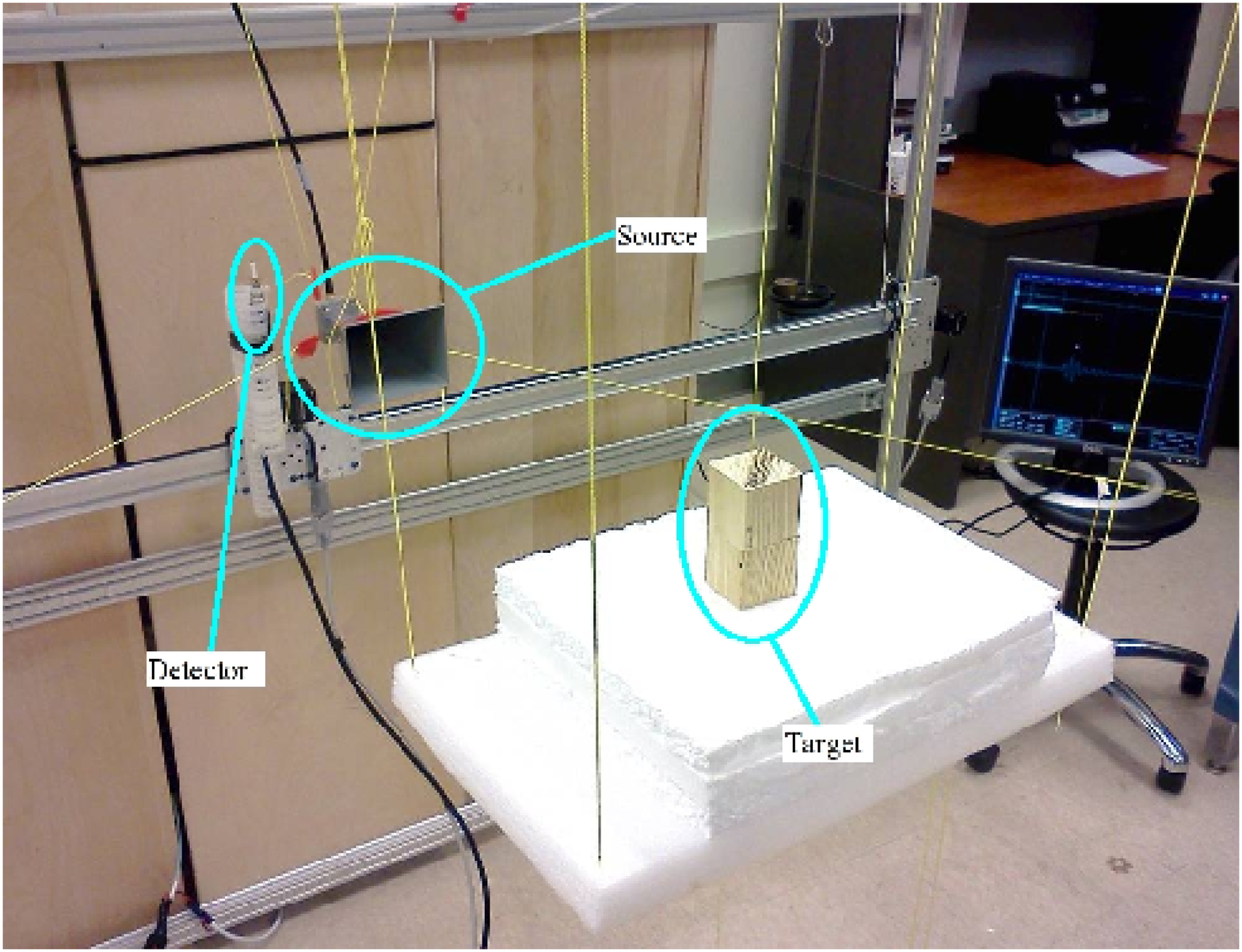}} \hspace{%
0.3truecm} 
\subfloat[]{\includegraphics[width =
0.34\textwidth,height=0.3\textwidth]
{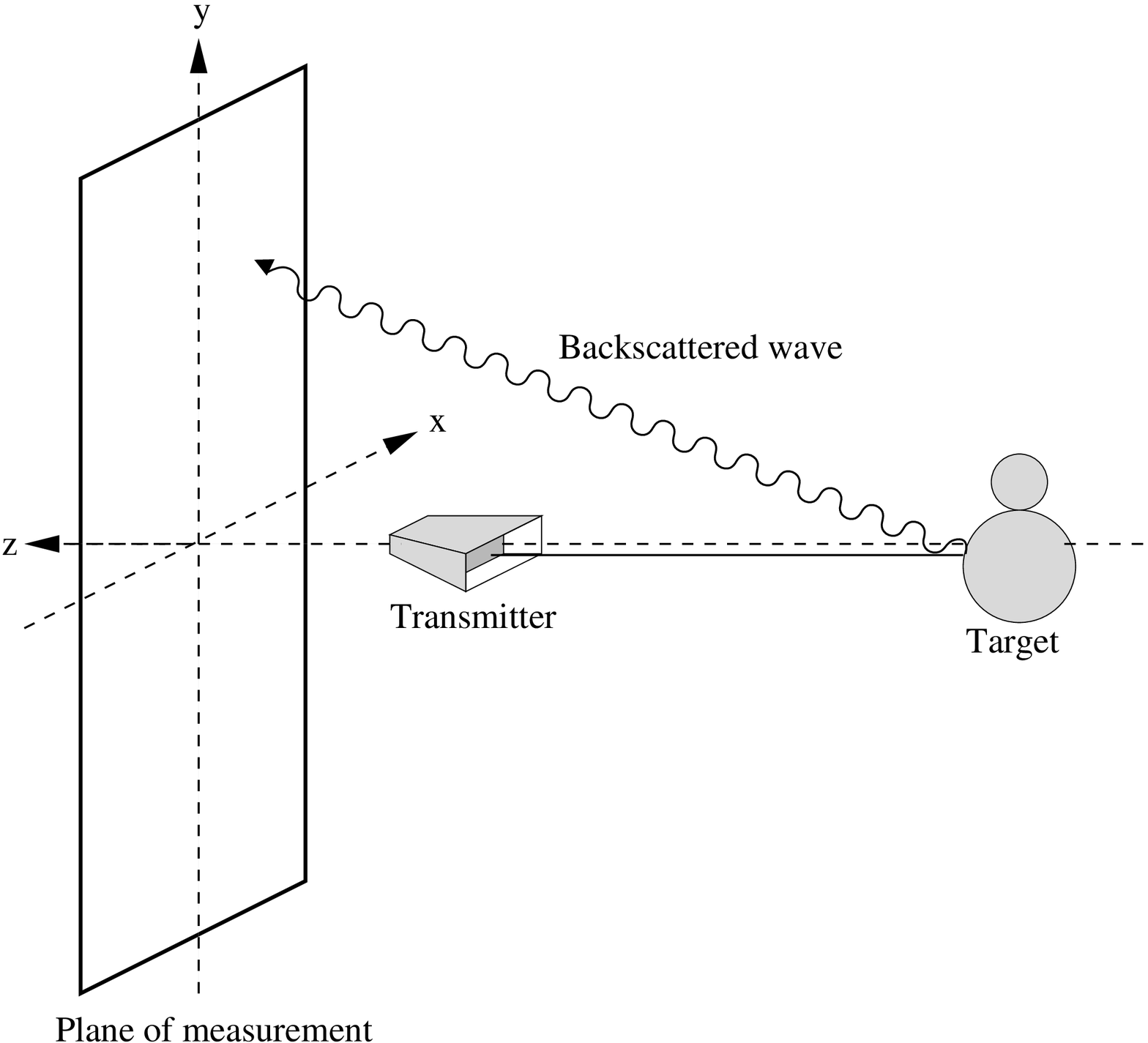}
} 
\caption{(a): Experimental setup; (b) Our data collection  scheme.}
\label{fig:setup}
\end{figure}

We reconstruct these three components simultaneously as parts of an
unknown coefficient, which is the spatially varying dielectric constant
$\varepsilon _{r}\left( \mathbf{x}\right),\ \bx\in \mathbb R^3, $ in the Maxwell's
equations. 
Below $\mathbf{x}=\left(
x,y,z\right) \in \mathbb{R}^{3},$ where $x$ is the horizontal axis,
$y$ is the vertical axis and $z$ is the axis which points from the
target towards the measurement plane, see Figure \ref{fig:setup}. Even though only one component $E_{2}$ of
the electric field $E=\left( E_{1},E_{2},E_{3}\right) $ was measured by our
experimental device, we numerically solve here a CIP for the three-dimensional (3-d) Maxwell's
equations. The boundary data for two other components $E_{1},E_{3}$ are obtained via computational simulations. 

Experimental data were collected by a microwave device which was
recently assembled at the University of North Carolina at Charlotte,
USA. Our desired application is imaging of explosives. In this paper
we consider only targets located in air. The work on real data for the
case when targets are buried under the ground is reported in \cite{TBKF:2013-2}. Note that
explosives may be located in air \cite{KBKSNF}, e.g., improvised
explosive devices (IEDs). We image both homogeneous and heterogeneous targets. Heterogeneous targets model IEDs.

 To collect those data, a single location of the source of electric pulses was
 used. Hence, we used the minimal amount of the information. The
 use of more sources was both hard to arrange experimentally and
 undesirable for our target application. The backscattering time dependent signal was
 measured at a number of detectors covering a part of a plane,
 i.e., over a narrow range of backscattering angles, see Figure
 \ref{fig:setup}. That plane was placed behind the source. Experimental
 data of this paper were collected for targets located in air on the
 distance of 80 centimeters from the measurement plane, which is 20
 wavelengths, i.e., in a far field zone. The distance between
 neighboring detectors was 2 centimeters. 


We refer to, e.g., \cite{Lak1,Lak2} for treatments of experimental data in
the frequency domain by other numerical methods for CIPs for Maxwell's
equations. In particular, blind real data were considered in \cite{Lak2}. As
to the adaptivity technique for inverse problems, we refer to, e.g., \cite%
{BJ,Li}. There are many works dedicated to inverse problems of shape
reconstruction; we refer to some most recent ones, e.g., \cite%
{Liu1,Liu2,T1,T2}. We also refer to \cite{Isakov} for a survey about inverse
problems of shape reconstruction.

An outline of this paper is as follows. In section \ref{sec:stages} we
describe the two-stage reconstruction procedure. In section \ref{sec:3} we state
the forward and inverse problems. In section \ref{sec:invprobl} we present
Tikhonov functional and optimality conditions. In section \ref{sec:fem} we
describe the finite element method used in our computation. In section \ref%
{subsec:ad_alg} we present the mesh refinement recommendation and the
adaptive algorithm. Some details of numerical implementation are described in section \ref{sec:numex}. In section \ref{sec:8.2} we present reconstruction results. Finally, a summary is given in section \ref{sec:summary}.

\section{Two-stage reconstruction procedure}

\label{sec:stages}

In \cite{BTKF,NBKF} we have considered the problem of the reconstruction of
the spatially distributed dielectric constant $\varepsilon _{r}(\mathbf{x}),%
\mathbf{x}\in \mathbb{R}^{3}$ from experimental data, which were the
same as in the current paper. In \cite{BTKF,NBKF} this function was the
unknown coefficient in a wave-like PDE 
\begin{equation}
\varepsilon _{r}(\mathbf{x})\frac{ \partial ^2 E_{2}}{\partial t^2} = \Delta E_{2}.  \label{1}
\end{equation}%
To reconstruct $\varepsilon _{r}(\mathbf{x}),$ we have used the
approximately globally convergent algorithm of \cite{BK}. The notion of the
approximate global convergence (\textquotedblleft global convergence" in
short) was introduced in \cite{BK,BKJIIP12}. Indeed, conventional least
squares cost functionals for CIPs are non convex and typically have many
local minima. Hence, given a CIP, the first question to address in its
numerical treatment is: \emph{How to obtain a good approximation for the
exact solution without any a priori knowledge of a small neighborhood of
this solution?} We call a numerical method addressing this question \emph{%
globally convergent}.

It is well known that it is tremendously difficult to address this question.
For this reason, a certain reasonable approximation was made in \cite%
{BK,BKJIIP12}.\ This approximation is used only on the first iteration of
that method. Because of this approximation, we call the technique of \cite%
{BK,BKJIIP12} approximately globally convergent. Due to that approximation,
a room is left for a refinement of results.


An important point here is that there exists a rigorous guarantee
within the framework of that approximate model that the solution
resulting from the globally convergent method is located in a small
neighborhood of the exact solution. This is achieved without any a
priori knowledge of a small neighborhood of the exact solution.  Thus,
a locally convergent numerical method can be used for a refinement of
the solution obtained by the globally convergent technique of
\cite{BK,BKJIIP12}.

The latter is the main goal of the current paper. We synthesize here the
adaptive finite element method of \cite{BMaxwell2} with the globally
convergent numerical method of \cite{BK} in order to improve the
reconstruction of shapes of objects imaged in \cite{BTKF,NBKF}.
 The idea of
this synthesis was first introduced in \cite{BK2}. The synthesis represents
the following two-stage reconstruction procedure:

\textbf{Stage 1}. In this stage the approximately globally convergent method
of \cite{BK} is applied and a good first approximation for the exact
solution is obtained.

\textbf{Stage 2}. This stage refines the solution obtained in the
first stage. The locally convergent adaptivity technique of
\cite{BMaxwell2} is applied. The solution obtained in \cite{BTKF,NBKF}
in the first stage is taken as the starting point in the minimization
procedure of the Tikhonov functional.

An important advantage of using the two-stage reconstruction procedure follows
from Theorem 1.9.1.2 of \cite{BK}. This theorem states that the minimizer
of the Tikhonov functional (i.e., the regularized solution) is closer to the
exact solution than the first guess, provided, however, that the first guess
is sufficiently close to the exact solution. Therefore, that first guess
should be delivered by the globally convergent method of the first stage.
The adaptive finite element method of \cite{BMaxwell2}, which we use in the
second stage, minimizes the Tikhonov functional on a sequence of locally
refined meshes, which is the main attractive point of the adaptivity. It
enables one to maintain a reasonable compromise between not using an
exceedingly large number of finite elements and a good accuracy of the resulting
solution. It follows from Theorem 4.9.3 of \cite{BK} and Theorem 5.2 of \cite{BKrelax}
that if the first guess is sufficiently close to the exact solution, then the accuracy of its
reconstruction monotonically improves with local mesh refinements. On the other hand, it
was shown in section 5.8.4 of the book \cite{BK} that a locally convergent
numerical method taken alone does not work for transmitted time dependent
experimental data generated by a single source. The same conclusion was
drawn for a different type of experimental data in \cite{Liu}.

\section{Statement of Forward and Inverse Problems}

\label{sec:3}

We model the electromagnetic wave
propagation in an isotropic, non-magnetic space $\mathbb R^3$ with the dielectric constant coefficient $\varepsilon_r(\bx)$. The electric field $E = (E_1,E_2,E_3)$ satisfies the following Cauchy problem:
\begin{equation}\label{E_gauge1}
\begin{split}
\varepsilon_r(\bx) \frac{\partial^2 E}{\partial t^2} + \nabla \times (\nabla \times E) 
 &=  (0,\delta(z-z_0)f(t),0),~ \mbox{in}~~ \mathbb{R}^{3} \times (0,T) , \\
  E(\mathbf{x},0) = 0, ~~~E_t(\mathbf{x},0) &=0~ \mbox{in}~~  \mathbb{R}^{3}.   
\end{split}
\end{equation}
where $f\left( t\right) \not\equiv 0$ is the time-dependent waveform
of the component $E_{2}$ of the incident plane wave, which is
originated at the plane $\left\{ z=z_{0}\right\} $ and propagates
along the $z$-axis. In our experiment the component $E_{2}$
corresponds to the electromagnetic pulse which is sent into the
medium. Thus, in (\ref{E_gauge1}) as well as in our computer simulations of
section \ref{sec:numex} the incident field has only one non-zero
component $E_{2}\left( \mathbf{x},t\right) $. This component
propagates along the $z$-axis until it reaches the target, where it is
scattered.


Let $\Omega \subset \mathbb{R}^{3}$ be a convex bounded domain. 
We impose the following conditions on the coefficient $\varepsilon _{r}\left( \mathbf{x}\right)$: 
\begin{equation}
\varepsilon _{r}(\mathbf{x})\in C^{\alpha }\left( \mathbb{R}^{3}\right)
,\varepsilon _{r}(\mathbf{x})\in \lbrack 1,d],~~\varepsilon _{r}(\mathbf{x}%
)=1\text{ for }\mathbf{x}\in \mathbb{R}^{3}\diagdown \Omega ,  \label{2.20}
\end{equation}
where $d=const.>1.$ We a priori assume the knowledge of the
constant $d$.  This means the knowledge of the set of admissible
coefficients in (\ref{2.20}). However, we do not impose small-value
assumptions on the unknown coefficient $\varepsilon
_{r}(\mathbf{x})$, i.e., we do not assume that $d$ is small. 
 Here $C^{\alpha },\alpha \in
\left( 0,1\right),$ is the H\"{o}lder space. Let $\Gamma \subset \partial
\Omega $ be the part of the boundary $\partial \Omega $ on which the backscattered data are measured.

\textbf{Coefficient Inverse Problem (CIP):} \emph{Suppose that the
  coefficient }$\varepsilon _{r}\left( \mathbf{x}\right)
$\emph{ satisfies conditions (\ref{2.20})  and that $\overline{\Omega }\cap \left\{
z=z_{0}\right\} =\varnothing$.  Determine the function
}$\varepsilon _{r}\left( \mathbf{x}\right) $\emph{\ for
}$\mathbf{x}\in \Omega $\emph{, assuming that the following function
}$g(\mathbf{x},t)$\emph{\ is known for a single incident plane wave}
\begin{equation}
g\left( \mathbf{x},t\right) = E\left( \mathbf{x},t\right) ,\forall \left( 
\mathbf{x},t\right) \in \Gamma \times \left( 0,\infty \right) .  \label{2.5}
\end{equation}


\subsection{Domain decomposition  finite element/finite difference method}

It is impossible to solve the problem (\ref{E_gauge1}) in the whole space 
$\mathbb{R}^{3}$. Hence, we solve it in a bounded domain $%
G$ which contains our domain of interest $\Omega$. For the convenience of our local mesh refinement procedure, we use the domain decomposition finite element/finite difference
method of \cite{BM}. For this purpose, we decompose $G$ as $G =\Omega _{FEM}\cup \Omega _{FDM}$ with $\Omega _{FEM}=\Omega$. Here we use a finite element mesh in $\Omega _{FEM}$, and in $\Omega _{FDM}$ we use a finite
difference mesh. These two domains have a thin layer of structured overlapping
nodes where we use an exchange procedure between computational solutions
obtained by finite element and finite difference methods, see details in 
\cite{BM}.
By (\ref{2.20}) 
\begin{equation}
\begin{split}
\varepsilon _{r}(\mathbf{x})& \geq 1,\text{ for }\mathbf{x}\in \Omega _{FEM},
\\
\varepsilon _{r}(\mathbf{x})& =1,\text{ for }\mathbf{x}\in \Omega _{FDM}.
\end{split}
\label{coefic}
\end{equation}

In our computation we use the following model problem in the computational
domain $G$ for the electric field $E$ with the
stabilizing divergence condition \cite{Ass} with $s\geq 1$ and with boundary
conditions specified in this section below:
\begin{eqnarray}
\varepsilon_r \frac{\partial^2 E}{\partial t^2} +
 \nabla \times ( \nabla \times E)  - s\nabla  ( \nabla \cdot(\varepsilon_r E))   &=& 0,~ \mbox{in}~~ G \times (0,T),    \label{model1_1} \\
  E(\mathbf{x},0) = 0, ~~~E_t(\mathbf{x},0) &=& 0~ \mbox{in}~~  G.      \label{model1_3}  
\end{eqnarray}
For simplicity, we choose the domains $\Omega$ and $G$ by
\begin{equation}
\Omega =\Omega _{FEM}=\left\{ \mathbf{x}=\left( x,y,z\right):-a<x<a,-b<y<b,-c<z<c_{1}\right\} ,  \label{4}
\end{equation}
\begin{equation}
G=\left\{ \mathbf{x}=\left( x,y,z\right) :-X<x<X,-Y<y<Y,-Z<z<z_{0}\right\} ,
\label{5}
\end{equation}
with positive numbers $a,b,c$, $X>a,\ Y>b,\ -Z < -c < c_1 < z_{0}$ and $\Omega
_{FDM}=G\diagdown \Omega _{FEM}.$ Denote by
\begin{equation}
\partial_{1}G :=\overline{G}\cap \left\{ z=z_{0}\right\} ,\quad \partial_{2}G :=%
\overline{G}\cap \left\{ z=-Z\right\} ,\quad \partial_{3} G :=\partial G\diagdown
\left( \partial _{1}G\cup \partial _{2}G\right) .  \label{6}
\end{equation}
The backscattering side of $\Omega$ is $\Gamma =\partial \Omega \cap \left\{ z=c_{1}\right\}.$ Next, define $%
\partial _{i}G_{T} :=\partial _{i}G\times \left( 0,T\right) ,i=1,2,3.$ Let $%
t_{1}\in \left( 0,T\right) $ be a number, and we assume that the function $f\left( t\right)
\in C\left[ 0,t_{1}\right]$ and $f(t) = 0, $ for $t > t_1$.  We impose the following boundary conditions%
\begin{equation}
E\left( \textbf{x},t\right) =  (0,f(t),0) \text{ }\mbox{on}\text{ }\partial
_{1}G\times \left( 0,t_{1}\right] ,  \label{70}
\end{equation}
\begin{equation}
\partial _{n} E(\mathbf{x},t)=-\partial _{t} E(\mathbf{x},t)~\mbox{on}%
~~\partial _{1}G\times \left( t_{1},T\right) ,  \label{7}
\end{equation}
\begin{equation}
\partial _{n} E (\mathbf{x},t)=-\partial _{t} E(\mathbf{x},t)~\mbox{on}%
~~\partial _{2}G_{T},  \label{8}
\end{equation}
\begin{equation}
\partial _{n} E(\mathbf{x},t)=0~\mbox{on}~~\partial _{3}G_{T},  \label{9}
\end{equation}
where $\partial _{n}$ is the normal derivative. Conditions (\ref{7}),
(\ref{8}) are first order absorbing boundary conditions \cite{EM} at the planes $\partial _{1}G$ and $\partial _{2}G$ of the
rectangular prism $G$, and (\ref{9}) is the zero Neumann condition at
the lateral part $\partial _{3}G$ of the boundary $\partial G.$
Condition (\ref{70}) means that the incident plane wave is emitted only up to the time $t =t_1$ and then propagates inside of the domain $G$.

It was demonstrated numerically in \cite{BM} that the solution
of the problem (\ref{model1_1})--(\ref{9}) approximates well the solution of the original
Maxwell's equations for $s=1.$ The energy estimate of Theorem 4.1 of \cite{BM}
guarantees the stability of the forward problem  (\ref{model1_1})--(\ref{9}) for $s\geq 1$.

Using the transformation $\nabla \times (\nabla \times E) = \nabla
(\nabla \cdot E) - \nabla \cdot (\nabla E)$, the model problem
(\ref{model1_1}), (\ref{model1_3}),  (\ref{70}) -- (\ref{9}) can be rewritten  
as
\begin{eqnarray}
\varepsilon_r \frac{\partial^2 E}{\partial t^2} + \nabla (\nabla \cdot E) 
- \nabla \cdot (\nabla E) - s\nabla  ( \nabla \cdot(\varepsilon_r E))  &=& 0,~ \mbox{in}~~ G\times (0,T),    \label{model3_1} \\
  E(\mathbf{x},0) = 0, ~~~E_t(\mathbf{x},0) &=&0 ~ \mbox{in}~~ G,      \label{model3_2}  \\
E\left( \textbf{x},t\right) &= &(0,f\left( t\right),0) \text{ }\mbox{on}\text{ }\partial
_{1}G\times \left( 0,t_{1}\right] , \label{model3_3}  \\
\partial _{n} E(\mathbf{x},t) &=&-\partial _{t} E(\mathbf{x},t)~\mbox{on}~~\partial _{1}G\times \left( t_{1},T\right) ,   \label{model3_4} \\
\partial _{n} E (\mathbf{x},t) &=& -\partial _{t} E(\mathbf{x},t)~\mbox{on}~~\partial _{2}G_{T},  \label{model3_5}\\
\partial _{n} E(\mathbf{x},t) &=& 0~\mbox{on}~~\partial _{3}G_{T}. \label{model3_6}
\end{eqnarray}
 We refer to \cite{BM} for details of the
numerical solution of the forward problem (\ref{model3_1})-(\ref{model3_6}).

\section{Tikhonov functional and optimality conditions}
\label{sec:invprobl}

Let $\Gamma _{1}$ be the extension of the backscattering side $\Gamma $ up to the boundary $\partial _{3}G$ of the domain $G,$ i.e.,
\begin{equation}
\Gamma _{1}=\left\{ \mathbf{x}=\left( x,y,z\right)
: - X<x< X,- Y<y< Y, z=c_{1}\right\} .  \label{10}
\end{equation}
Let $G_{b}$ be the part of the rectangular prism $G$ which lies between the two planes $\Gamma _{1}$ and $\{z = -Z\}$:
\begin{equation}
G_{b}=\left\{ \mathbf{x}=\left( x,y,z\right)
:- X<x< X,- Y<y< Y,- Z<z<c_{1}\right\} .  \label{11}
\end{equation}
Denote by $Q_{T}=G_{b}\times \left( 0,T\right) ,S_{T}=\partial
G_{b}\times \left( 0,T\right) .$ Even though we have the data $g\left(
\mathbf{x} ,t\right) $ in (\ref{2.5}) only on $\Gamma ,$ we show in
subsection \ref{subsec:immersing} below how we complement these
data on the rest of the boundary $\partial G_{b}$ of the domain
$G_{b},$ i.e., on $\partial G_{b}\diagdown \Gamma .$ This way we
approximately obtain the function $\widetilde{g}\left(
\mathbf{x},t\right)$:
\begin{equation}
\widetilde{g}\left( \mathbf{x},t\right) = E\left( \mathbf{x},t\right)
,\forall \left( \mathbf{x},t\right) \in S_{T}.  \label{12}
\end{equation}%
We reformulate our inverse problem as an optimization problem.
Thus, we find $\varepsilon_r$ by minimizing  the Tikhonov functional:
\begin{equation}
F(E, \varepsilon_r) := \frac{1}{2} \int_{S_T}(E - \tilde{g})^2 z_{\delta }(t) d\bx dt +
\frac{1}{2} \gamma \int_{G}( \varepsilon_r -  {\varepsilon_r}_{glob})^2~~ d\bx ,
\label{functional}
\end{equation}
where $\gamma > 0$ is the regularization parameter, and $\varepsilon _{r,glob}\left( \mathbf{x}\right) $ is the computed coefficient via the globally convergent method. 

Let $E_{glob}\left( \mathbf{x},t\right) $ be the solution
of the forward problem  (\ref{model3_1})--(\ref{model3_6}) with $%
\varepsilon _{r}\left( \mathbf{x}\right) :=\varepsilon _{r,glob}\left( 
\mathbf{x}\right).$ Denote by $p\left( \mathbf{x},t\right) =\partial
_{n} E_{glob}\left( \mathbf{x},t\right) \mid _{S_{T}}.$ In addition to the
Dirichlet condition (\ref{12}), we set the Neumann boundary condition as 
\begin{equation}
\partial _{n} E\left( \mathbf{x},t\right) =p\left( \mathbf{x},t\right)
,\forall \left( \mathbf{x},t\right) \in S_{T}.  \label{120}
\end{equation}

To formulate the Fr\'{e}chet derivative of
the Tikhonov functional (\ref{functional}) (see formula (\ref{22}) below), we make use of the adjoint method. This method is based on the state and adjoint problems. The state problem in the domain $G_{b}$ is given by 
\begin{eqnarray}
\varepsilon_r \frac{\partial^2 E}{\partial t^2} + \nabla (\nabla \cdot E) 
- \nabla \cdot (\nabla E) - s\nabla  ( \nabla \cdot(\varepsilon_r E))  &=& 0,~ \mbox{in}~~ Q_T,    \label{model4_1} \\
  E(\mathbf{x},0) = 0, ~~~E_t(\mathbf{x},0) &=&0 ~ \mbox{in}~~ G_b,      \label{model4_2}  \\
\partial _{n} E\left( \textbf{x},t\right) &=& p\left( \mathbf{x},t\right) \mbox{on}\text{ } S_T. \label{model4_3}  
\end{eqnarray}
The adjoint problem is: 
\begin{eqnarray}
\varepsilon_r \frac{\partial^2 \lambda}{\partial t^2} + \nabla (\nabla \cdot \lambda) 
- \nabla \cdot (\nabla \lambda) - s \varepsilon_r \nabla  ( \nabla \cdot  \lambda)  &=&0,~%
\mbox{in}~~Q_{T},  \label{16} \\
\lambda (\mathbf{x},T)=0,~~~\lambda _{t}(\mathbf{x},T) &=&0~\mbox{in}~~G_{b},
\label{17} \\
\partial _{n}\lambda (\mathbf{x},t) &=&z_{\delta }\left( t\right) \left( 
\widetilde{g}- E\right) \left( \mathbf{x},t\right) \mbox{on}~~S_{T}.
\label{18}
\end{eqnarray}
Here, $z_{\delta}(t)$ is used to ensure the compatibility conditions at
 $\overline{Q}_{T}\cap \left\{ t=T\right\} $ for the adjoint problem and $
\delta >0$ is a small number. The function $z_{\delta }(t)$ is chosen such
that 
\begin{equation*}
z_{\delta }\in C^{\infty }\left[ 0,T\right] ,\text{ }z_{\delta
}\left( t\right) =\left\{ 
\begin{array}{ll}
1 & \text{ for} t \in \left[ 0,T-\delta \right],  \\ 
0\ & \text{ for }t\in \left( T-\frac{\delta }{2},T\right],  \\ 
 0 < z_{\delta } < 1 & \text{ for }t\in \left( T-\delta ,T-\frac{\delta }{2}\right),
\end{array}
\right.
\end{equation*}

Weak solutions $E,\lambda \in H^{1}\left( Q_{T}\right) $ to problems
(\ref{model4_1})--(\ref{model4_3}) and (\ref{16})--(\ref{18}) 
are
defined similarly with the case of only one hyperbolic equation in Chapter 4
of the book \cite{Lad}, also see formula (34) in \cite{BM}. The weak solution to the state problem (\ref{model4_1})-(\ref{model4_3}) is the solution to the following equation:
\begin{equation}\label{weakforw}
\begin{split}
& \int_{Q_T} (-\varepsilon_r \frac{\partial E}{\partial t}\frac{\partial v}{\partial t} ~ d\bx dt -  \int_{Q_T}( \nabla \cdot E)( \nabla \cdot v)~d\bx dt  
+  \int_{Q_T}( \nabla  E)( \nabla  v)~d\bx dt  \\
&+ s  \int_{Q_T}(\nabla \cdot(\varepsilon_r E))( \nabla \cdot v)~d\bx dt - \int_{\partial S_T} v p ~d\sigma dt = 0, \forall v\in
H^{1}\left( Q_{T}\right) ,v\left( x,T\right) =0,
\end{split}
\end{equation}
The weak solution to the adjoint problem (\ref{16})-(\ref{18}) is the solution to the following equation:
\begin{equation} \label{control}
\begin{split}
&- \int_{S_T}(\tilde{g} - E)~ w~ z_{\delta}~ d \sigma dt 
-  \int_{Q_T} \varepsilon_r \frac{\partial \lambda}{\partial t} 
\frac{\partial w}{\partial t}~ d\bx dt -  \int_{Q_T}( \nabla \cdot \lambda)( \nabla \cdot w)~d\bx dt  \\
&+  \int_{Q_T}( \nabla  \lambda)( \nabla w)~d\bx dt  
+ s  \int_{Q_T} (\nabla \cdot \lambda)( \nabla \cdot (\varepsilon_r w))~ d\bx dt,\forall w\in H^{1}\left( Q_{T}\right), w\left( x,0\right) =0.
\end{split}
\end{equation}

Introduce the following spaces of real valued vector functions
\begin{equation*}
H_{E}^{1}(Q_{T})=\left\{ f\in [H^{1}(Q_{T})]^3:f(\mathbf{x},0)=0\right\} ,
\end{equation*}
\begin{equation*}
H_{\lambda }^{1}(Q_{T})=\left\{ f\in [H^{1}(Q_{T})]^3:f(\mathbf{x},T)=0\right\} ,
\end{equation*}
\begin{equation*}
U^{1}=H_{E}^{1}\left( G_{T}\right) \times H_{\lambda }^{1}\left(
G_{T}\right) \times B\left( G\right) ,
\end{equation*}%
where $B\left( G\right) $ is the space of functions bounded on $G$
with the norm $\left\Vert f\right\Vert _{B\left( G\right)
}=\sup_{G}\left\vert f\right\vert .$ 
To minimize the functional  (\ref{functional}) we introduce the Lagrangian
\begin{equation}\label{lagrangian}
\begin{split}
L(E,\lambda, \varepsilon_r) &= F(E, \varepsilon_r) 
-  \int_{Q_T} \varepsilon_r  \frac{\partial
 \lambda }{\partial t} \frac{\partial E}{\partial t}  ~d\bx dt  -  \int_{Q_T}( \nabla \cdot E)( \nabla \cdot \lambda)~d\bx dt  \\
& +  \int_{Q_T}( \nabla  E)( \nabla  \lambda)~d\bx dt  + s  \int_{Q_T} ( \nabla \cdot (\varepsilon_r E))( \nabla \cdot \lambda)~d\bx dt 
  - \int_{\partial S_T} \lambda p ~d\sigma dt, 
\end{split}
\end{equation}
where $E$ and $\lambda$ are weak solutions of problems
(\ref{model4_1})-(\ref{model4_3}) and (\ref{16})-(\ref{18}), respectively.

Clearly, (\ref{weakforw}) implies that the sum of integral terms in
 (\ref{lagrangian}) equals zero.\ Hence, $L\left(E,\lambda
 ,\varepsilon _{r}\right) =F\left(E,\varepsilon _{r}\right).$ In
 (\ref{lagrangian}) $\left(E,\lambda ,\varepsilon _{r}\right) =w\in
 U^{1}$ and functions $E$ and $\lambda $ depend on the coefficient
 $\varepsilon _{r}.$ To get the Fr\'{e}chet derivative
 $L^{\prime }$ of the Lagrangian (\ref{lagrangian}) rigorously, one
 should assume that variations of functions $u$ and $\lambda $ depend
 on variations of the coefficient $\varepsilon _{r}$, similarly with
 section 4.8 of \cite{BK}, and we will do that in our future
 publications. In this work, to derive the Fr\'{e}chet derivative
 of the Lagrangian (\ref{lagrangian}) we assume for brevity that in
 (\ref{lagrangian}) the vector function $(E,\lambda ,\varepsilon_{r})$ can be varied
 independently of each other. 
Thus, we search a point  $w \in U^1$  such that
\begin{equation}
L^{\prime }(w)\left( \overline{w}\right) =0,~~~\forall \overline{w}\in U^{1}.
\label{scalar_lagr}
\end{equation}
To find the Fr\'{e}chet derivative $L^{\prime }(w),$ we consider $L\left( w+\overline{w}\right)
-L\left( w\right) ,~\forall \overline{w}\in U^{1}$ and single out the
linear, with respect to $\overline{w},$ part of the obtained
  expression. Thus, using (\ref{weakforw}), (\ref{control}) and
  (\ref{functional}), we obtain
\begin{equation}
L^{\prime }(w)\left( \mathbf{x}\right) =\gamma \left( \varepsilon
_{r}-\varepsilon _{r,glob}\right) \left( \mathbf{x}\right)
 -\int_0^T \frac{\partial \lambda }{\partial t} 
\frac{\partial E}{\partial t}(\mathbf{x},t)~dt 
+  s\int_0^T (\nabla \cdot  E) (\nabla \cdot \lambda)(\mathbf{x},t)  ~dt,
\mathbf{x}\in G_{b}.  \label{22}
\end{equation}

\section{Finite element discretization}
\label{sec:fem}

Consider a partition $K_{h}=\{K\}$ of $G_{b}$ which consists of tetrahedra with a mesh function $h$ defined as $h|_{K}=h_{K}$ --- the local diameter of the
element $K$. Let $J_{\tau }=\left\{ J\right\} $ be a partition of the
time interval $(0,T)$ into subintervals $J=(t_{k-1},t_{k}]$ of uniform
length $\tau =t_{k}-t_{k-1}$. We also assume  the minimal angle condition on
the $K_{h}$ \cite{Brenner}.


To formulate the finite element method for solving the state problem (\ref{model4_1})--(\ref{model4_3}) and the adjoint problem (\ref{16})--(\ref{18}), and to compute the gradient of the Lagrangian via (\ref{22}),  we define the finite element
spaces $V_{h}\subset L_{2}\left( G_{b}\right) $, $W_{h}^{E}\subset
H_{E}^{1}\left( Q_{T}\right) $ and $W_{h}^{\lambda }\subset H_{\lambda
}^{1}\left( Q_{T}\right) $. First, we introduce the finite element trial
space $W_{h}^{E}$ for every component of the electric field $E$ defined by 
\begin{equation}
W_{h}^{E}:=\{w\in H_{E}^{1}:w|_{K\times J}\in P_{1}(K)\times
P_{1}(J),\forall K\in K_{h},\forall J\in J_{\tau }\},  \notag
\end{equation}%
where $P_{1}(K)$ and $P_{1}(J)$ denote the set of linear functions on $K$
and $J$, respectively. We also introduce the finite element test space $%
W_{h}^{\lambda }$ defined by 
\begin{equation}
W_{h}^{\lambda }:=\{w\in H_{\lambda }^{1}:w|_{K\times J}\in P_{1}(K)\times
P_{1}(J),\forall K\in K_{h},\forall J\in J_{\tau }\}.  \notag
\end{equation}%
Hence, the finite element spaces $W_{h}^{E}$ and $W_{h}^{\lambda }$ consist
of continuous piecewise linear functions in space and time. 
To approximate the function $\varepsilon _{r}(\mathbf{x}),$ we use the space
of piecewise constant functions $V_{h}\subset L_{2}\left( \Omega \right) $, 
\begin{equation}
V_{h}:=\{u\in L_{2}(\Omega ):u|_{K}\in P_{0}(K),\forall K\in K_{h}\},
\label{p0}
\end{equation}%
where $P_{0}(K)$ is the set of piecewise constant functions on $K$. In our
computations we truncate computed functions $\varepsilon _{r}\left( \mathbf{x%
}\right) $ to unity outside of the domain $\Omega _{FEM}=\Omega $ using (\ref%
{coefic}) as%
\begin{equation}
\overline{\varepsilon }_{r}\left( \mathbf{x}\right) =\left\{ 
\begin{array}{c}
\varepsilon _{r}\left( \mathbf{x}\right) ,\mathbf{x}\in \Omega _{FEM}, \\ 
1,\mathbf{x}\in \Omega _{FDM}.%
\end{array}%
\right.   \label{23}
\end{equation}


Next, we set $U_{h}=W_{h}^{E}\times W_{h}^{\lambda }\times V_{h}$.
Obviously $\dim U_{h}<\infty $ and $U_{h}\subset U^{1}$ as a
set. Because of this, we consider $U_{h}$ as a discrete analogue of
the space $U^{1}.$ We introduce the same norm in $U_{h}$ as the one in
$U^{0},\left\Vert \cdot \right\Vert _{U_{h}}:=\left\Vert \cdot
\right\Vert _{U^{0}}$, with 
\begin{equation*}
U^{0}=L_{2}\left(G_{T}\right) \times L_{2}\left(G_{T}\right) \times
L_{2}\left( \Omega \right). 
\end{equation*}
The finite element method for solving equation
(\ref{scalar_lagr}) now reads: \emph{Find }$u_{h}\in U_{h}$\emph{,
  such that}
\begin{equation}
L^{\prime }(u_{h})(\bar{u})=0, ~\forall \bar{u}\in U_{h}.  \label{varlagr}
\end{equation}

\section{Mesh refinement recommendation and the  adaptive algorithm}
\label{subsec:ad_alg}

From Theorem 5.1 and Remark 5.1  of \cite{BMaxwell2} it follows that the  finite element mesh  should be
locally refined  in such subdomain of $\Omega$ where the maximum norm of the Fr\"echet derivative of the objective functional is large. 
For each mesh we first linearly interpolate the coefficient
$\varepsilon_{r,glob}(\mathbf{x}) $ on it, and use the interpolated coefficient as an initial
guess for on the current mesh. Our algorithm consists of two loops: the outer loop deals with the locally adaptive mesh refinement. In the inner loop, i.e., on each mesh, we iteratively update
the approximations $\varepsilon_{h}^{m}$ of the function $\varepsilon_{
  h}$ by solving (\ref{scalar_lagr}) using an optimization procedure, where $m$ is the index of iteration in the optimization procedure.
Denote by
\begin{equation}\label{Bhm}
\begin{split}
L_h^{\prime,m}(\mathbf{x}) = - {\int_0}^T  \frac{\partial \lambda_h^m}{\partial t}
\frac{\partial E_h^m}{\partial t}~ dt  +  s\int_{0}^T \nabla \cdot E_h^m \nabla \cdot \lambda_h^m  dt + \gamma (\bar{\varepsilon_h}^m - \bar{\varepsilon}_{r,glob}).
\end{split}
\end{equation}

\textbf{Adaptive algorithm}

\begin{itemize}
\item[Step 0.]  Choose an initial mesh $K_{h}$ in $\Omega$ and an
initial time partition $J_{0}$ of the time interval $\left( 0,T\right) .$
Start from the initial guess $\varepsilon_{h}^{0}= \varepsilon_{r,glob}$, we compute the approximations $\varepsilon_{h}^{m}$ via the following steps:

\item[Step 1.]  Compute the solutions $E_{h}\left(
  x,t,\varepsilon_{h}^{m}\right) $ and $\lambda _{h}\left(
  x,t,\varepsilon_{h}^{m}\right) $ of the state problem  
(\ref{model3_1})--(\ref{model3_4}) and the adjoint problem (\ref{16})--(\ref{18}) on $K_{h}$ and $J_{k}$, and compute the Fr\"echet derivative $L^{\prime,m}_h$ via (\ref{Bhm}).

\item[Step 2.]  Update the coefficient  
 on $K_{h}$ and $J_k$ using the conjugate gradient method:
\begin{equation*}
\varepsilon_h^{m+1} : =  \varepsilon_h^{m}  + \alpha d^m(\mathbf{x}),
\end{equation*}
where $\alpha > 0$ is a step-size in the conjugate gradient method, given by a line search procedure, see, e.g., \cite{Peron}, and
\begin{equation*}
\begin{split}
 d^m(\mathbf{x})&=  -L^{\prime,m}_h(\mathbf{x})  + \beta^m  d^{m-1}(\mathbf{x}),
\end{split}
\end{equation*}
with
\begin{equation*}
\begin{split}
 \beta^m &= \frac{|| L^{\prime,m}_h ||^2}{|| L^{\prime,m-1}_h ||^2},
\end{split}
\end{equation*}
where $d^0(\mathbf{x}) = -L^{\prime,0}_h(\mathbf{x})$.
\item[Step 3.]  Stop updating the coefficient and set
  $\varepsilon_h : = \varepsilon_h^{m+1}$, $M:= m+1$, if either $||L_h^{\prime,m}||_{L_{2}( \Omega)}\leq
  \theta$ or norms 
$||\varepsilon_h^{m} ||_{L_{2}(\Omega)}$ 
are stabilized. Here
  $\theta$ is a tolerance number.
  Otherwise, set $m:=m+1$ and go to step~1.

\item[Step 4.] Compute $L^{\prime,M}_h$ via (\ref{Bhm}). Refine the mesh at all grid points $\bx$ where 
\begin{equation}
|L^{\prime,M}_h\left( \bx\right) |  \geq \beta _{1}\max_{\overline{\Omega }} |L_{h}^{\prime,M}\left(
\bx\right)|.  \label{62}
\end{equation}
Here the tolerance number $\beta _{1}\in \left( 0,1\right) $ is chosen by
the user.

\item[Step 5.]  Construct a new mesh $K_{h}$ in $\Omega$ and a new time
partition $J_{k}$ of the time interval $\left( 0,T\right) $. On $J_{k}$ the
new time step $\tau $ should be chosen in such a way that the CFL condition
is satisfied. 
Interpolate the initial approximation $\varepsilon _{r,glob}$ from the
previous mesh to the new mesh. Next, return to step 1 at $m=1$ and perform
all above steps on the new mesh.

\item[Step 6.] Let $\left\Vert L_{h,prev}^{\prime,M}\right\Vert _{L_{2}\left(
  \Omega \right) }$ and $\left\Vert L_{h,current}^{\prime,M}\right\Vert
  _{L_{2}\left( \Omega \right) }$ be norms defined in step 4 on the
  previous and current mesh, respectively. Stop mesh refinements if
  $\left\Vert L_{h,current}^{\prime,M}\right\Vert _{L_{2}\left( \Omega
    \right) }\geq \left\Vert L_{h,prev}^{\prime,M}\right\Vert _{L_{2}\left(
    \Omega \right) }.$

\end{itemize}


\section{Some details of numerical implementation}

\label{sec:numex}

In this section we present results
of reconstruction of dielectric constants/refractive indices and shapes of
some targets using the adaptivity algorithm of section \ref{subsec:ad_alg}.
 One of the discrepancies between our mathematical model
 (\ref{model3_1})- (\ref{model3_6}) and the measured experimental data
 is that formally equation (\ref{model3_1}) is invalid for the case
 when metallic targets are present (we refer to \cite{BTKF} for
 the description of other discrepancies).  However, it was shown
 computationally in \cite{KBKSNF} that one can treat metallic targets
 as dielectrics with large dielectric constants.

 We
call these  \emph{effective} (or ``appearing'') dielectric constants and values for them are in the interval
\begin{equation}
\varepsilon _{r}\left( \text{metallic target}\right) \in \left( 10,30\right)
.  \label{2.51}
\end{equation}
Modeling metallic targets as integral parts of the unknown coefficient
$\varepsilon _{r}\left( \mathbf{x}\right) $ is convenient for our
practical computations in order to image IEDs.
  Since IEDs usually consist of mixtures of some
dielectrics with a number of metallic parts, these targets are
heterogeneous ones, and we consider three heterogeneous cases in
section \ref{sec:8.2}.  However, modeling metallic parts separately from
dielectric ones is impractical for our application because of those
mixtures.

 Using (\ref{2.51}), we define in all our tests the upper value
 of the function $\varepsilon_{r}\left(\mathbf{x}\right)$ as $d=25,$
 see (\ref{2.20}). Thus, we set lower and upper bounds for the
 reconstructed function $\varepsilon_{r}(\mathbf{x})$ in $\Omega$ as
\begin{equation}
M_{\varepsilon _{r}}=\{\varepsilon _{r}(\mathbf{x}):\varepsilon _{r}\left( 
\mathbf{x}\right) \in \left[ 1, 25\right] \}.  \label{4.30}
\end{equation}
In our computation we ensure lower and upper bounds via truncating those
values of $\varepsilon _{r}\left( \mathbf{x}\right) $ which are outside of
the interval (\ref{4.30}).

\subsection{Dirichlet and Neumann boundary conditions}

In our experiments only one component $E_{2}\left( \mathbf{x},t\right) $ of
the electric field $E\left( \mathbf{x},t\right) =\left(
E_{1},E_{2},E_{3}\right) \left( \mathbf{x},t\right) $ is both sent into the
medium and measured. Thus, only the second component of the function $g$ in the Dirichlet boundary condition (\ref{2.5}) on $\Gamma$ is available. We
approximate the other two components of $g$ on $\Gamma$ by the numerical solution of
the forward problem (\ref{model3_1})--(\ref{model3_6}), with the coefficient given by $\varepsilon_r = \varepsilon _{r,glob}(\mathbf{x})$ --- the solution of the globally convergent method.  The Dirichlet  data on the rest of the boundary $\partial G_{b},$ i.e., on $\partial G_b\setminus \Gamma$, as well as the Neumann condition $%
p\left( \mathbf{x},t\right) $ in (\ref{120}) at the entire boundary $\partial G_{b}$ are taken from the numerical solution of that forward problem. 

\subsection{Computational domains}

To generate the boundary data (\ref{12}), (\ref{120}) for all three
components of the electric field $E,$ as specified in the previous section,
we solve the forward problem in the computational domain $G,$ which we
choose as 
\begin{equation*}
G=\left\{ \mathbf{x=}(x,y,z)\in (-0.56,0.56)\times (-0.56,0.56)\times
(-0.16,0.1)\right\}.
\end{equation*}
The boundary of the domain $G$ is $\partial
G=\partial _{1}G\cup \partial _{2}G\cup \partial _{3}G.$ Here, $\partial
_{1}G$ and $\partial _{2}G$ are front and back sides of the
domain $G$ at $\{z=0.1\}$ and $\{z=-0.16\}$, respectively, and $\partial _{3}G$ is the union of left, right, top and bottom
sides of this domain.

We use a stabilized domain decomposition method of \cite{BM}
implemented in the software package WavES \cite{waves}. The FEM domain $\Omega_{FEM}$ is chosen as
\begin{equation}
\Omega_{FEM}=\Omega =\left\{ \mathbf{x=}(x,y,z)\in (-0.5,0.5)\times
(-0.5,0.5)\times (-0.1,0.04)\right\} .  \label{8.0}
\end{equation}
After the data propagation procedure, see discussions in section \ref{sec:2}, the
data $g\left( \mathbf{x},t\right) $ in (\ref{2.5}) are given at the front
side $\Gamma $ of the domain $\Omega $ which is defined as 
\begin{equation}
\Gamma =\{\mathbf{x}\in \partial \Omega :z=0.04\}.  \label{8.1}
\end{equation}

The waveform function  $f(t)$ in our simulated incident plane wave is chosen as 
\begin{equation*}
f(t)=\sin \omega t,~0\leq t\leq t_{1}:=\frac{2\pi }{\omega }.
\end{equation*}
Here, we initialize the plane wave at $\left\{ z=0.1\right\} .$
We use $\omega =30$, $T=1.2$  and $s=1$.  We solve the problem (\ref{model3_1})--(\ref{model3_6}) using the
explicit scheme of \cite{BM}  with the time step size $\tau =0.003,$ which satisfies the
CFL condition. Note that this time step is dimensionless, which corresponds to the time step of 10 ps in our real experiment. Here we use the dimensionless time step in order to normalize the coefficient $\varepsilon_r$ to be unity outside of $\Omega$. The dimensionless time $T = 1.2$ corresponds to $4$ ns in our real time. We do not use the whole 10 ns recorded data since after data preprocessing, they are shifted earlier in time and they do not contain any target's signal after 4 ns. 

\subsection{Data Preprocessing}

\label{sec:2}

In the previous section, we implicitly assumed that our model is comparable to the experimental data. Moreover, the data is available at the plane $\Gamma$ of $\Omega$, which is quite close to the targets. Unfortunately, this is not the case in practice. In fact, there is a huge misfit bewteen our experimental data and the simulated ones, see \cite{NBKF}. Therefore, data preprocessing is required in order to prepares the experimental data to become an input for our inversion algorithm. In our experience, data preprocessing of experimental data is always a
heuristic procedure. That procedure for the globally convergent method was
described in detail in \cite{NBKF}. Since the globally convergent method works with the PDE
which is obtained by the Laplace transform of the original wave-like
equation (\ref{1}), and we work here directly in time domain, we need less number of data preprocessing steps than in \cite{NBKF}. We only describe the following preprocessing steps which are different from those of \cite{NBKF}. 

\subsubsection{Data propagation} 

\begin{figure}[tph]
 \centering
  \subfloat[]{\includegraphics[width=0.45\textwidth,height=0.4\textwidth]{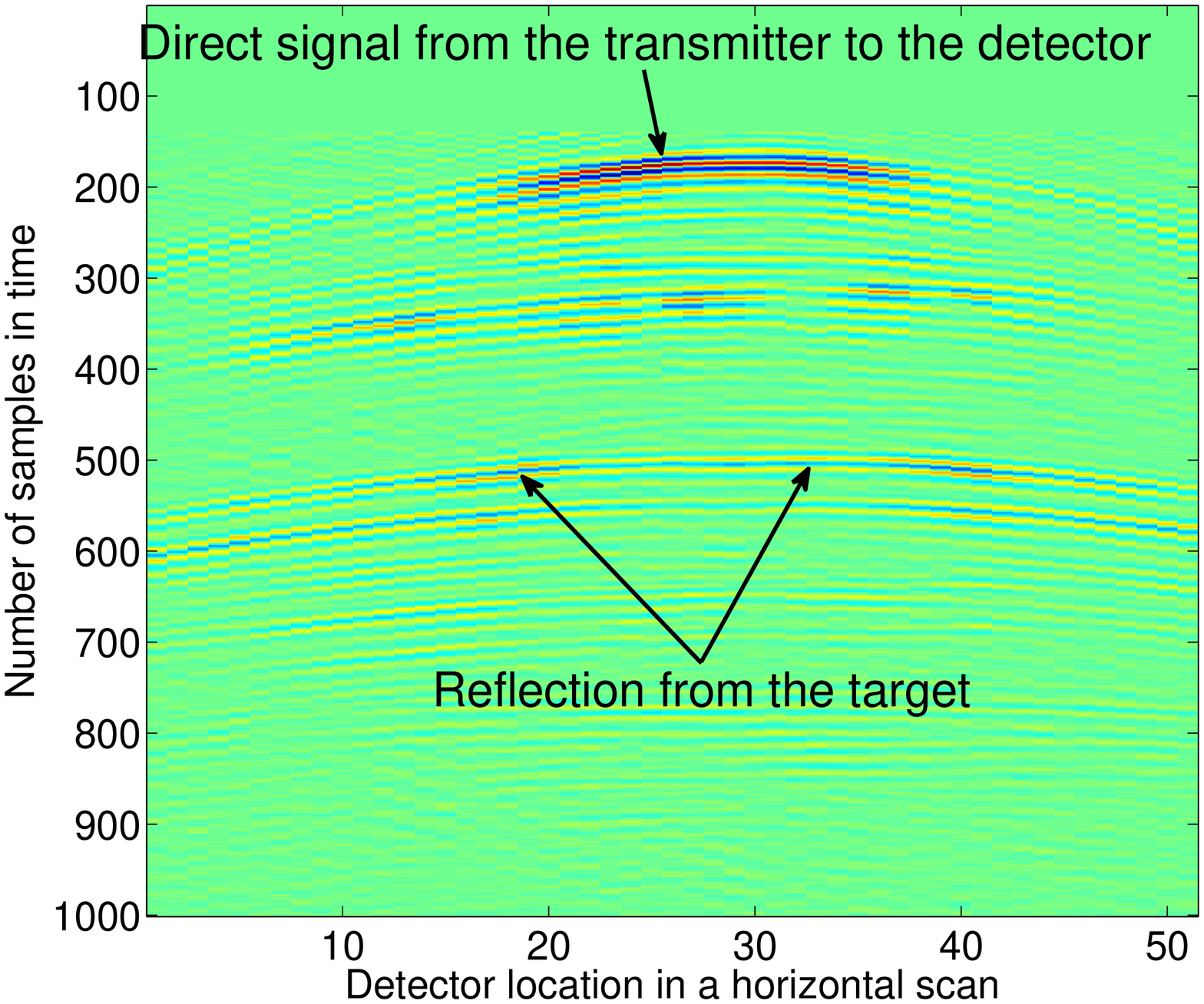}}
  \subfloat[]{\includegraphics[width=0.45\textwidth,height=0.4\textwidth]{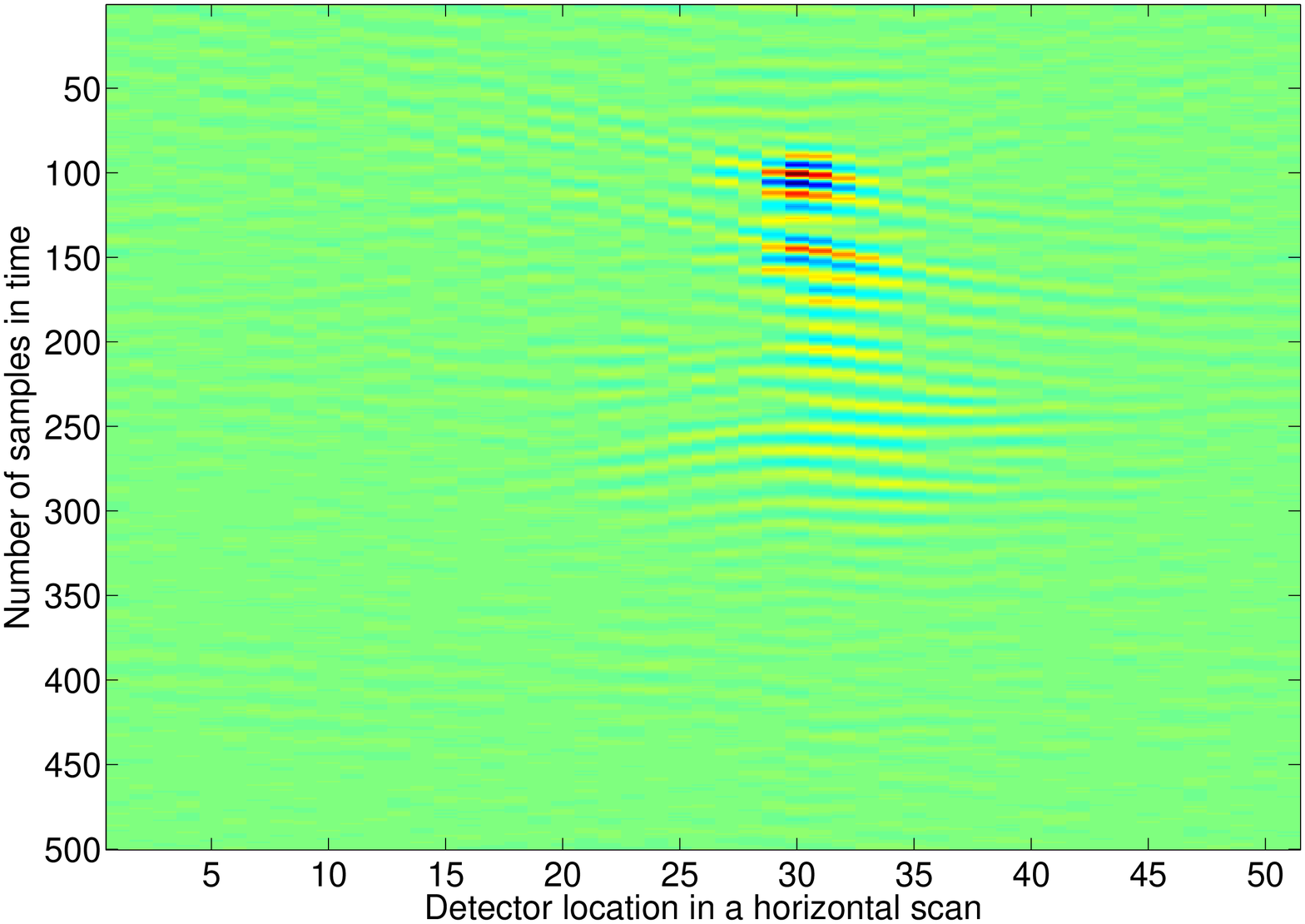}}
  \caption{Result of data propagation for target number 1 (Table \ref{tab:table1}).
Horizontal axis is the indices of the detector's locations and vertical axis is number of samples in time. a) data at the measurement plane, b) data at the propagated
plane. One can see that the propagated data are focused at the target,
whereas original data are smeared out. }
\label{fig:dataprop}
\end{figure}

In \cite{NBKF}, we used a time-reversal data propagation method in order to migrate our data from the measurement plane, which is in the far field zone, to the plane $\Gamma$, which is at about 4 cm far from the targets. In this paper, we use another data propagation method based on the Fourier transform. This technique is basically the same as the Stolt migration  in Geophysics, see \cite{Stolt:1978,Yilmaz:1987}. However, in the
standard Stolt migration the wave at the initial time is calculated in the
whole spatial domain of interest, whereas we calculate the wave only at a
plane parallel to the measurement plane but in the whole time interval. The
technique is described in detail in \cite{TBKF:2013-2}. 

A result of the data propagation is illustrated in Figure \ref{fig:dataprop}%
. The figure shows a horizontal scan of target number 1, see Table \ref{tab:table1}. The horizontal axis is the indices of the detector's
locations and the vertical axis is the number of samples in time. Time increases from the top to
the bottom. 
Figure \ref{fig:dataprop}(a) shows the original data while Figure \ref%
{fig:dataprop}(b) shows the data after the propagation. As can be seen from
these figures, the target's signas in the original data is smeared out. On
the other hand, it is focused after the data propagation.

%
%

\subsubsection{Data calibration} 

Since the amplitude of experimental data are very different from that of computational simulations, we must multiply the experimental data
by a calibration factor $\theta$ so that they have similar amplitude as the simulated data. This factor is not easy to obtain in practice and we should somehow find an appropriate one.
The calibration of the measured data was done in the Laplace transform domain in \cite{NBKF}. However, here we work with the time domain data only. Therefore, we use a new data calibation procedure we described below.

Let the function $g_{\text{exper}}\left( \mathbf{x},t\right)
,\mathbf{x}\in \Gamma ,t\in \left( 0,T\right) $ be our propagated
experimental data. This function is given only on grid points $\left(
\mathbf{x}_{i},t_{j}\right).$ We compute the maximal value of this
function, 
\begin{equation}
g_{max}=\max_{\left(  \mathbf{x}_{i},t_{j}\right) }g_{\text{exper}}\left( \mathbf{x}%
_{i},t_{j}\right).
\end{equation} 
Usually the number $g_{max}$ is quite
large. Next, let $E_{2,sim}\left( \mathbf{x},t\right) $ be the
function which is computed
via solving the problem (\ref{model3_1})-(\ref{model3_6}).  We compute
the maximal value of this function on $\Gamma ,$
\begin{equation*}
E_{2, max} =\max_{\Gamma \times \left[ 0,T\right] }E_{2,sim}\left(
\mathbf{x},t\right).
\end{equation*}
 Define $r= E_{2, max}/g_{max}.$
Next, we assign 
\begin{equation*}
g_{incl}\left( \mathbf{x}_{i},t_{j}\right) :=r\cdot
g_{\text{exper}}\left( \mathbf{x}_{i},t_{j}\right)
\end{equation*}
 and use the
  function $g_{incl}\left( \mathbf{x}_{i},t_{j}\right) $ as the second
  component of the vector function $g\left( \mathbf{x},t\right) $ in
  (\ref{2.5}). Two other steps of data preprocessing are due to (\ref{12}) and (\ref{120}) in section \ref{sec:invprobl}. The final step of data preprocessing, the so-called ``immersing
procedure'' is done as follows.

\subsubsection{Data immersing}

\label{subsec:immersing}

In this section we describe a heuristic immersing procedure of the
time-dependent propagated experimental data $g_{incl}\left( \mathbf{x}%
,t\right) =E_{2}\left( \mathbf{x},t\right) \mid _{\mathbf{x}\in \Gamma}.$
This procedure does two things:
\begin{itemize}

\item   immerses the data $g_{incl}\left( 
\mathbf{x},t\right) $ into computationally simulated ones;

\item  extends the data $g_{incl}\left( \mathbf{x},t\right) $ from $\Gamma $
to $\Gamma _{1}.$
\end{itemize}

By (\ref{8.1}) the rectangle $\Gamma $ is smaller than the rectangle 
\begin{equation*}
\Gamma _{1}=\left\{ \mathbf{x}:\left( x,y\right) \in \left(
-0.56,0.56\right) \times \left( -0.56,0.56\right) ,z=0.04\right\} .
\end{equation*}%
It is clear from the adjoint problem (\ref{16})-(\ref{18}) that we need to
get a proper data for the function $E_{2}\left( \mathbf{x},t\right) $ for $%
\left( \mathbf{x},t\right) \in \Gamma _{1}\times \left( 0,T\right) $ while
having the data $E_{2}\left( \mathbf{x},t\right) =g_{incl}\left( \mathbf{x}%
,t\right) $ only for $\mathbf{x}\in \Gamma .$ We now describe how do we
extend the data from $\Gamma $ to $\Gamma _{1}$. Let $E_{2}\left( \mathbf{x}%
,t\right) $ be the $E_{2}$-component of the solution $E\left( \mathbf{x}%
,t\right) $ of the forward problem (\ref{model3_1}--\ref{model3_6}) with $\varepsilon_r := \varepsilon_{r,glob}$. Then we define our immersed function 
$E_{2}^{immers}\left( \mathbf{x},t\right) $ for $\left( \mathbf{x},t\right)
\in \Gamma _{1}\times \left( 0,T\right) $ as 
\begin{equation}
E_{2}^{immers}\left( \mathbf{x},t\right) =\left\{ 
\begin{array}{ll}
g_{incl}\left( \mathbf{x},t\right) , & \text{ if }\mathbf{x}\in \Gamma \text{
and }g_{incl}\left( \mathbf{x},t\right) \geq \beta \max_{\overline{\Gamma }%
}g_{incl}\left( \mathbf{x},t\right) , \\ 
E_{2}\left( \mathbf{x},t\right) , & \text{ if }\mathbf{x}\in \Gamma \text{
and }g_{incl}\left( x,t\right) <\beta \max_{\overline{\Gamma }%
}g_{incl}\left( \mathbf{x},t\right) , \\ 
E_{2}\left( \mathbf{x},t\right) , & \text{ if }\mathbf{x}\in \Gamma
_{1}\diagdown \Gamma \text{.}%
\end{array}%
\right.   \label{5.109}
\end{equation}%
We choose the parameter $\beta \in \left( 0,1\right) $ in (\ref{5.109}) in
numerical experiments of section \ref{sec:8.2}. It follows from (\ref{5.109}%
) that $E_{2}^{immers}\left( \mathbf{x},t\right) =E_{2}\left( \mathbf{x}%
,t\right) $ for $\mathbf{x}\in \Gamma _{1}\diagdown \Gamma .$ 

Figures \ref{fig:figure1_1}, \ref{fig:figure1} show that, depending on the parameter 
$\beta $ in (\ref{5.109}), the data immersing procedure not only allows to
extend the data from $\Gamma $ to $\Gamma _{1}\diagdown \Gamma $ but also
make the experimental data usable in our inverse algorithm. Indeed, we note that the experimental data is measured at a very high frequency, say, $\omega \approx 170$, whereas our simulations are done at $\omega = 30$ in order to reduce the computational cost. Therefore, the experimental data are not compatible with the simulations. Our immersing procedure helps to avoid solving the problem at a very high frequency. After this
immersing procedure we solve the inverse problem using the algorithm of
section \ref{subsec:ad_alg}.

\subsection{Postprocessing of results}

Results of the globally convergent algorithm of the first stage procedure
have demonstrated that this algorithm provides accurate locations of targets
as well as accurate values of refractive indices $n=\sqrt{\varepsilon _{r}}$
of dielectric targets and large values of appearing dielectric constants $%
\varepsilon _{r}$ for metallic targets \cite{BTKF,NBKF}. However, it does
not reconstruct shapes of targets well, especially the size in the $z$-direction. The latter is the reason why we
apply the second stage to refine results of the first.

 Let $\varepsilon_{r}\left( \mathbf{x}\right) $ be the function
 obtained in the adaptive algorithm of section
 \ref{subsec:ad_alg}.  We form the image of the dielectric targets  based
 on the function $\varepsilon_{r,diel}\left( \mathbf{x}\right),$
\begin{equation}\label{8.6}
\varepsilon _{r,diel}\left( \mathbf{x}\right) =\left\{ 
\begin{array}{l}
\varepsilon_{r}\left( \mathbf{x}\right) \text{ if } \varepsilon_{r}\left( \mathbf{x}\right) \geq 0.85\max_{\overline{\Omega }}%
\varepsilon_{r}\left( \mathbf{x}\right) , \\ 
1\text{ otherwise.}
\end{array}
\right.
\end{equation}
As to the metallic targets (i.e., the ones with large computed maximal values
of $\varepsilon _{r}\left( \mathbf{x}\right) ),$ we use the function $%
\varepsilon _{r,metal}\left( \mathbf{x}\right) ,$ 
\begin{equation}\label{8.7}
\varepsilon _{r,metal}\left( \mathbf{x}\right) =\left\{ 
\begin{array}{l}
\varepsilon_{r}\left( \mathbf{x}\right) \text{ if } \varepsilon_{r}\left( \mathbf{x}\right) \geq 0.3 \max_{\overline{\Omega }}%
\varepsilon_{r}\left( \mathbf{x}\right) , \\ 
1\text{ otherwise.}
\end{array}
\right.
\end{equation}

\begin{table}[p]
\begin{center}
\begin{tabular}{|c|l|}
\hline
\textbf{\ Target number} & \textbf{Specification of the target}   \\ \hline
1 & a piece of oak, rectangular prism   \\ 
2 & a piece of pine   \\ 
3 & a metallic sphere   \\ 
4 & a metallic cylinder   \\ 
5 & a piece of oak   \\ 
6 & a metallic rectangular prism \\  
7 & a wooden doll, air inside, heterogeneous target \\ 
8 & a wooden doll, metal inside, heterogeneous target \\ 
9 & a wooden doll, sand inside, heterogeneous target\\ 
 \hline
\end{tabular}
\end{center}
\caption{\emph{Names of targets.}}
\label{tab:table1}
\end{table}

\begin{table}[p]
\begin{center}
\begin{tabular}{|l|l|l|l|l|l|c|}
\hline
Target number & 1 & 2 & 5 & 7 & 9 & Average error \\ \hline
Measured $n$, error & 2.11, 19\% & 1.84, 18\% & 2.14, 28\% & 
1.89, 30\% & 2.1, 26\% & 24\% \\ \hline
$n$ in glob.conv, error & 1.92, 9\% & 1.8, 2\% & 1.83, 15\% & 
1.86, 2\% & 1.92, 9\% & 8\% \\ \hline
$n$, coarse mesh, error & 1.94, 8\% &  1.82, 1\% &  1.84, 14\% & 
1.88, 0.5\% & 1.93, 8\% & 6\% \\ \hline
$n$, 1 time ref. mesh, error & 1.94, 8\% & 1.82, 1\% & 1.85, 14 \% & 
1.89, 0\% & 1.93, 8\%  &  6\% \\ \hline
$n$, 2 times ref.mesh, error & &  1.84, 0\% &  & 
1.9, 0.5\% & 1.96, 7\% & 2\% \\ \hline
$n$, 3 times ref.mesh, error &  &  &   & 
1.89,0 \% &  & 0\% \\ \hline
\end{tabular}%
\end{center}
\caption{\emph{ Computed $n$(target) and directly measured refractive
    indices of dielectric targets together with both measurement and
    computational errors as well as the average error. Note that the
    average computational errors are at least 4 times less than the
    average error of direct measurements.  In all tests we have used
    the following values of above parameters: the regularization
    parameter $\gamma =0.01$ in (\ref{functional}), $\theta =10^{-9}$ in
    Step 3 of the adaptive algorithm, $\beta _{1}=0.7$ in (\ref{62}) and
    $s=1$ in (\ref{model3_1}). } }
\label{tab:table2}
\end{table}

\begin{table}[p]
\begin{center}
\begin{tabular}{|l|l|l|l|l|l|l|l|}
\hline
Target number & 3 & 4 & 6 & 8  \\ \hline
$\varepsilon_{r}(target)$ of glob.conv. & 14.4 & 15.0 & 25 & 13.6   \\ 
\hline
$\varepsilon_{r}(target)$ coarse mesh & 14.4      &  17.0  & 25  & 13.6   \\ 
\hline
$\varepsilon_{r}(target)$  1 time ref.mesh & 14.5 &  17.0  & 25 &  13.6  \\ 
\hline
$\varepsilon_{r}(target)$  2 times ref.mesh & 14.6 & 17.0  & 25 &  13.7  \\ 
\hline
$\varepsilon_{r}(target)$  3 times ref.mesh & 14.6 & 17.0  &  &  14.0  \\ 
\hline
$\varepsilon_{r}(target)$  4 times ref.mesh &  & 17.0  &  &   \\ 
\hline
\end{tabular}
\end{center}
\caption{\emph{Computed appearing dielectric constants }$\protect\varepsilon_{r}(target)$\emph{\ of metallic targets with numbers 3,4,6 as well as
of target number 8 which is a metal covered by a dielectric. The mesh refinement process for target number 6 has stopped on three (3)
mesh refinements.}}
\label{tab:table3}
\end{table}

\section{Reconstruction results} \label{sec:8.2}

\begin{figure}[h]
  \begin{center}
    \begin{tabular}{cc}
      \includegraphics[width=0.45\textwidth]{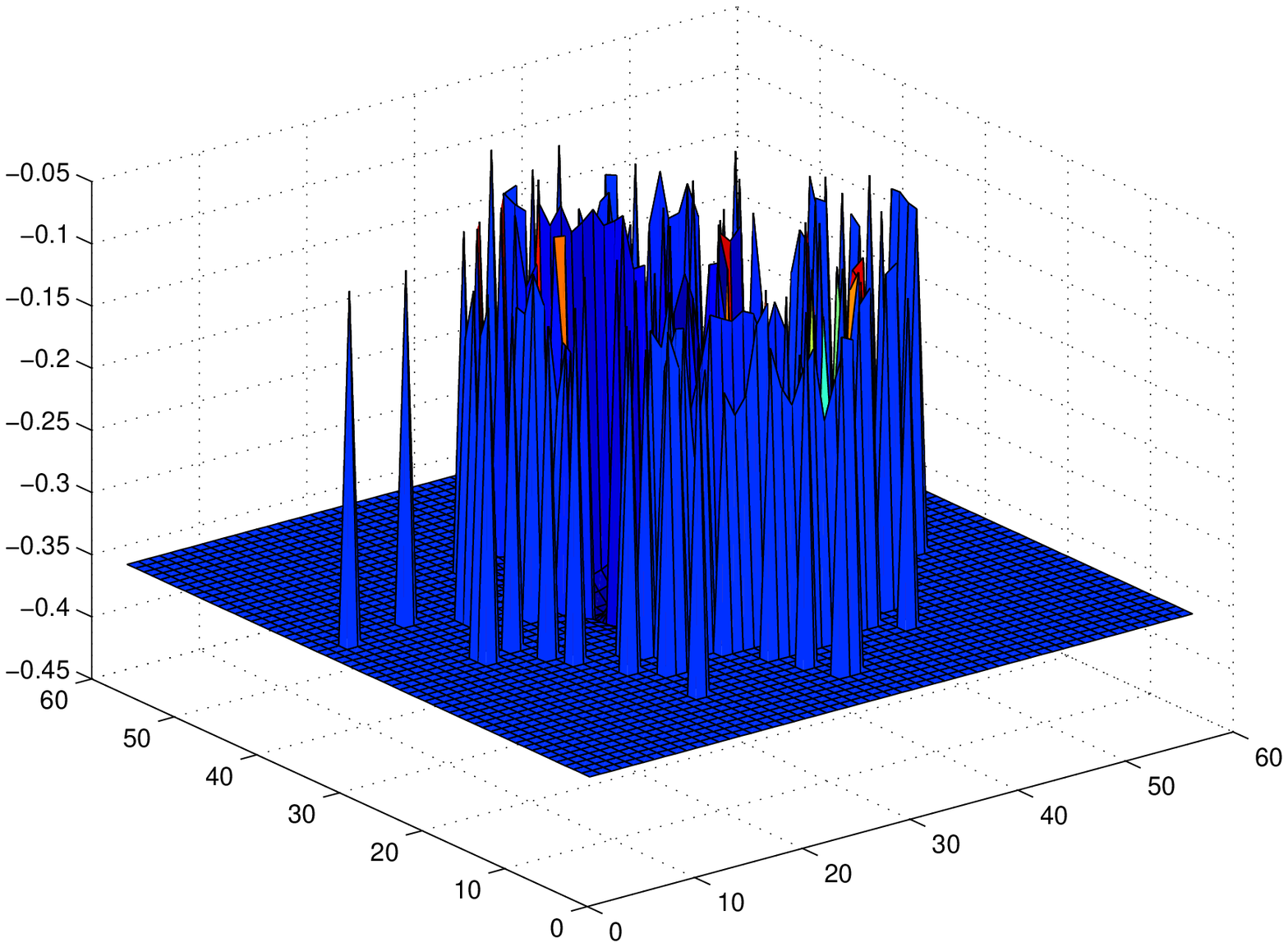}&
      \includegraphics[width=0.45\textwidth]{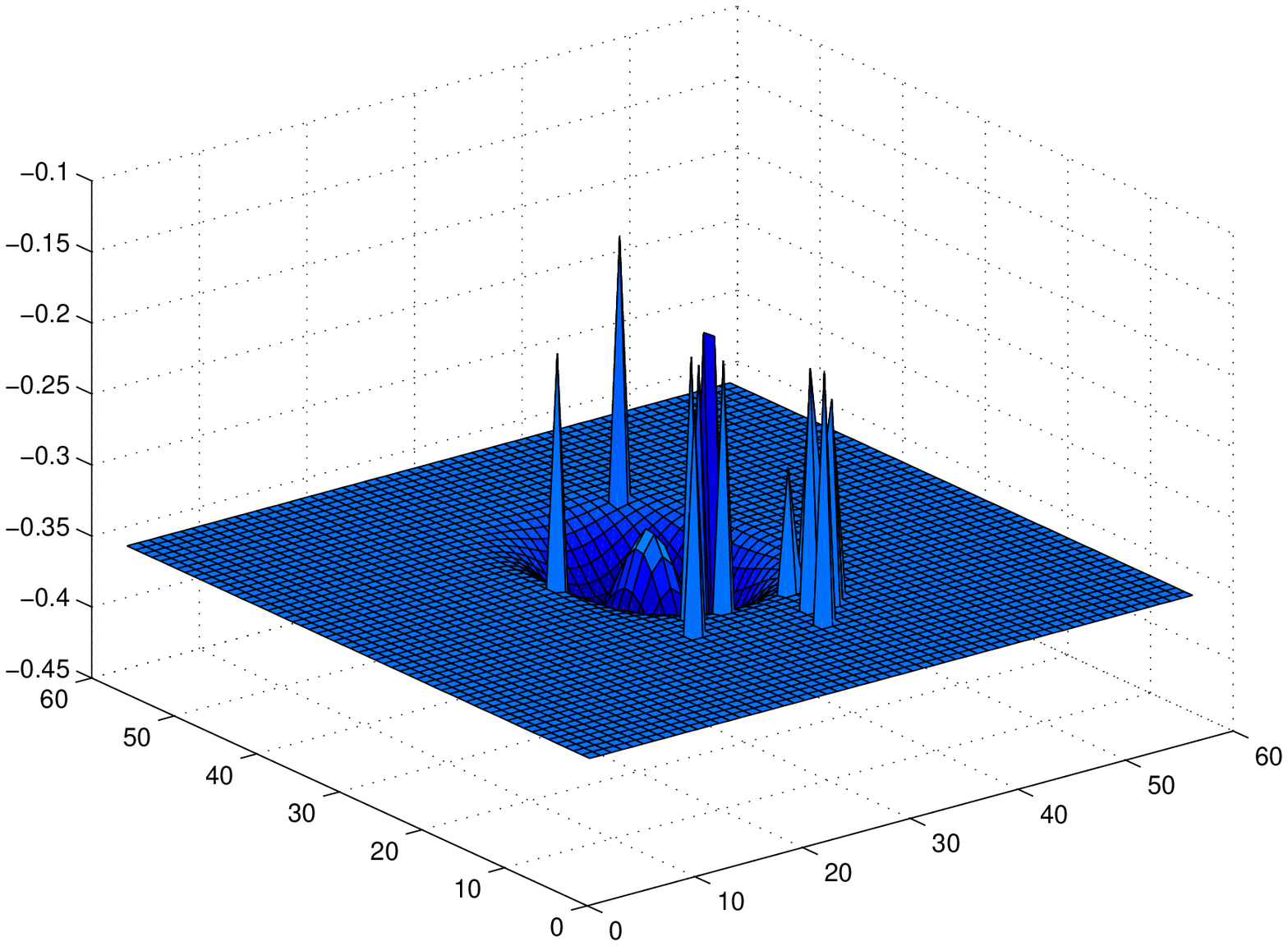}\\
      (a) $\beta = 0.2,t = 0.3 $ & (b) $\beta = 0.5, t = 0.3$\\
      \includegraphics[width=0.45\textwidth]{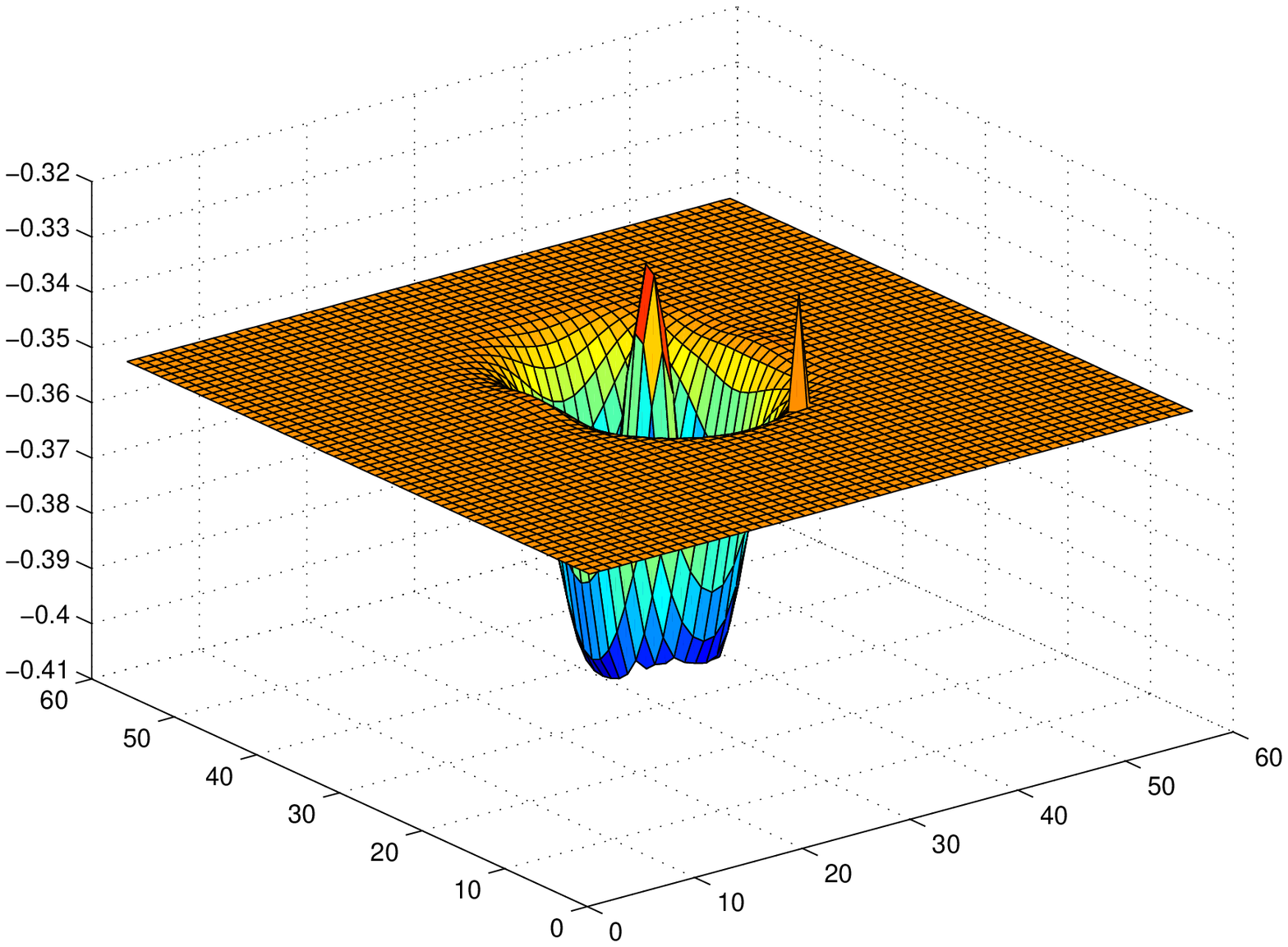}&
      \includegraphics[width=0.45\textwidth]{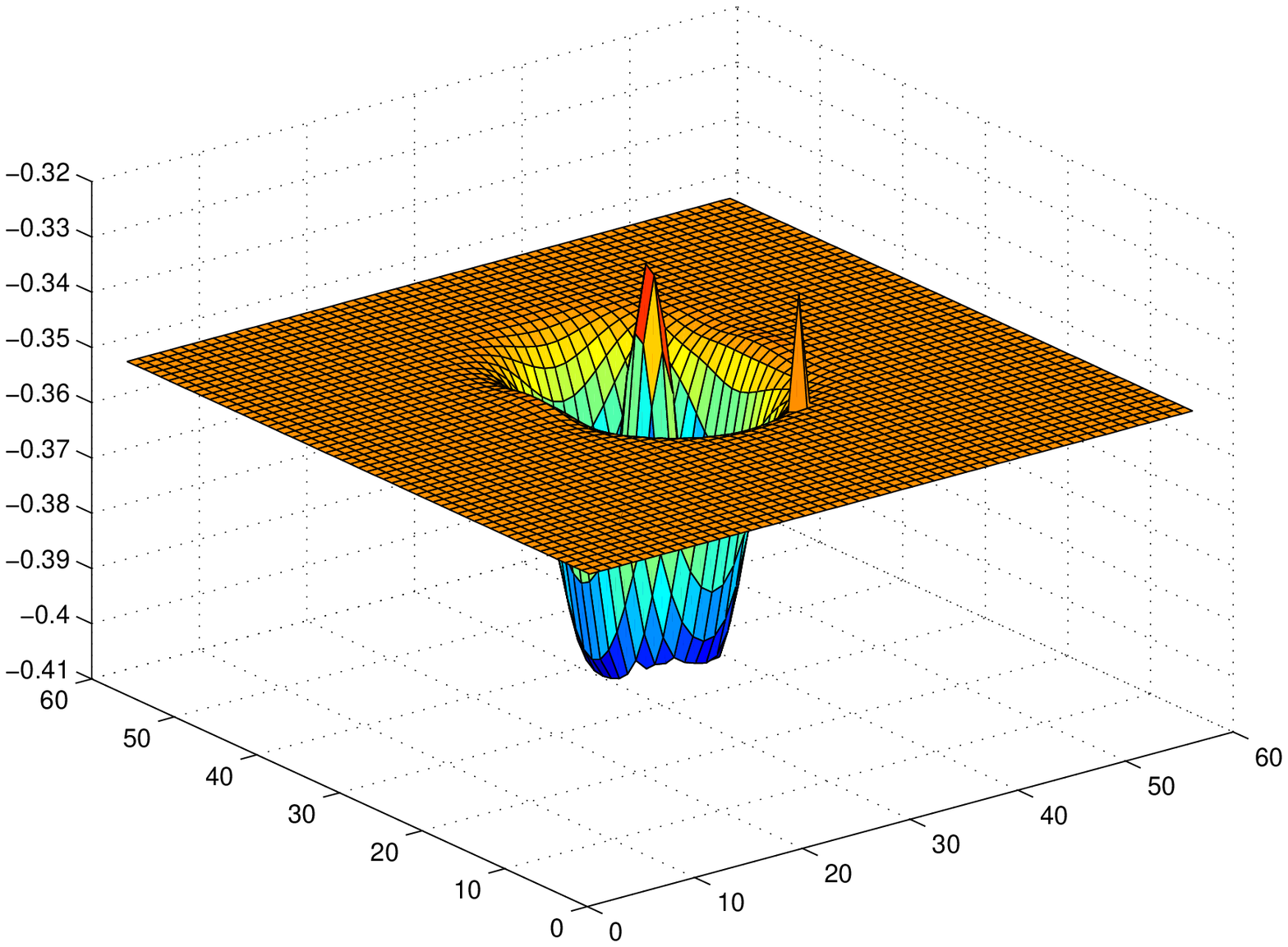}\\
      (c) $\beta = 0.8,t = 0.3 $ & (d) $\beta = 0.9, t = 0.3$
    \end{tabular}
 \end{center}
    \caption{Backscattered immersed data of the second component $E_2$ of electric field  for object number 7 (wooden doll, empty inside) of Table \ref{tab:table1} for different values of the parameter  $\beta$  in (\ref{5.109}). Recall that the final time is $T=1.2$.}
\label{fig:figure1_1}
\end{figure}

In our numerical studies we apply adaptive algorithm of section
\ref{subsec:ad_alg} to refine shape for nine (9) targets listed in Table
\ref{tab:table1}.  Three of them (targets number 1, 2, 5) were dielectrics,
three were metallic objects (targets number 4, 5, 6) and three (targets number 7, 8, 9)
were dolls with different objects placed inside them.  Heterogeneous
targets present models for explosive devices in which explosive
materials are masked by dielectrics.  Target number 7 was a wooden
doll which was empty inside, target number 8 was a piece of a metal
inserted inside that doll, and in target number 9 dry sand was partly
inserted inside the doll.

In all tests we have used the regularization parameter $\gamma =
0.01$.  The target number 1 was used for the calibration to choose
optimal parameters for dielectrics. We choose the cut-off number in
(\ref{8.6}) to be $0.85$ and the number $\beta =0.5$ in (\ref{5.109})  for all dielectrics of Table \ref{tab:table1}. The
target number 4 was used for calibration purpose of metallic targets:
we choose the cut-off number 0.3 in (\ref{8.7}).  In the case of
metals we have used the same numbers $\beta =0.5$ in (\ref{5.109}) as
 for dielectrics.

\begin{figure}[h]
  \begin{center}
    \begin{tabular}{cc}
      \includegraphics[width=0.45\textwidth]{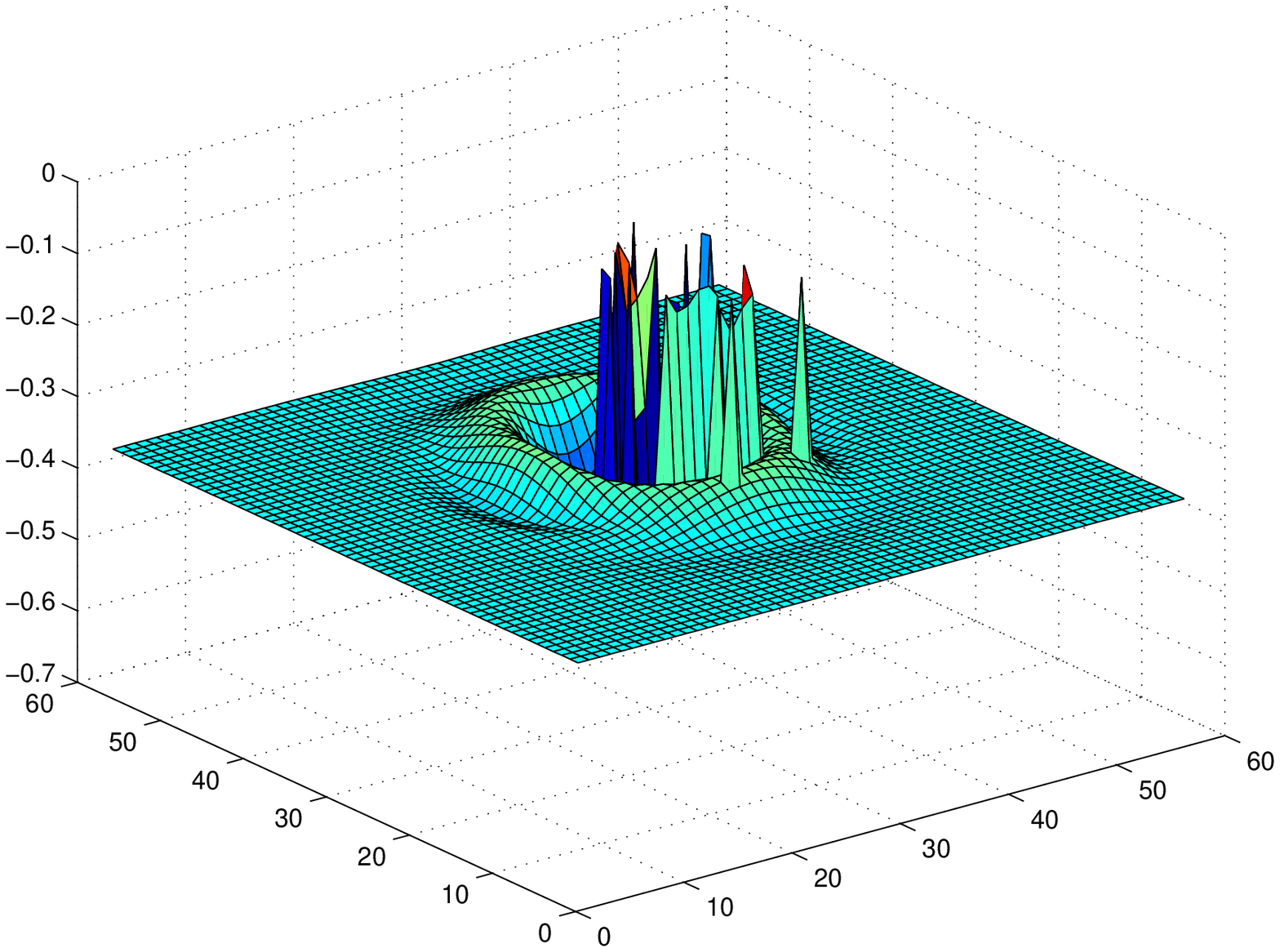}&
      \includegraphics[width=0.45\textwidth]{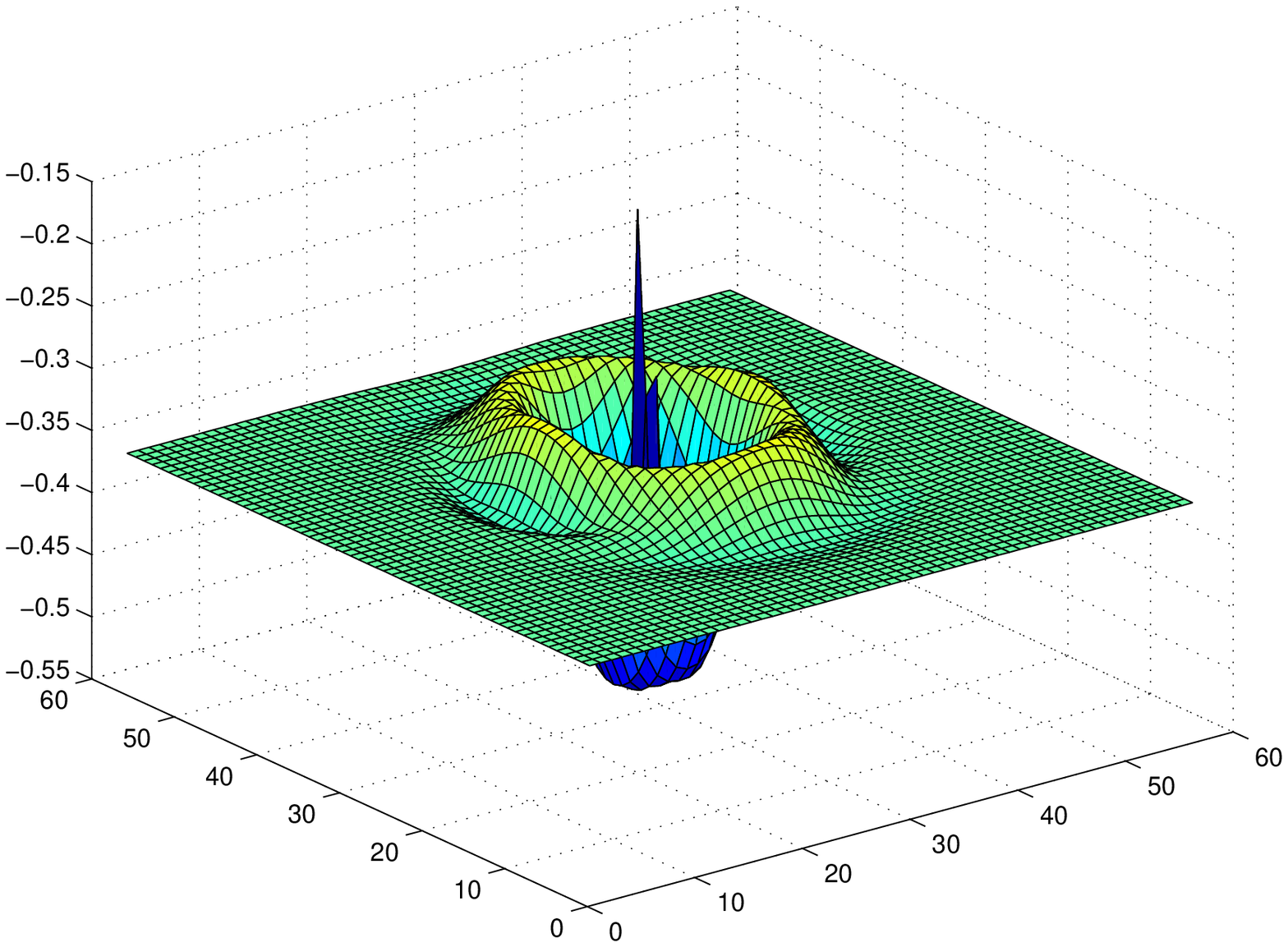}\\
      (a) $\beta = 0.2,t = 0.45 $ & (b) $\beta = 0.5, t = 0.45$\\
      \includegraphics[width=0.45\textwidth]{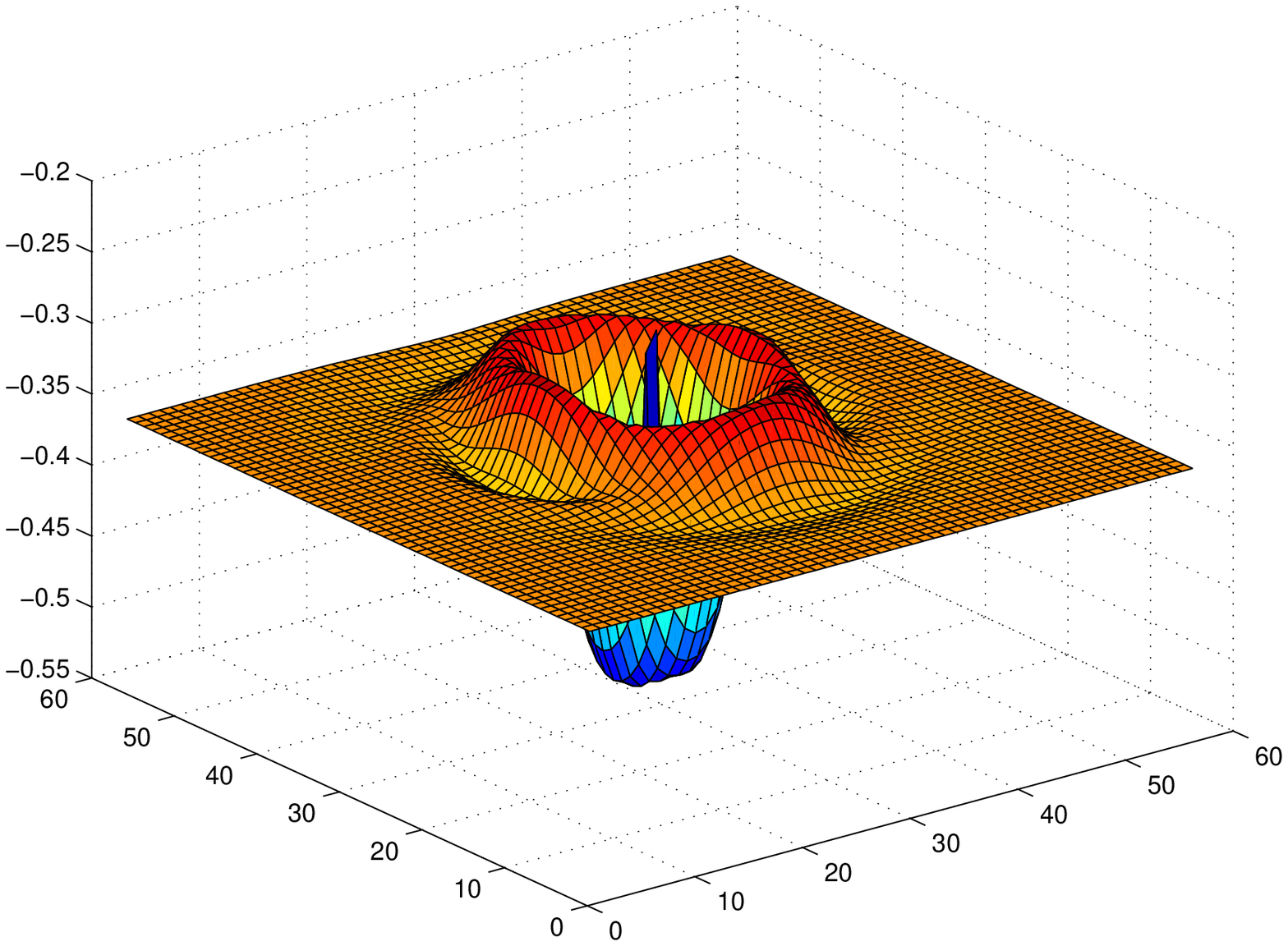}&
      \includegraphics[width=0.45\textwidth]{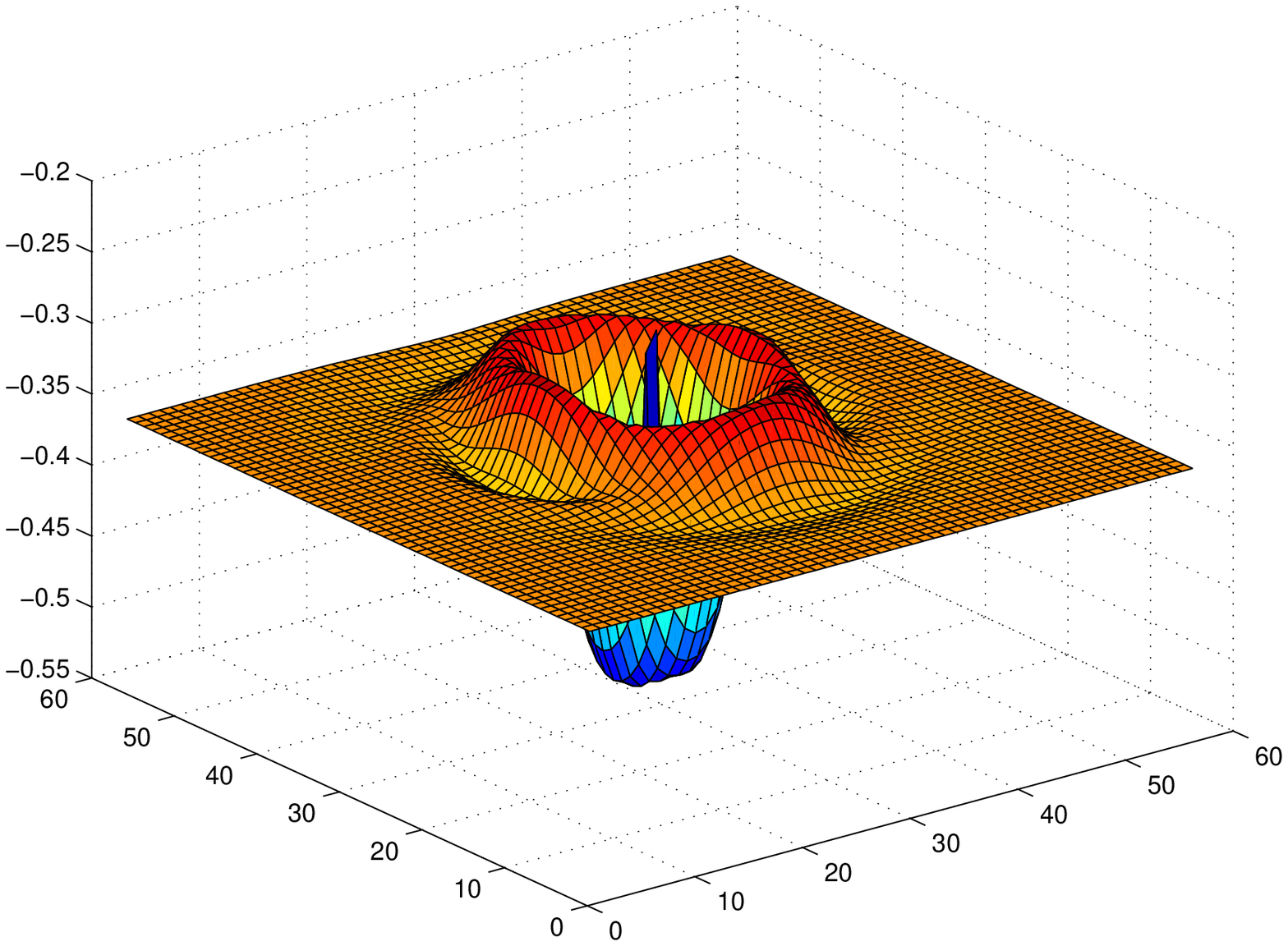}\\
      (c) $\beta = 0.8,t = 0.45 $ & (d) $\beta = 0.9, t = 0.45$
    \end{tabular}
 \end{center}
    \caption{Backscattered immersed data of the second component $E_2$ of electric field  for object 7 (wooden doll, empty inside) of Table \ref{tab:table1} for different values of the parameter  $\beta$  in (\ref{5.109}). Recall that the final time is $T=1.2$.}
\label{fig:figure1}
\end{figure}

  Figures \ref{fig:figure1_1}, \ref{fig:figure1}
show backscattered immersed data of the second component $E_2$ of the electric
field for object 7
at different times and with different immersing factor $\beta$ in
(\ref{5.109}).

In the case of dielectric targets we have \emph{a posteriori }directly
measured their refractive indices $n=\sqrt{\varepsilon _{r}}.$ Let
$\varepsilon _{r}^{comp}\left( \mathbf{x}\right) $ be the computed
coefficient. We consider
maximal values of functions $\varepsilon
_{r}^{comp}\left(\mathbf{x}\right).$ This means that in our Tables
\ref{tab:table2}, \ref{tab:table3}  we list values of dielectric
constants $\varepsilon _{r}\left( \text{target}\right) $ and
refractive indices $n\left( \text{target}\right) $ as
\begin{equation*}
\varepsilon _{r}\left( \text{target}\right) =\max_{\overline{\Omega }%
}\varepsilon _{r}^{comp}\left( \mathbf{x}\right) ,n\left( \text{target}%
\right) =\sqrt{\varepsilon _{r}\left( \text{target}\right) }.
\end{equation*}%

Tables \ref{tab:table2}, \ref{tab:table3} are quite informative
ones since they show the accuracy of our reconstruction of
either refractive indices (Table \ref{tab:table2}) or effective dielectric constants (Table \ref{tab:table3}).  Table \ref{tab:table2} lists refractive indices of dielectric
targets, both computed $n\left( \text{target}\right) $ and directly
measured ones $n$.  Computed numbers $n\left( \text{target}\right) $
are displayed for different locally refined meshes. This table also
shows the measurement error in direct measurements of $n$. Table
\ref{tab:table3} lists calculated effective dielectric constants
$\varepsilon _{r}\left( \text{target}\right) $ of metallic targets,
again for different locally refined meshes. In Table \ref{tab:table3}
 for all metallic
targets we have $\varepsilon \left( \text{target}\right) >10$, which means that (\ref{2.51}) is satisfied.

One can derive several important observations from Table \ref{tab:table2}.
First, for all targets and on all adaptively refined meshes the
computational error is significantly less than the error of direct
measurements. Thus, the average computational error is significantly less
than the average measurement error on all adaptively refined meshes. Second,
computed refractive indices are within reasonable error estimates in all
cases. The accuracy on all adaptively refined meshes is about the same.

We observe from Table \ref{tab:table3}  that for target number 8 we
have obtained effective dielectric constant $\varepsilon \left( \text{target}%
\right) \in \lbrack 13.6,14]$ on all adaptively refined meshes, which is
less than for other metallic objects. We can explain this by the fact that
target number 8 is a mixture of metal and dielectric. An important
observation, which can be deduced from Table \ref{tab:table3}, is that our
two-stage algorithm can still compute large inclusion/background contrasts
exceeding 10:1, just as the algorithm of the first stage.

Figures \ref{fig:fig3}--\ref{fig:fig14} display 3-d images of some
targets listed in Table \ref{tab:table1} as well as corresponding adaptively locally
refined meshes.  To have a better visualization, we have zoomed some
figures to $0.4\times 0.4$ square from $1\times 1$ square in the $x,y$
directions. So, in these Figures we display only the image in the domain $%
\Omega _{zoom},$ 
\begin{equation}
\Omega _{zoom}=\left\{ \mathbf{x=}(x,y,z)\in (-0.2,0.2)\times
(-0.2,0.2)\times (-0.1,0.04)\right\} ,  \label{zoom}
\end{equation}%
Figures \ref{fig:fig3}-- \ref{fig:fig14} also show estimates of sizes of the
targets in the $z$-direction. Locations of all targets as well as their
sizes in $x,y,z$ directions are well estimated.

\begin{figure}[h]
  \begin{center}
    \begin{tabular}{ccc}
      \includegraphics[angle=-90,width = 5.5cm]{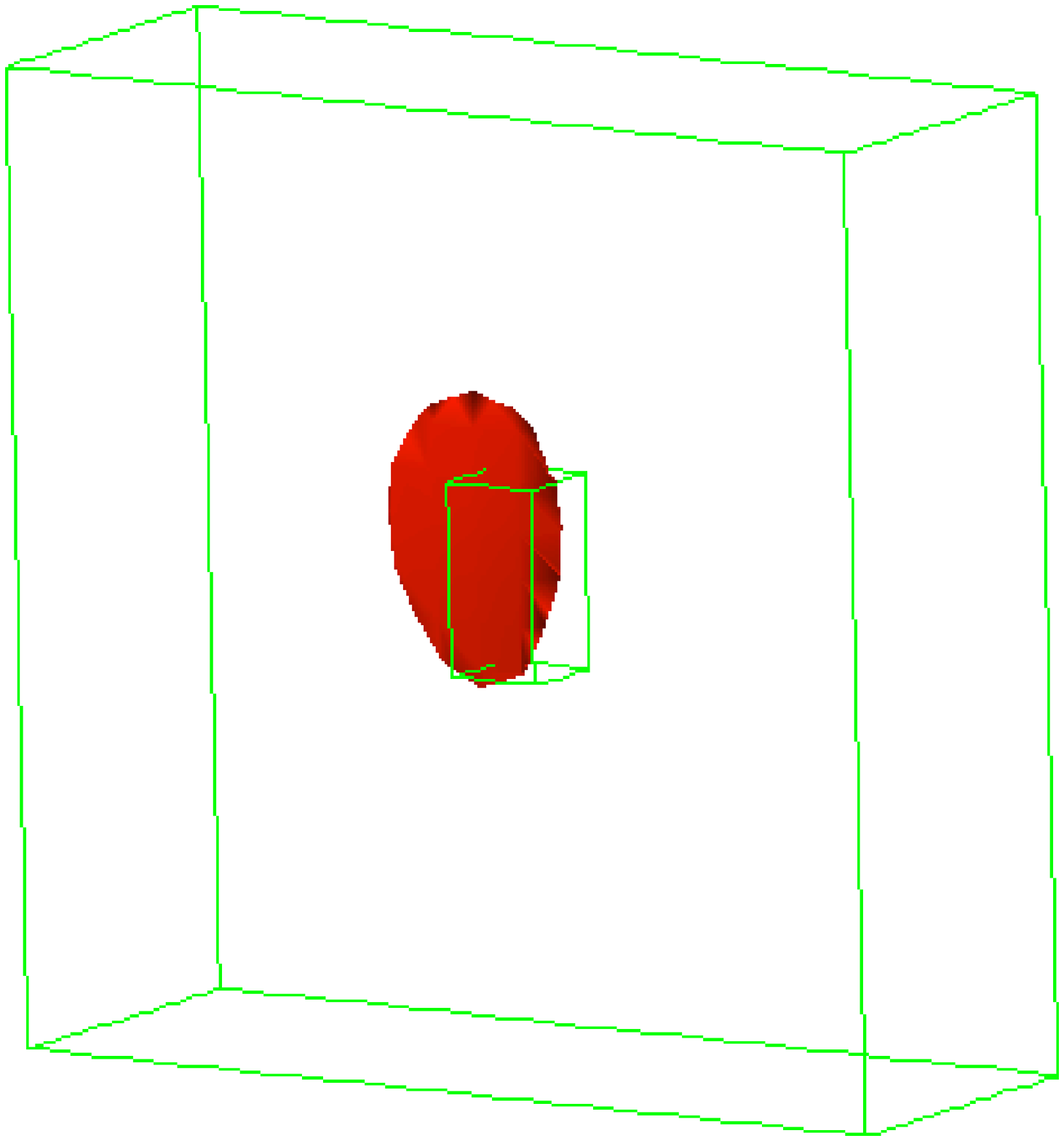}&
      \includegraphics[angle=-90,width = 5.5cm]{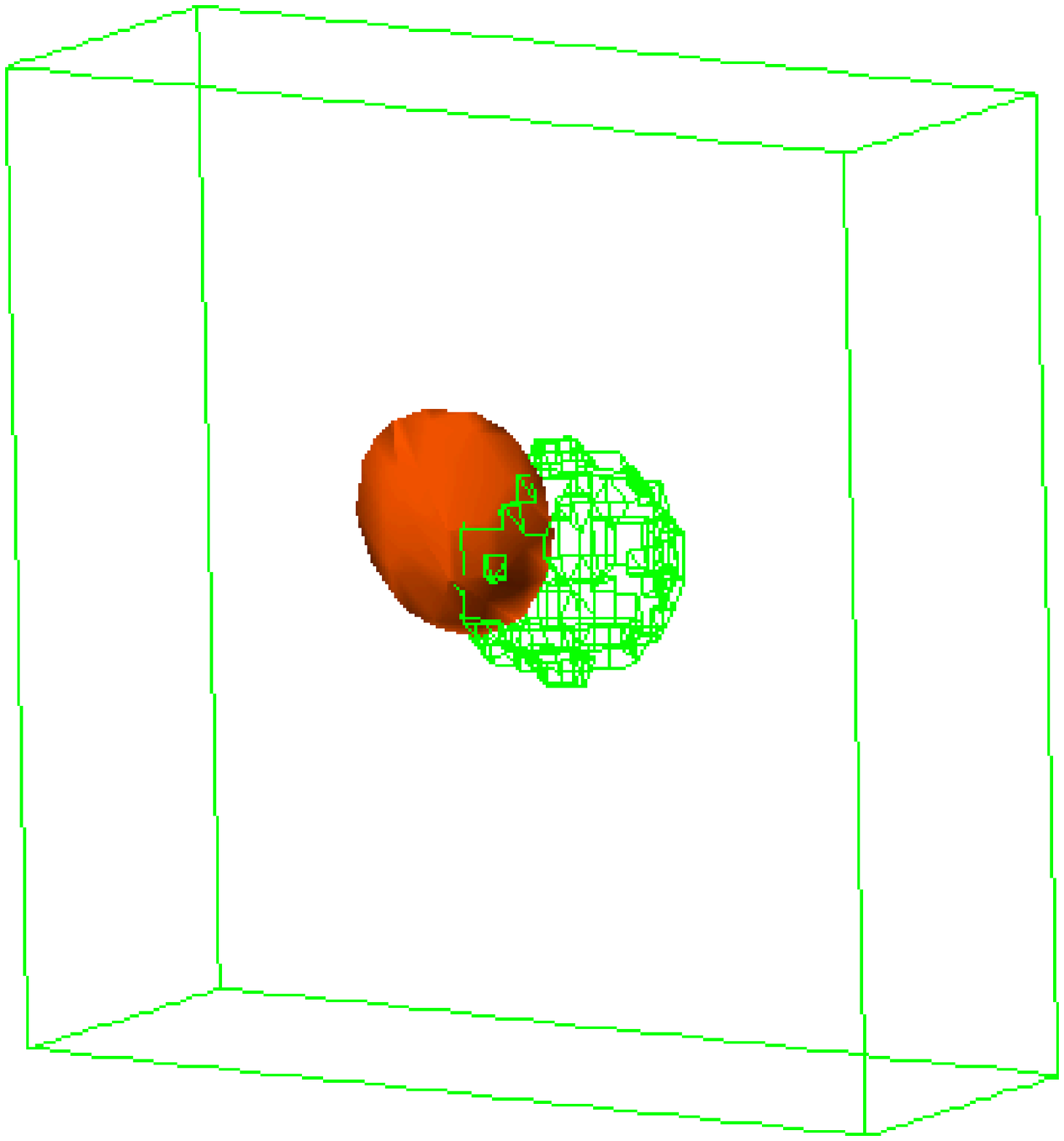} &
   \includegraphics[angle=-90,width = 5.5cm]{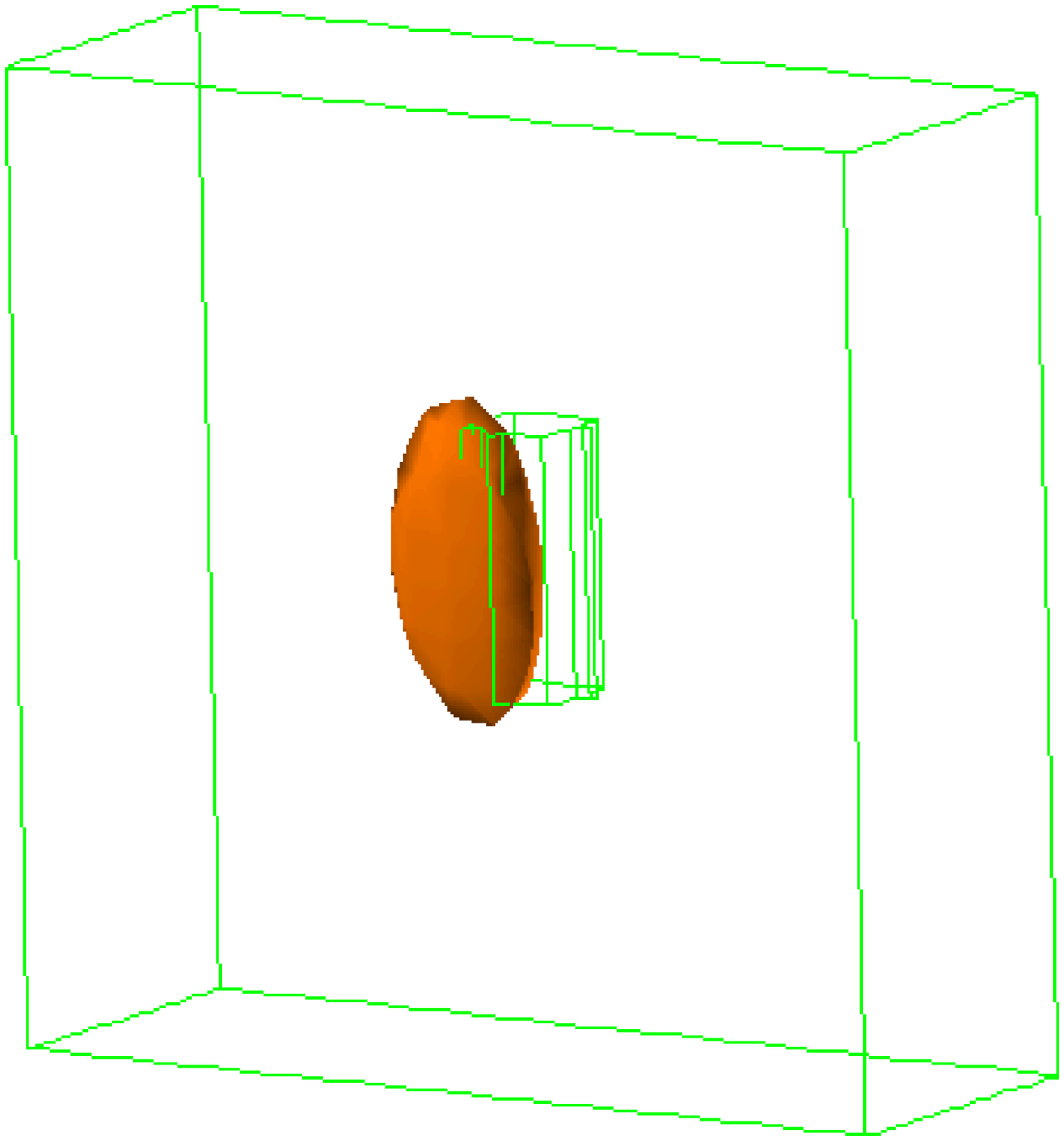} \\
a) target 1  &  b) target 3 & c)  target 4\\
  \includegraphics[angle=-90, width = 5cm]{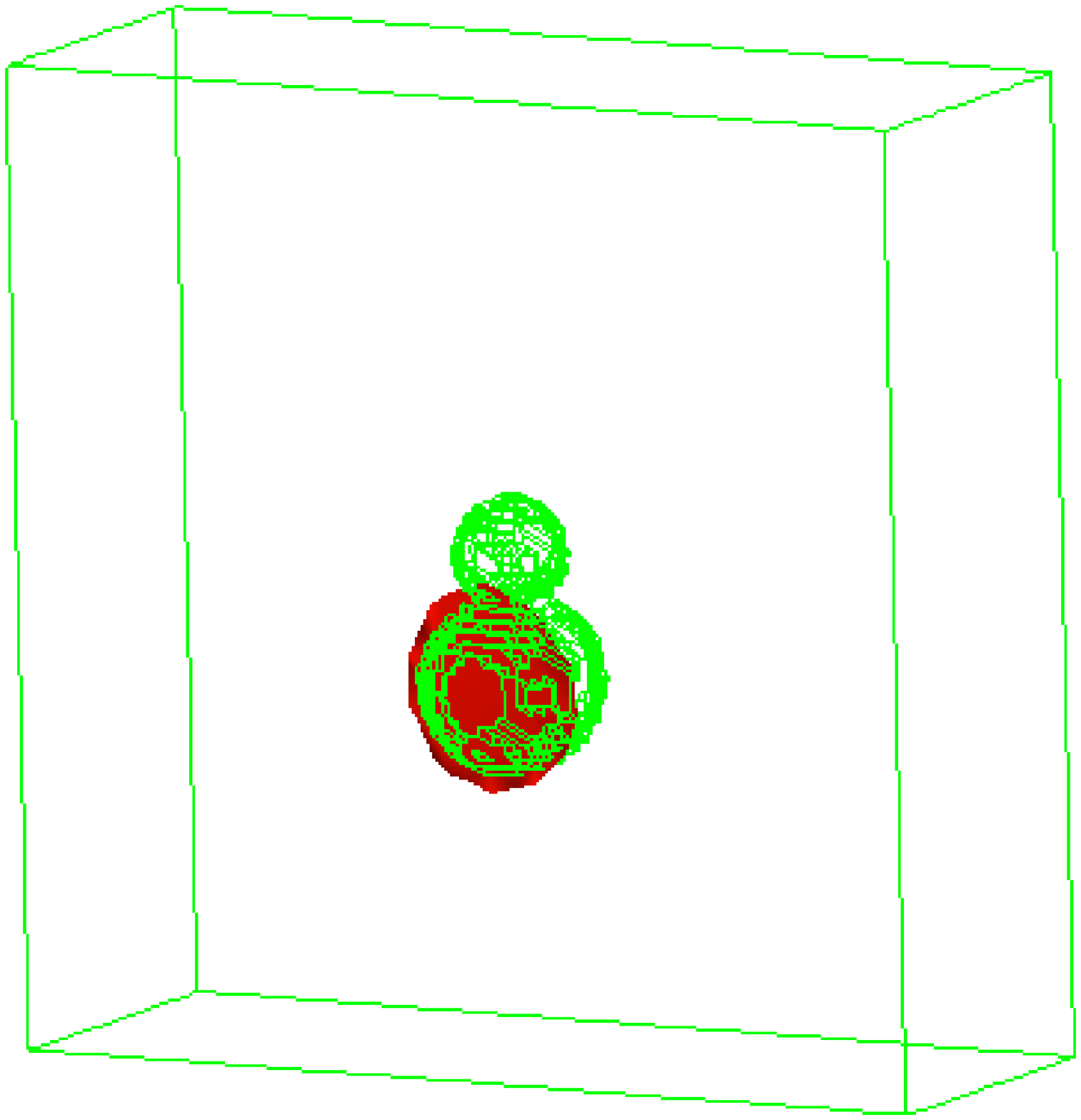} &
  \includegraphics[angle=-90,width = 5cm]{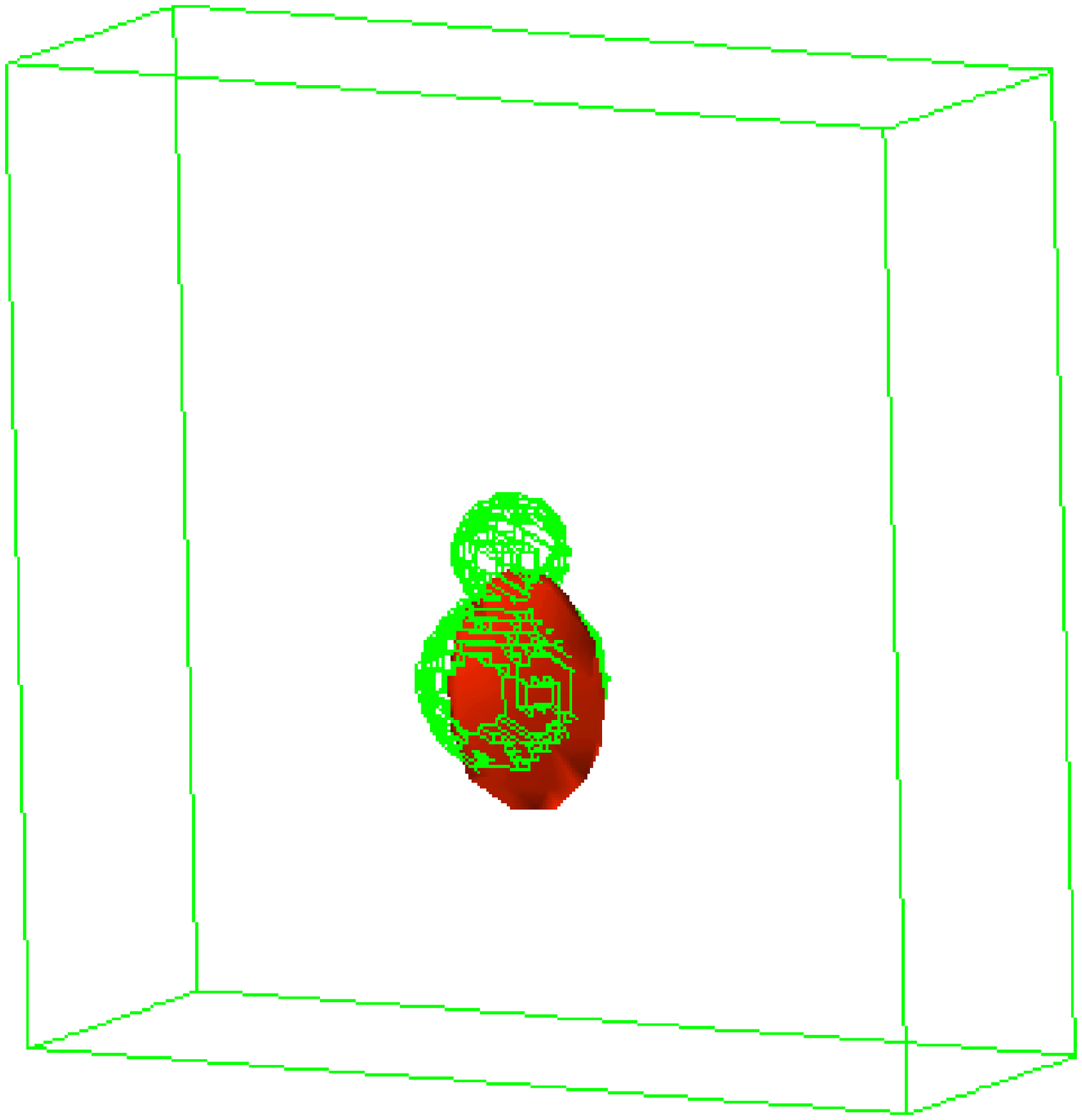} &
  \includegraphics[angle=-90,width = 5cm]{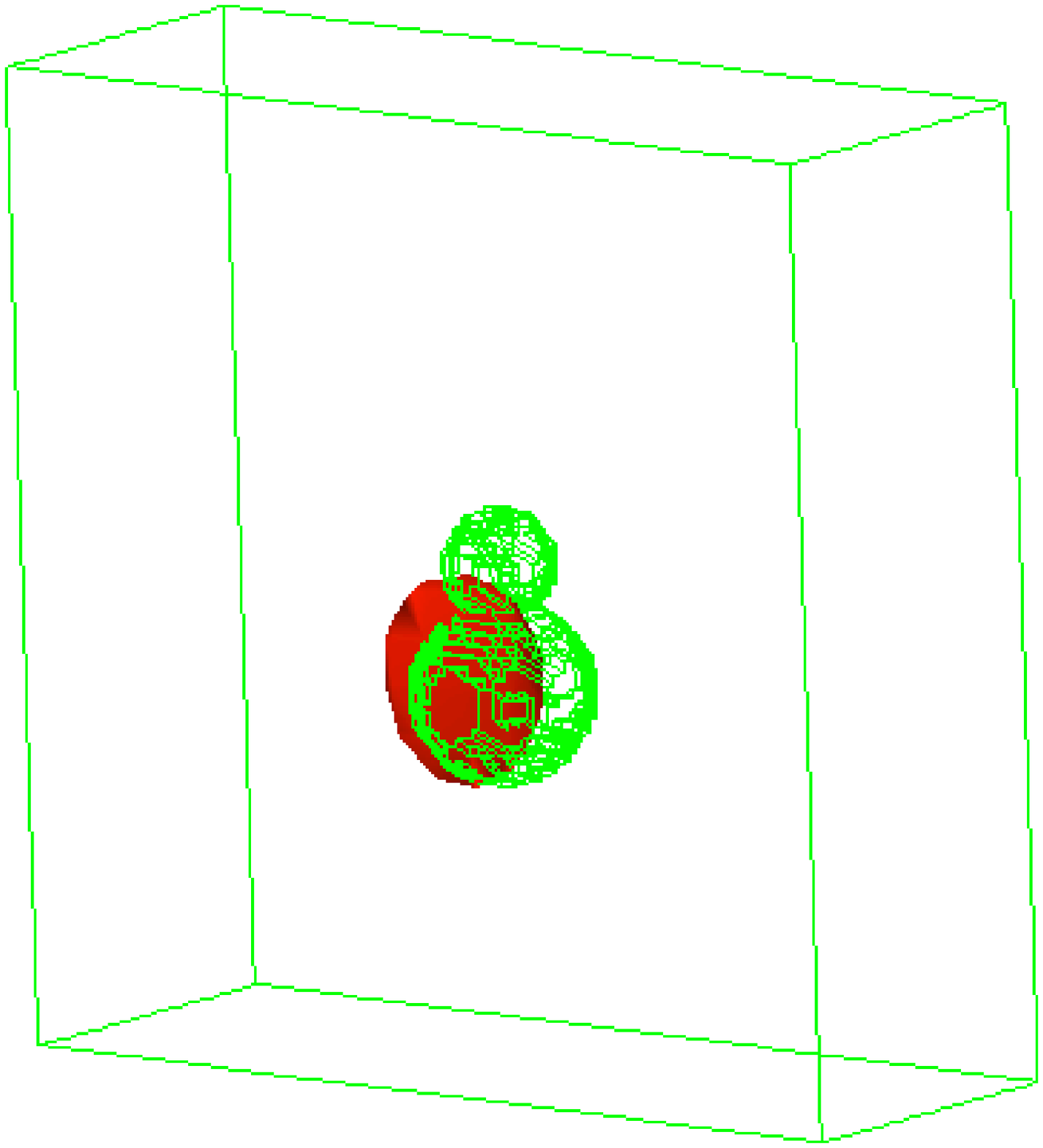} 
 \\
d) target 7  &  e) target 8  & f)  target 9
    \end{tabular}
 \end{center}
    \caption{Reconstructions of some targets of Table \ref{tab:table1}
      obtained in \cite{BTKF, NBKF} on the first stage of our two-stage numerical
      procedure.}
\label{fig:globconv}
\end{figure}

\section{Summary}\label{sec:summary}

We have used time dependent backscattering experimental data generated by a
single source of electric pulses to simultaneously reconstruct all three parameters
of interest of explosive-like targets: their refractive indices, shapes and
locations. To do this, we have used the two-stage reconstruction procedure, which
was first proposed in \cite{BK2}. On the first stage the globally
convergent  method of \cite{BK1,BK} was used. This method has
provided accurate estimates of refractive indices and locations of targets 
\cite{BTKF,NBKF}. On the second stage,
which is the focus of this paper, images were refined using the
adaptivity technique of  \cite{BMaxwell2}. The second stage has provided
accurate estimates of the third component: the shape, in addition to the
first two ones. Interestingly, even heterogeneous targets, which model
heterogeneous IEDs, were quite accurately imaged.

In all cases one can observe a significant improvement of
the image quality after the application of the adaptivity on the second
stage. Another observation here is that we can accurately image shapes of not
only targets with ''straight" boundaries, like a rectangular prism
(target number 1) or a cylinder (target number 4) but targets with
curved boundaries as well (targets number 3 and 7). Even shadow parts
of targets number $1, 3, 4, 8$ are imaged rather well. This is regardless
on the fact that we use the minimal possible information content and
on a narrow view angle: single location of the source and time
resolved backscattering data.



\begin{figure}
  \begin{center}
    \begin{tabular}{cccc}
      \includegraphics[width=5cm, clip=true, trim = 11cm 4cm 11cm 4cm]{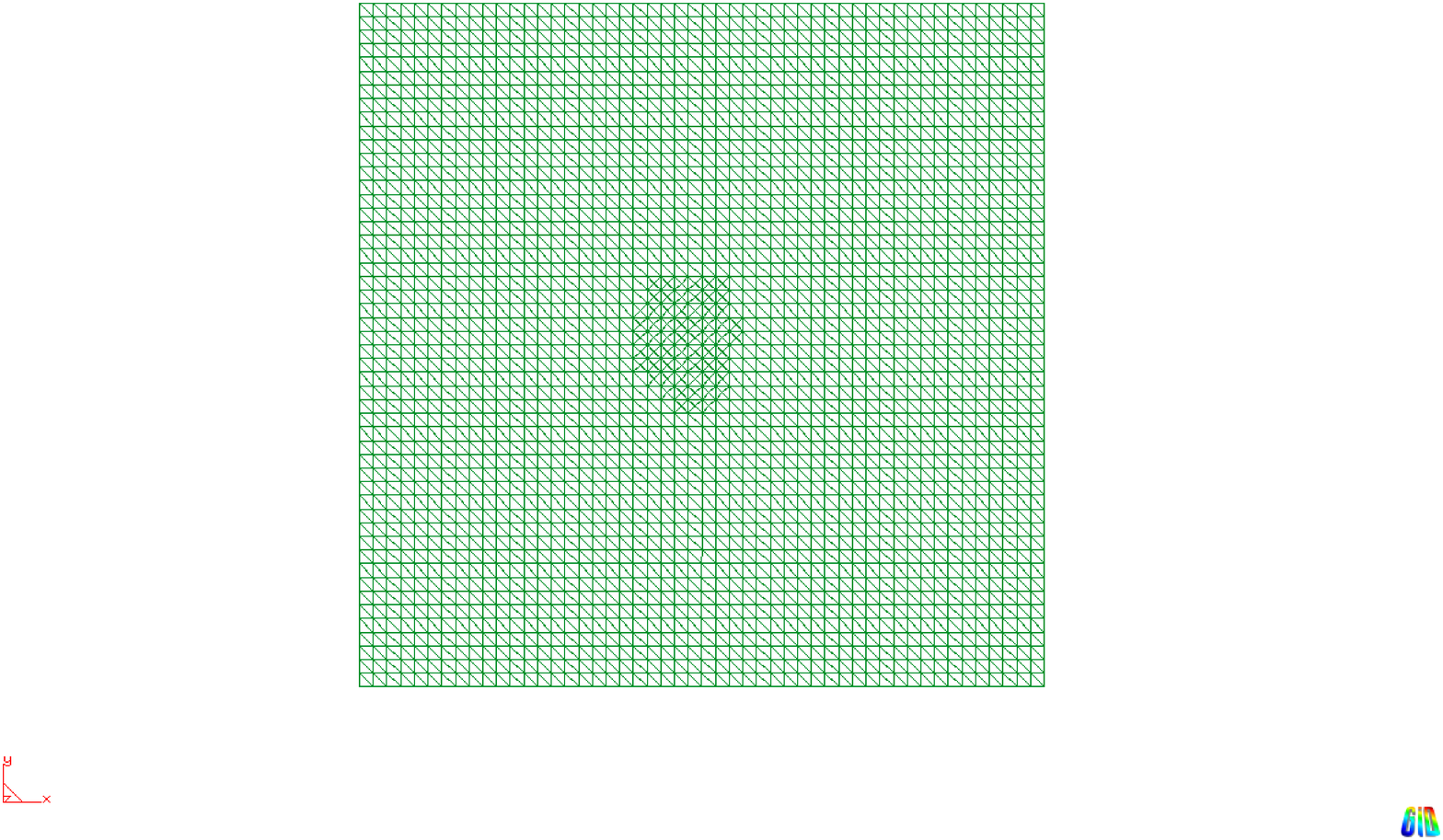}&
      \includegraphics[width=5cm, angle = 90, clip = true, trim = 11cm 12.5cm 11cm 12.5cm]{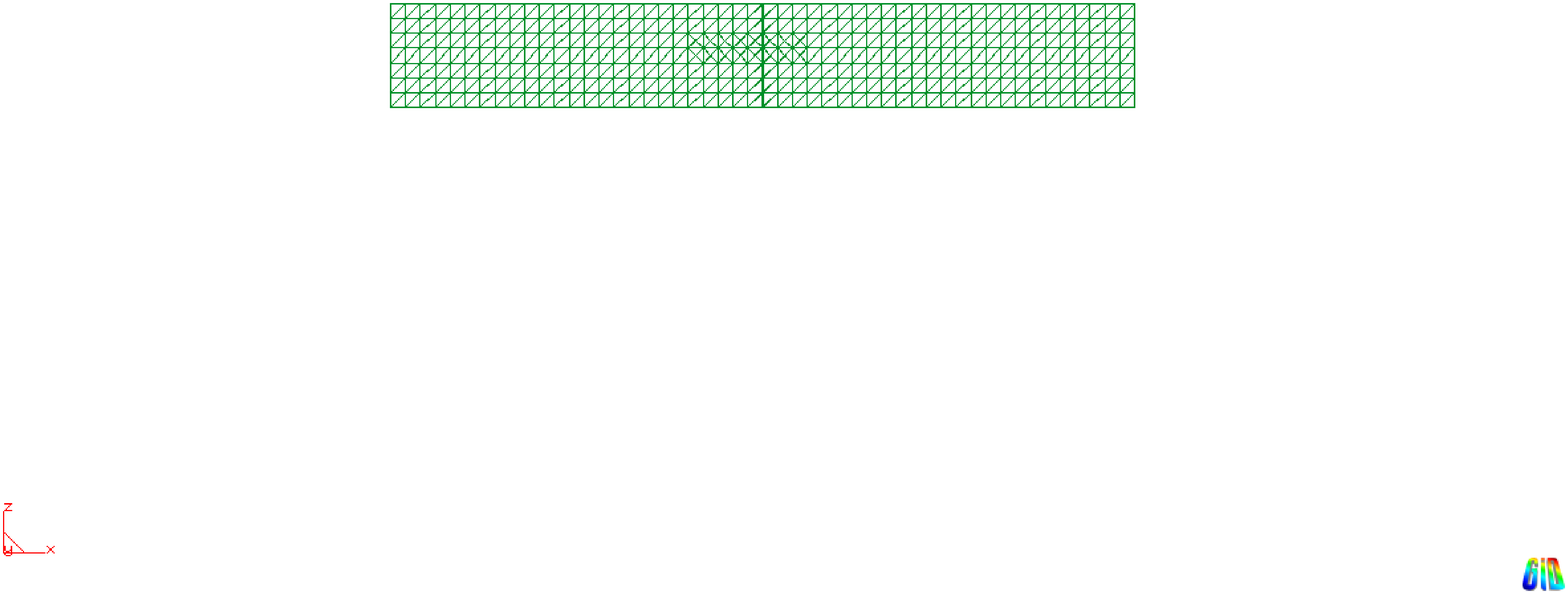}&
      \includegraphics[width=5cm, angle = 90, clip = true, trim = 11cm 12.5cm 11cm 12.5cm]{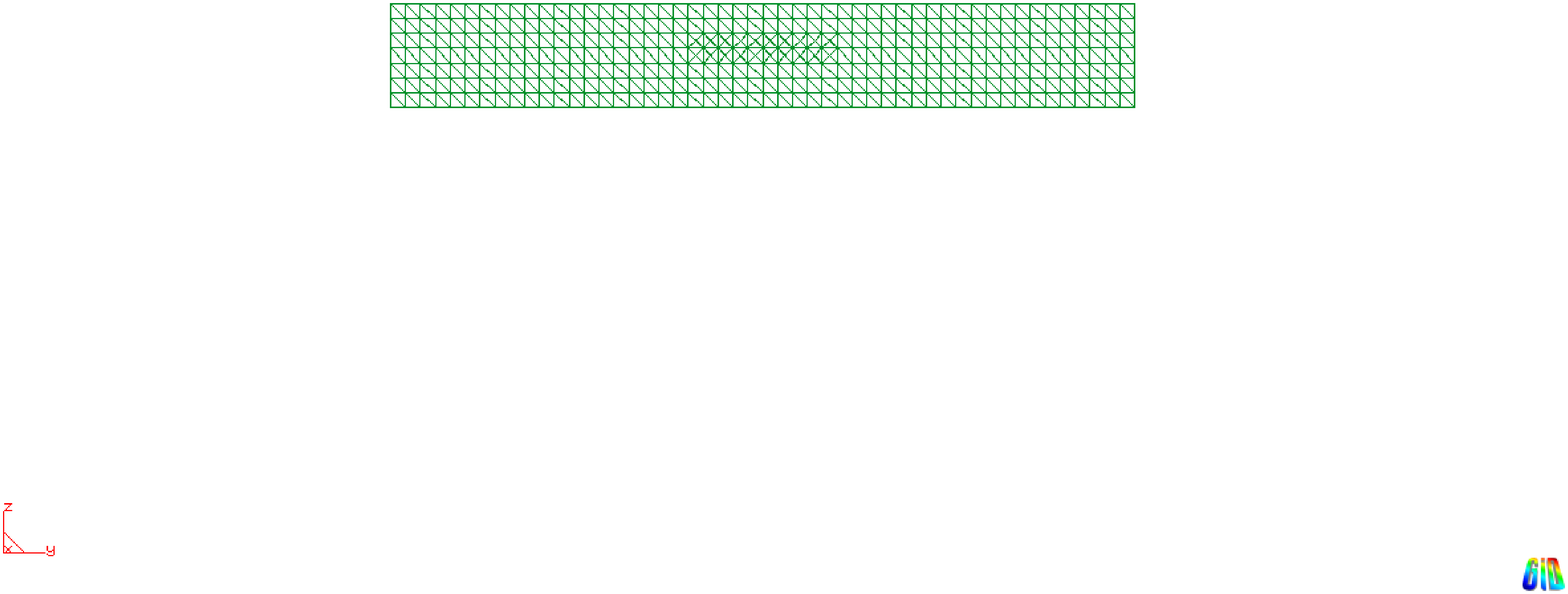}&
      \includegraphics[width=4.0cm, clip = true, trim = 1.8cm 0.0cm 4.0cm 0.0cm]{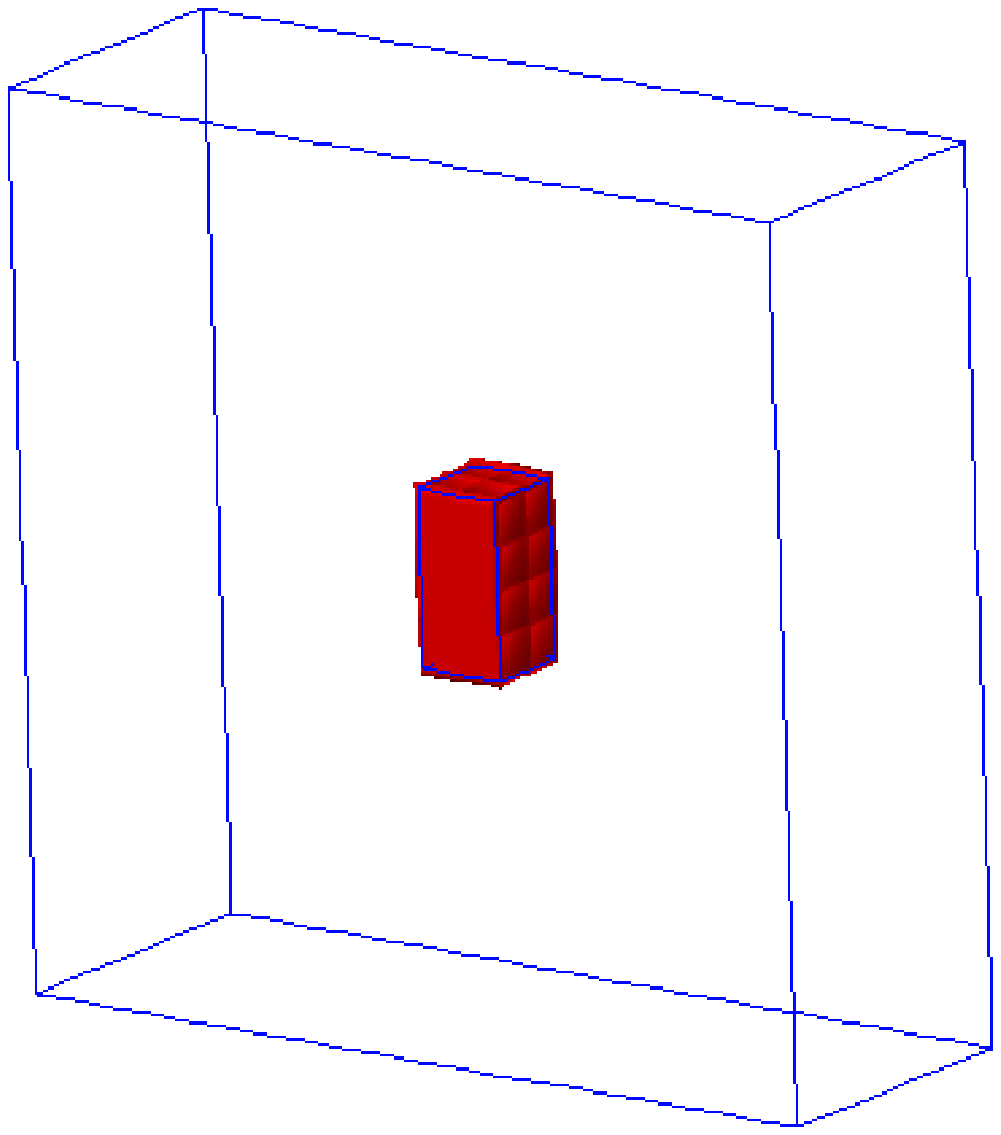}\\
      (a) &
      (b) &
      (c) &
      (d) 
    \end{tabular}
 \end{center}
    \caption{(a) $xy$-projection, (b) $xz$-projection, and (c)
      $yz$-projection of the once refined (optimal) mesh; d) Computed
      image of target number 1 of Table \ref{tab:table1} on that
      mesh. Thin lines indicate correct shapes.  To have a better 
      visualization we have zoomed the domain $\Omega$ in (\ref{8.0})
      in the domain $\Omega_{zoom}$ in (\ref{zoom}).  This target
      number 1 was used for the calibration purpose for the case of
      dielectric targets. A significant improvement of the image of
      d), compared with the image of Figure \ref{fig:globconv}-a), obtained on the
      first stage is evident.}
\label{fig:fig3}
\end{figure}

\begin{figure}[ht!]
  \begin{center}
    \begin{tabular}{cccc}
      \includegraphics[width=5cm, clip=true, trim = 11cm 4cm 11cm 4cm]{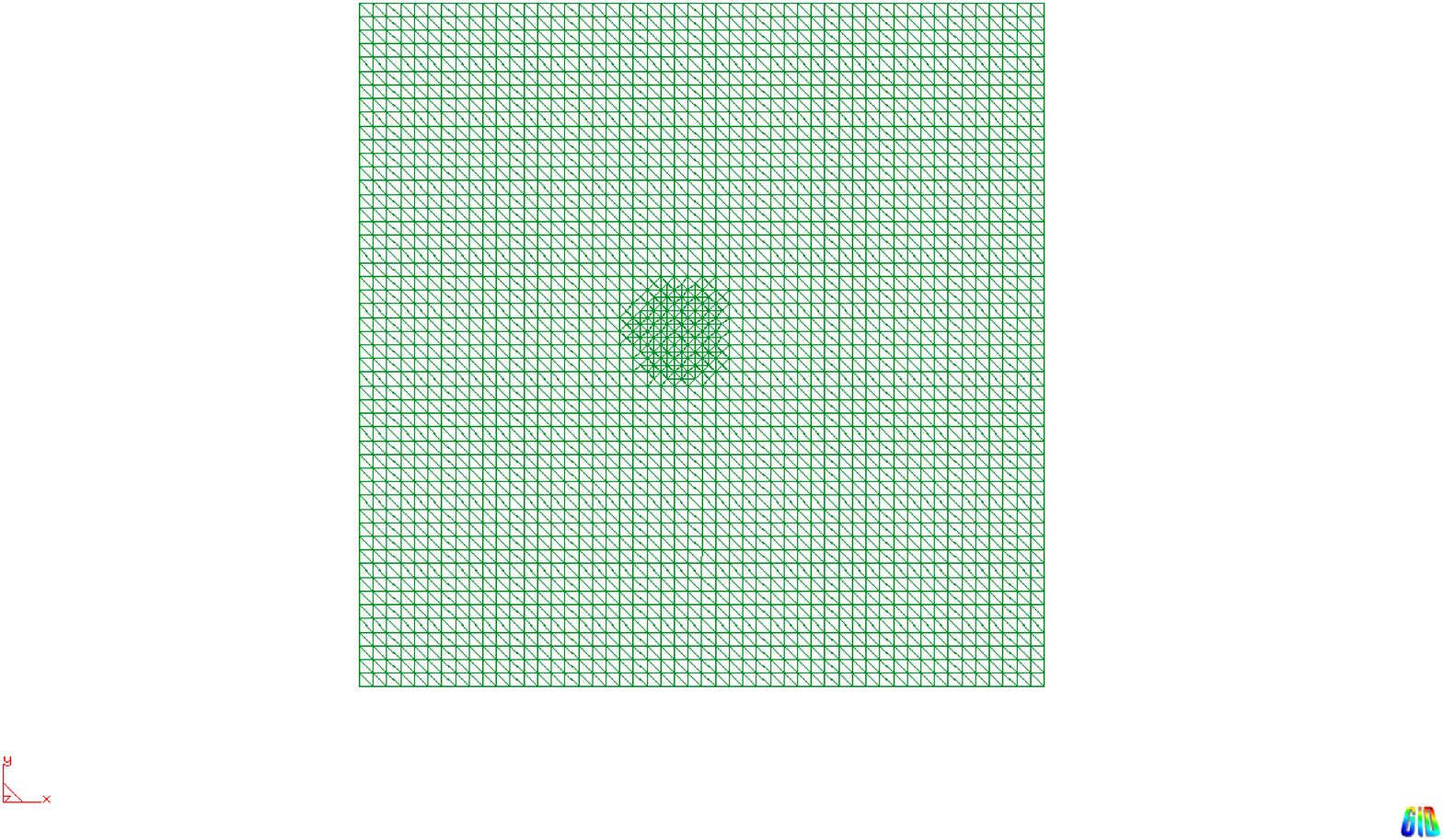}&
      \includegraphics[width=5cm, angle = 90, clip = true, trim = 11cm 12.5cm 11cm 12.5cm]{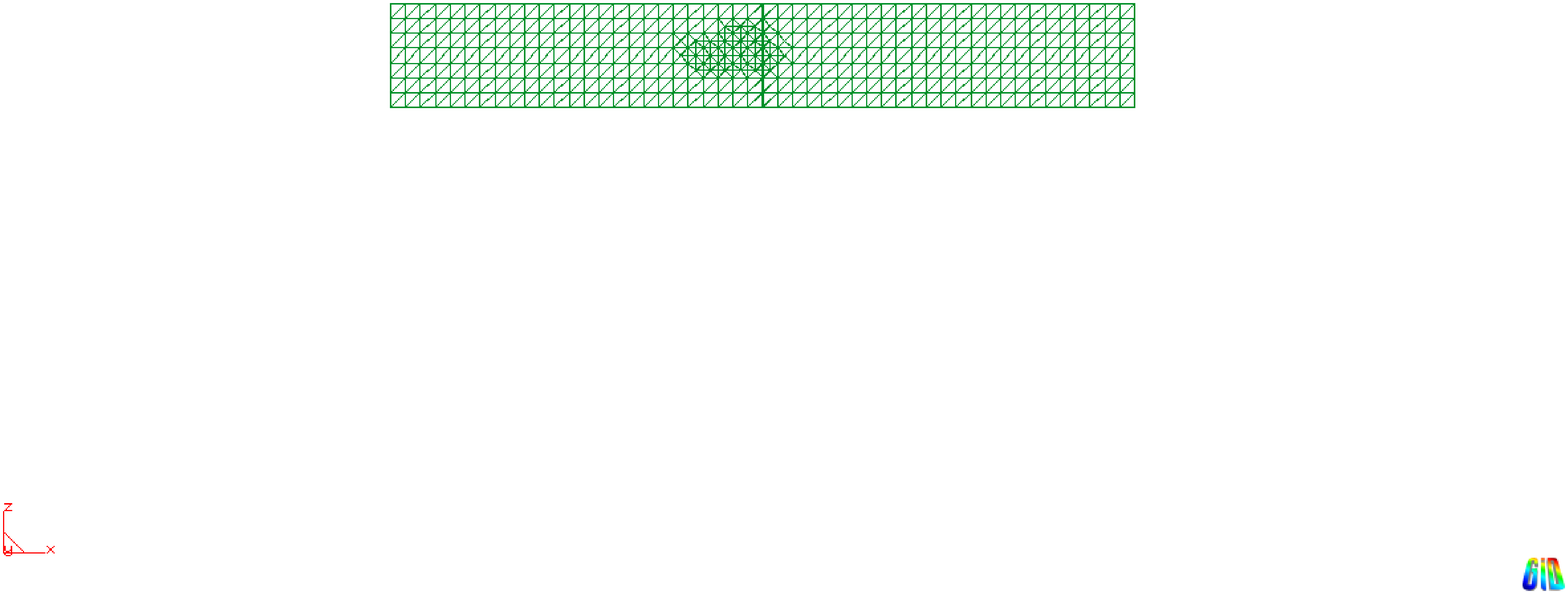}&
      \includegraphics[width=5cm, angle = 90, clip = true, trim = 11cm 12.5cm 11cm 12.5cm]{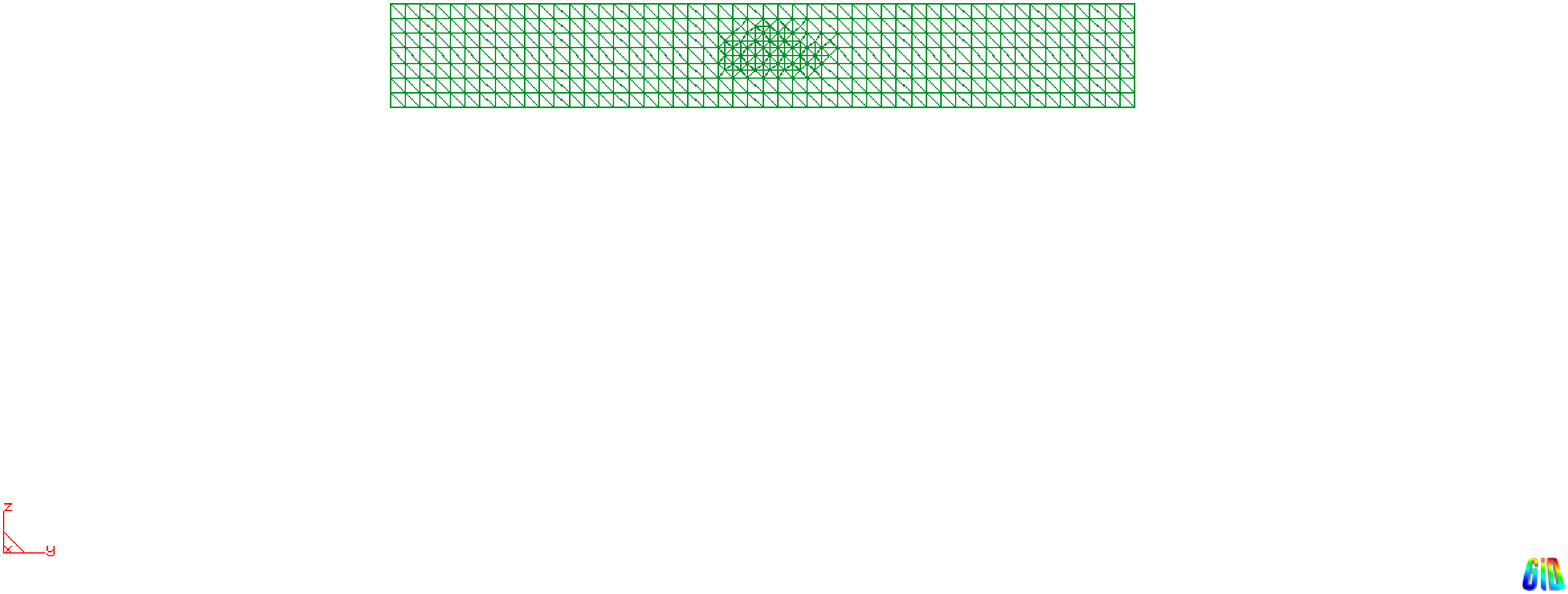}&
      \includegraphics[width=4.0cm, clip = true, trim = 1.8cm 0.0cm 3.5cm 0.0cm]{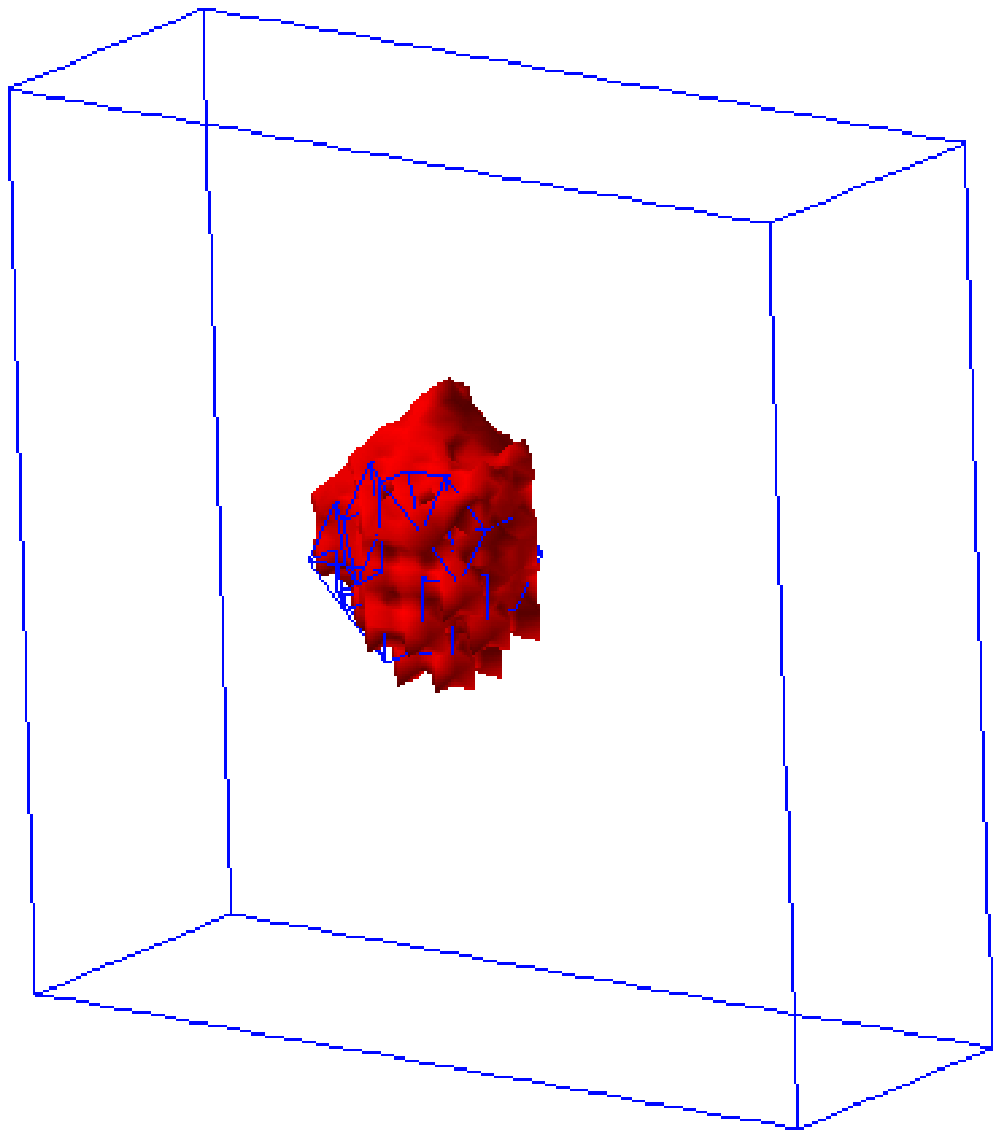}\\
      (a) &
      (b) &
      (c) &
      (d) 
    \end{tabular}
 \end{center}
    \caption{(a) $xy$-projection, (b) $xz$-projection, and (c)
      $yz$-projection of the three times refined (optimal) mesh; d) Computed
      image of target number 3 of Table \ref{tab:table1} on that
      mesh. Thin lines indicate correct shape.  To have better a
      visualization we have zoomed the domain $\Omega$ in (\ref{8.0})
      in the domain $\Omega_{zoom}$ in  (\ref{zoom}).  Comparison with
Figure \ref{fig:globconv}-b) shows a significant improvement compared with the first stage. }
\label{fig:fig5}
\end{figure}

\begin{figure}
  \begin{center}
    \begin{tabular}{ccc}
      \includegraphics[width=4.0cm, clip = true, trim = 1.8cm 0.0cm 4.0cm 0.0cm]{1/eps_perspective.eps} &
      \includegraphics[width=4.0cm, clip = true, trim = 1.8cm 0.0cm 4.0cm 0.0cm]{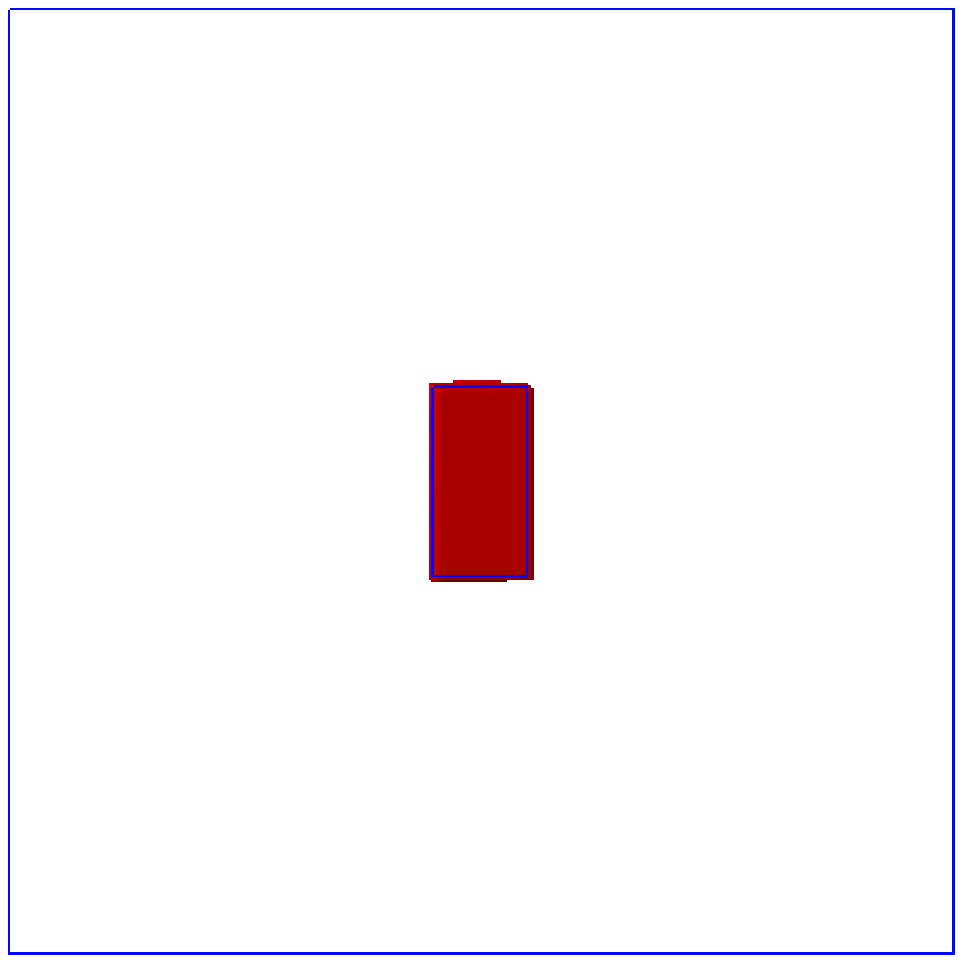} &
      \includegraphics[width=4.0cm, clip = true, trim = 1.8cm 0.0cm 4.0cm 0.0cm]{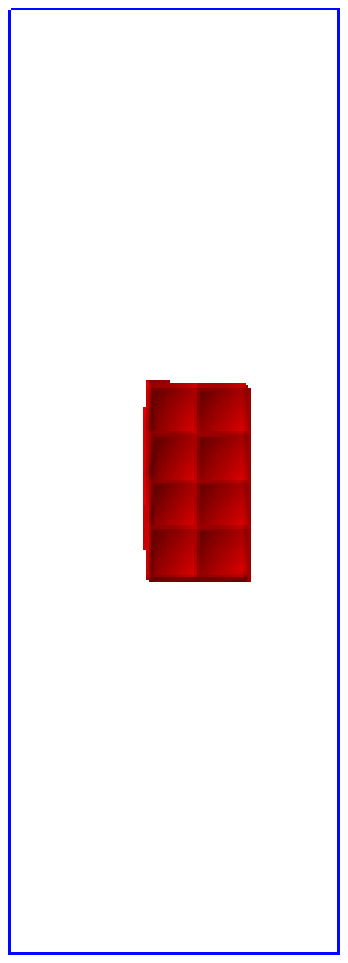}\\
      (a) Perspective view&
      (b) Front view&
      (c) Side view\\
      \includegraphics[width=4.0cm, clip = true, trim = 2.0cm 0.0cm 4.0cm 0.0cm]{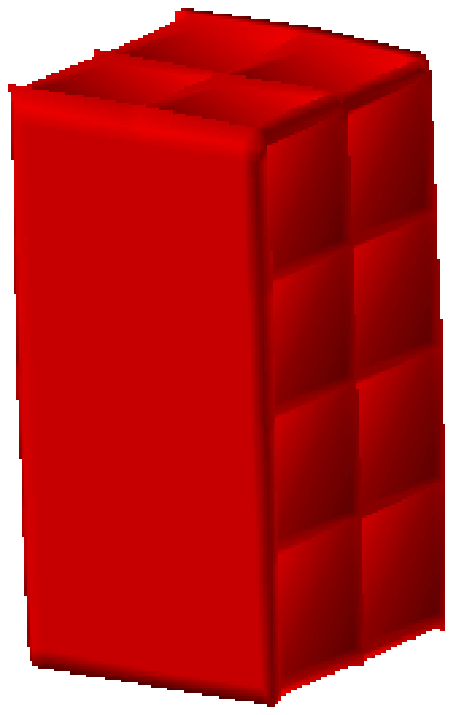} &
      \includegraphics[width=4.0cm, clip = true, trim = 2.0cm 0.0cm 4.0cm 0.0cm]{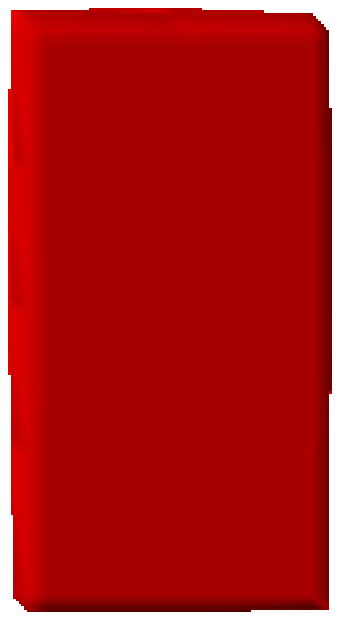} &
      \includegraphics[width=4.0cm, clip = true, trim = 2.0cm 0.0cm 4.0cm 0.0cm]{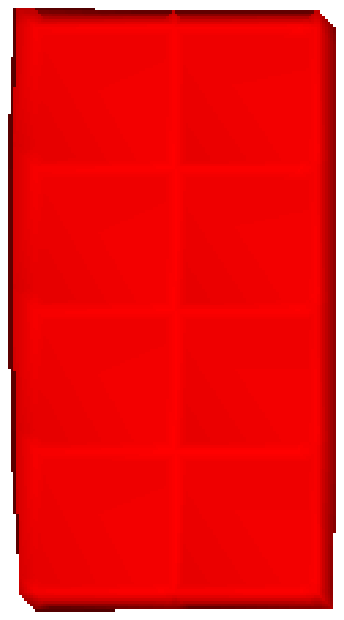}\\
      (d) Zoom, perspective&
      (e) Zoom, front&
      (f) Zoom, side\\
  \includegraphics[width=4.0cm, clip = true, trim = 1.8cm 0.0cm 3.5cm 0.0cm]{4/eps_perspective.eps} &
      \includegraphics[width=4.0cm, clip = true, trim = 1.8cm 0.0cm 4.0cm 0.0cm]{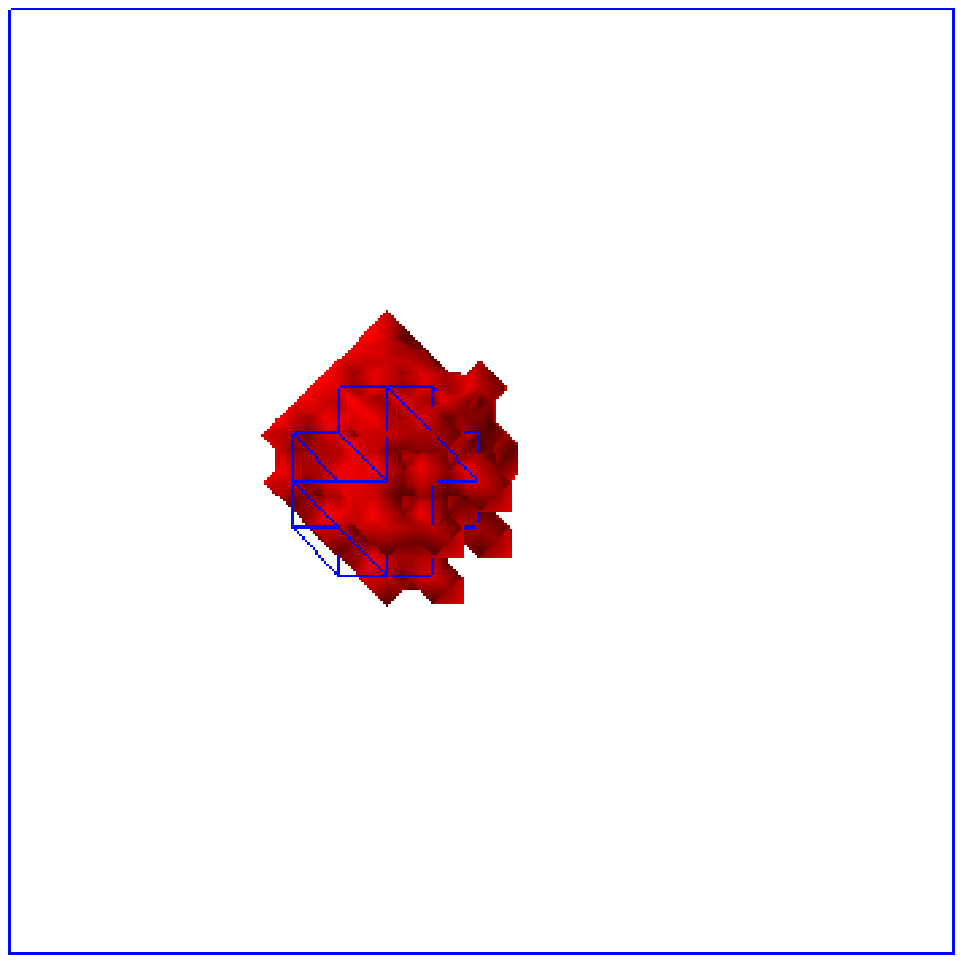} &
      \includegraphics[width=4.0cm, clip = true, trim = 1.8cm 0.0cm 4.0cm 0.0cm]{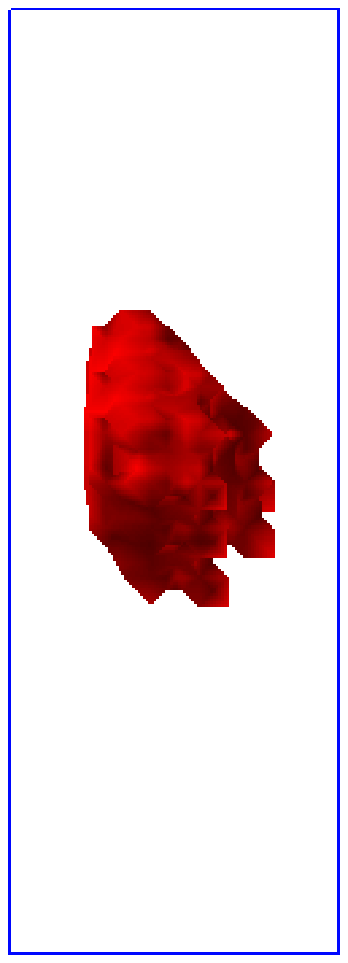}\\
      (g) Perspective view&
      (h) Front view&
      (i) Side view\\
      \includegraphics[width=4.0cm, clip = true, trim = 2.0cm 0.0cm 4.0cm 0.0cm]{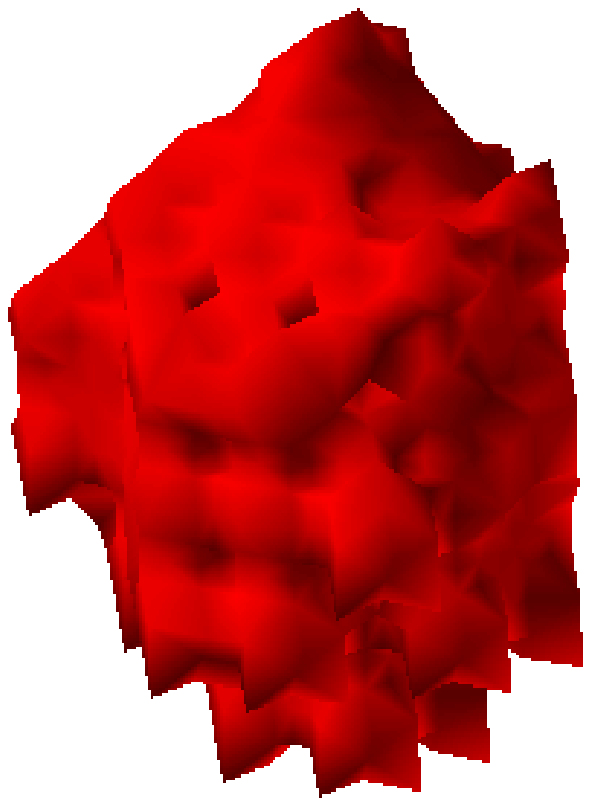} &
      \includegraphics[width=4.0cm, clip = true, trim = 2.0cm 0.0cm 4.0cm 0.0cm]{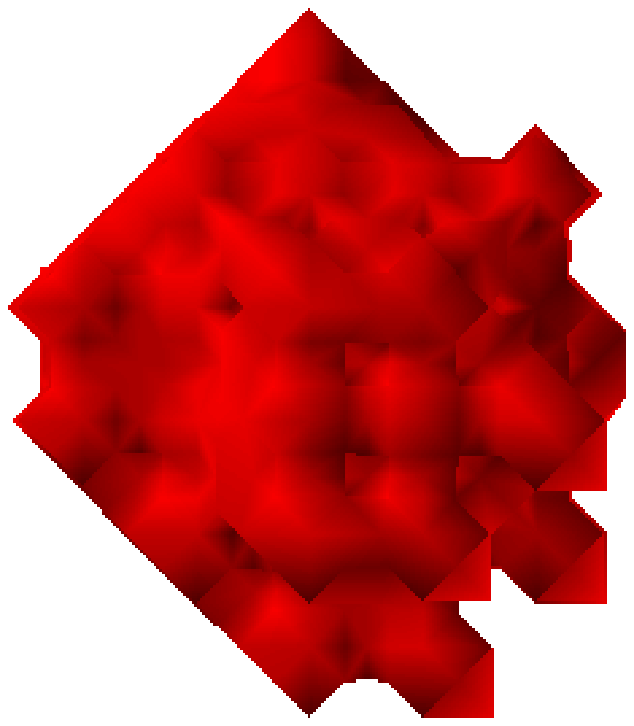} &
      \includegraphics[width=4.0cm, clip = true, trim = 2.0cm 0.0cm 4.0cm 0.0cm]{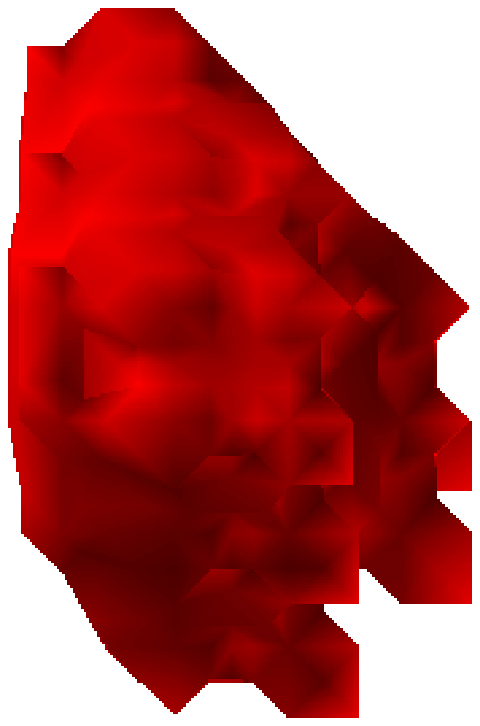}\\
      (j) Zoom, perspective&
      (k) Zoom, front&
      (l) Zoom, side
    \end{tabular}
 \end{center}
    \caption{Three views and zooms of the reconstruction of target number 1
      (figures a)-f)) on the once refined mesh. Three views and zooms
      of the reconstruction of target number 3 (figures g)-l)) of
      Table \ref{tab:table1} on three times refined mesh. Recall that
      target number 3 is a metallic sphere.}
\label{fig:fig4}
\end{figure}


\begin{figure}
  \begin{center} 
    \begin{tabular}{ccc}
      \includegraphics[width=5cm, clip=true, trim = 11cm 4cm 11cm 4cm]{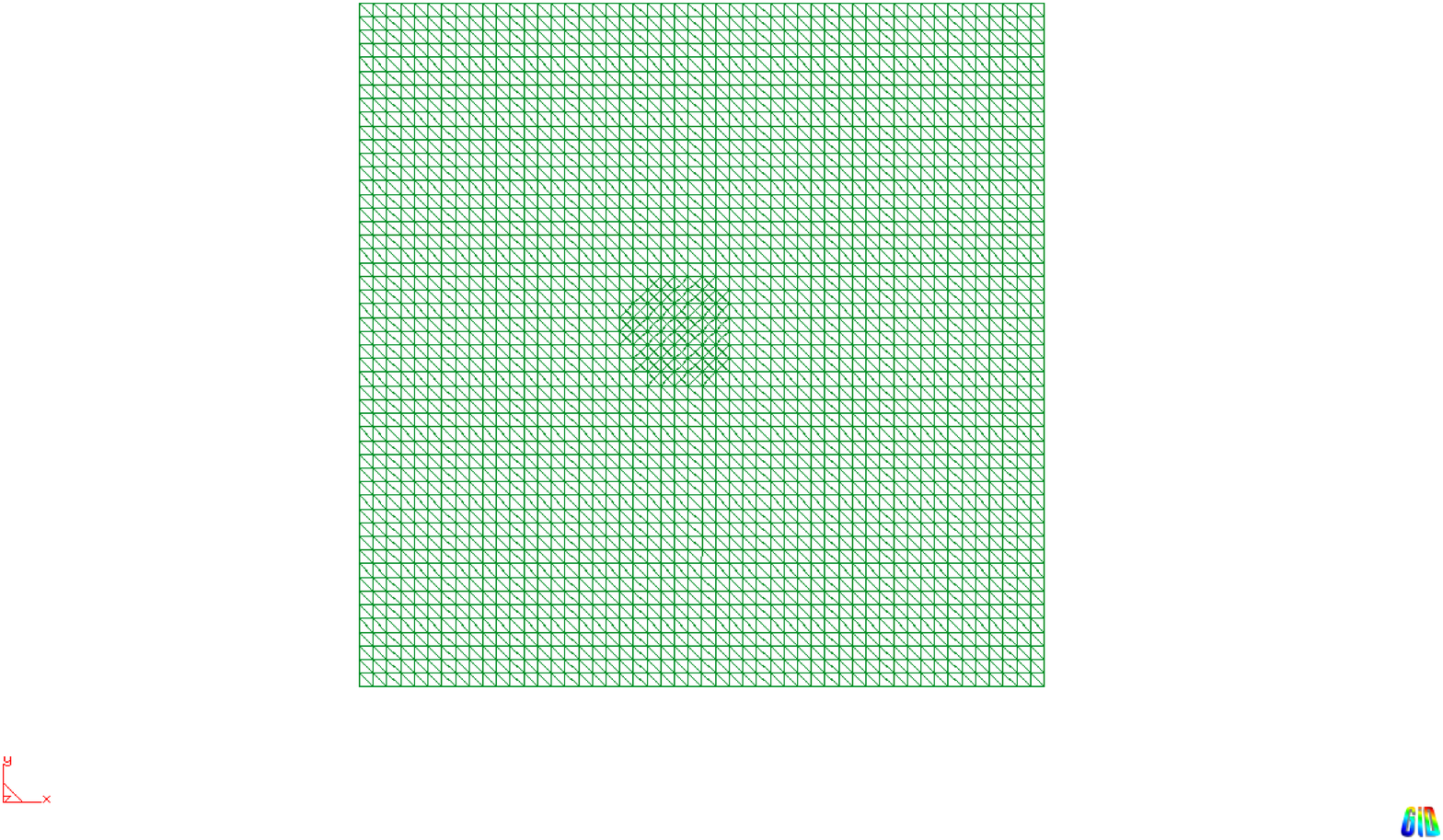}&
      \includegraphics[width=5cm, angle = 90, clip = true, trim = 11cm 12.5cm 11cm 12.5cm]{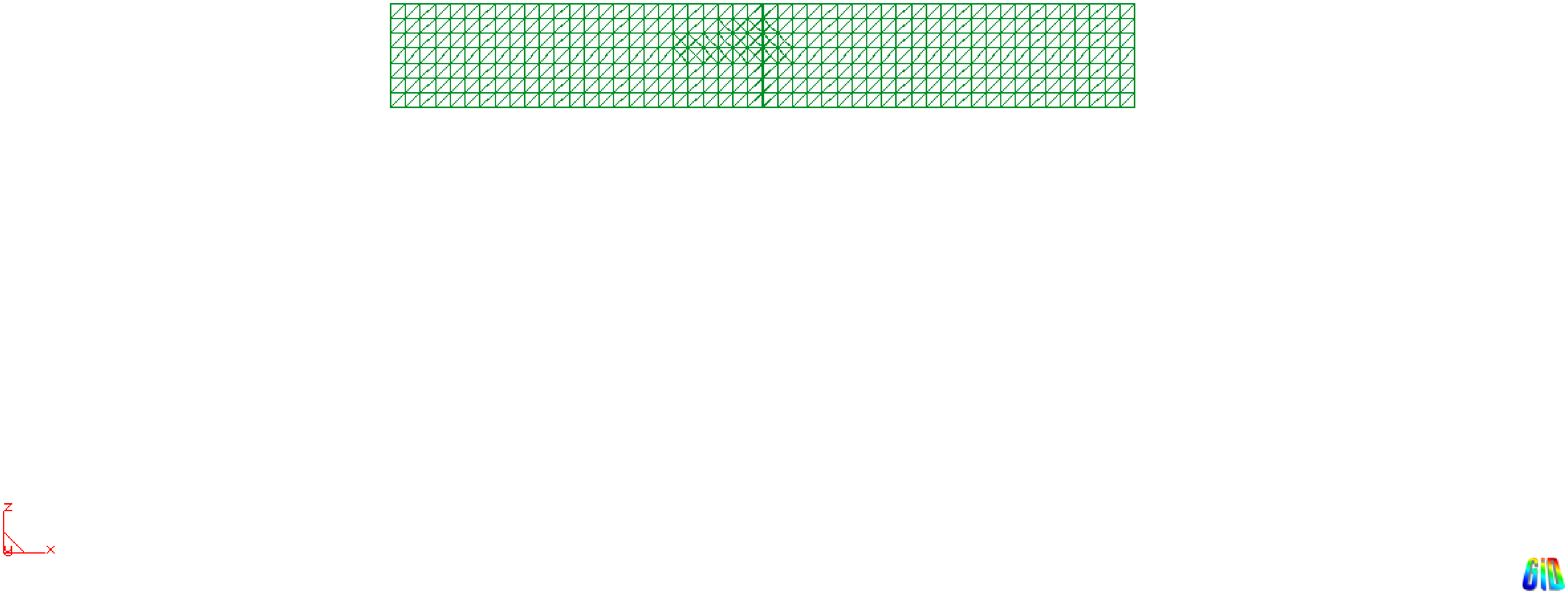}&
      \includegraphics[width=5cm, angle = 90, clip = true, trim = 11cm 12.5cm 11cm 12.5cm]{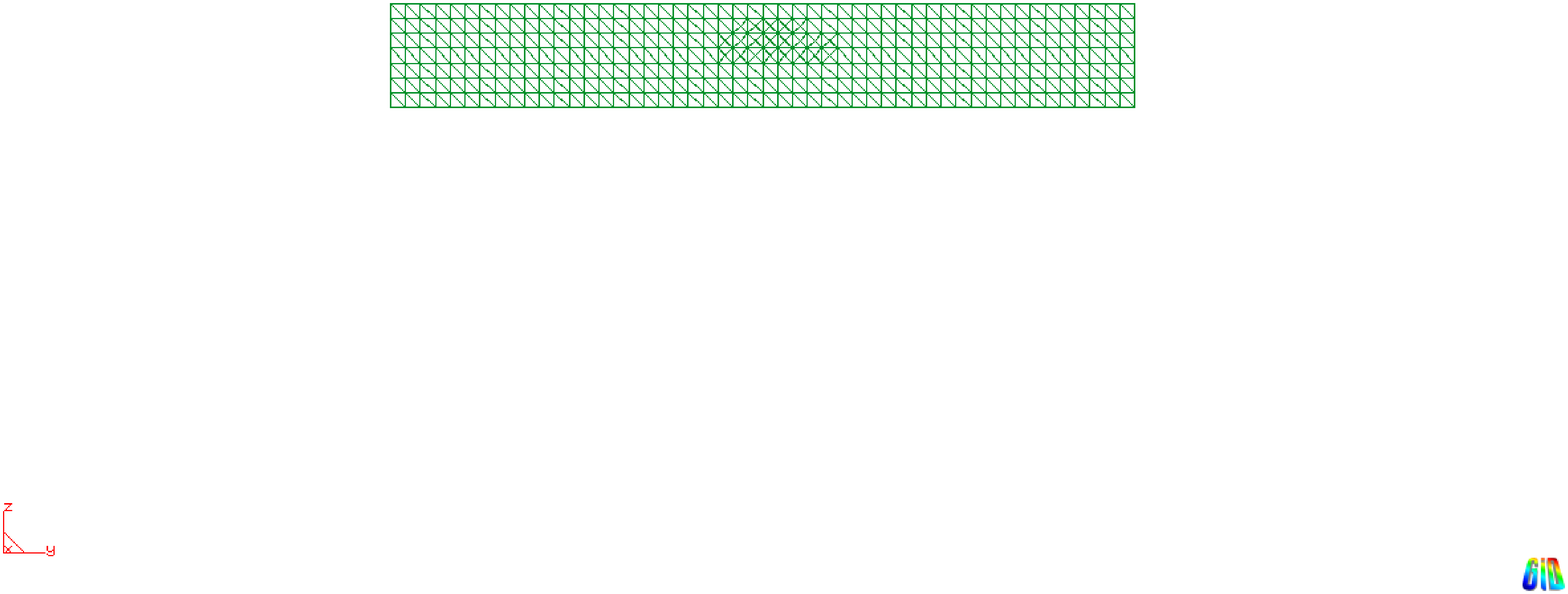}\\
      (a) $xy$-projection &
      (b) $xz$-projection &
      (c) $yz$-projection\\
      \includegraphics[width=5cm, clip=true, trim = 11cm 4cm 11cm 4cm]{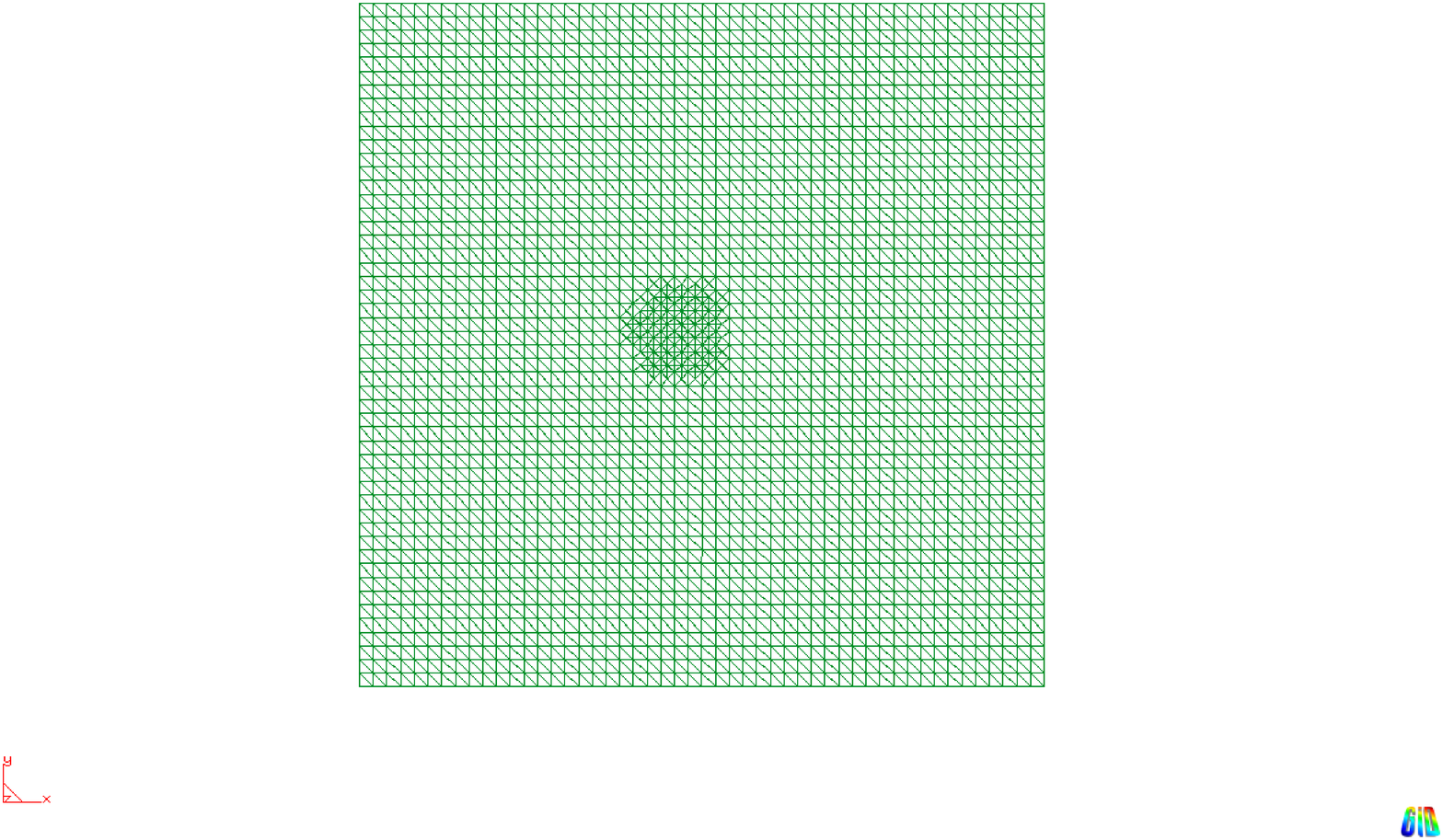}&
      \includegraphics[width=5cm, angle = 90, clip = true, trim = 11cm 12.5cm 11cm 12.5cm]{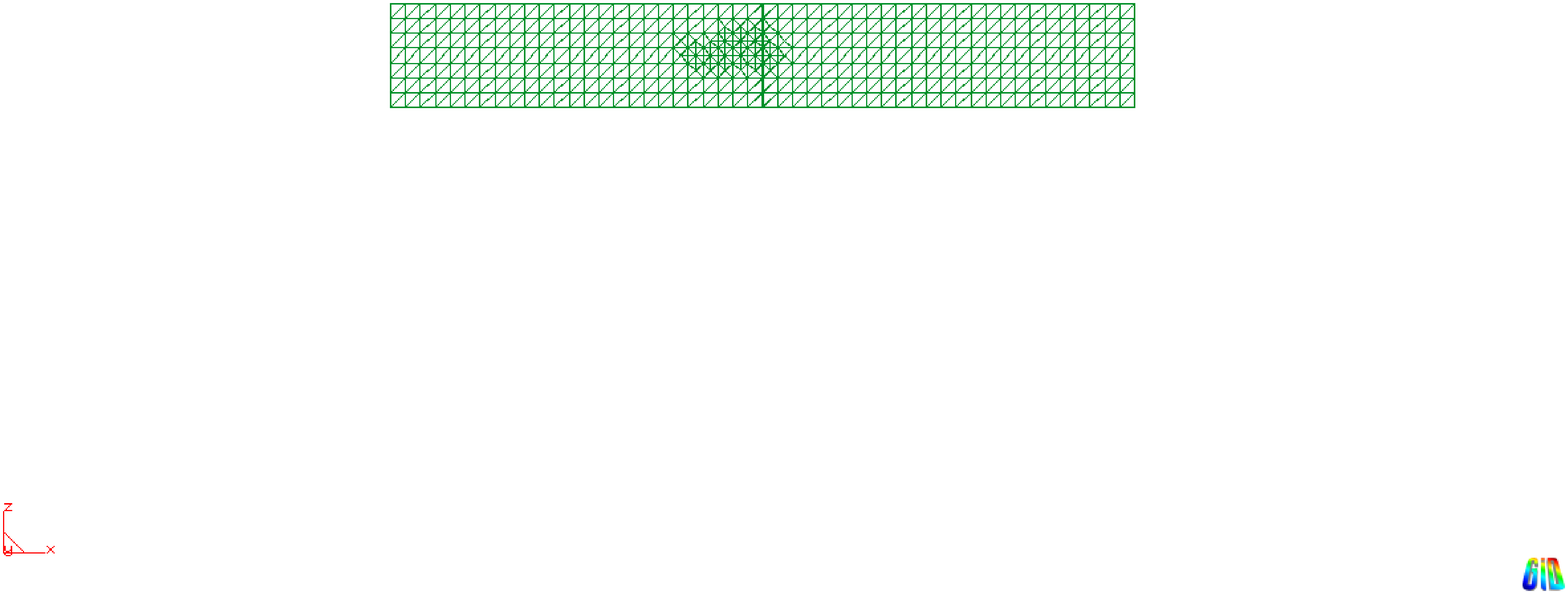}&
      \includegraphics[width=5cm, angle = 90, clip = true, trim = 11cm 12.5cm 11cm 12.5cm]{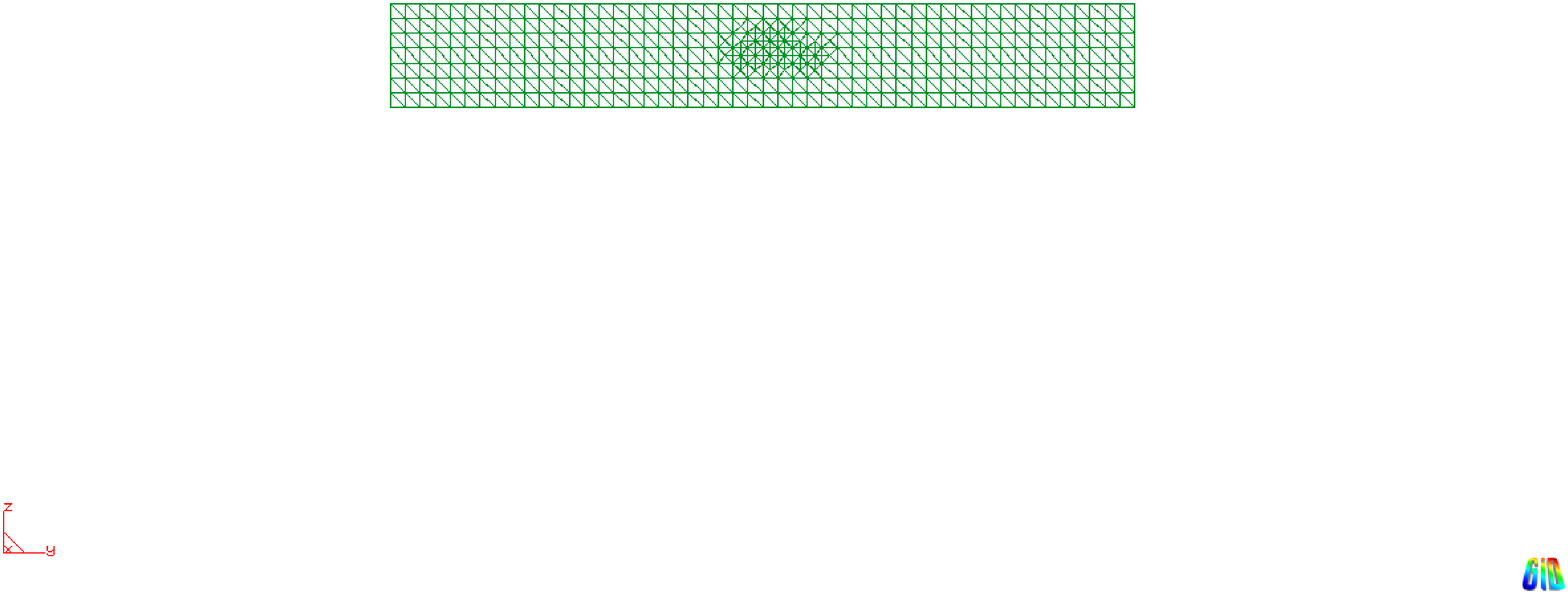}\\
      (d) $xy$-projection &
      (e) $xz$-projection &
      (f) $yz$-projection \\
      \includegraphics[width=5cm, clip=true, trim = 11cm 4cm 11cm 4cm]{4/ref3xy.eps}&
      \includegraphics[width=5cm, angle = 90, clip = true, trim = 11cm 12.5cm 11cm 12.5cm]{4/ref3xz.eps}&
      \includegraphics[width=5cm, angle = 90, clip = true, trim = 11cm 12.5cm 11cm 12.5cm]{4/ref3yz.eps}\\
      (g) $xy$-projection &
      (h) $xz$-projection &
      (i) $yz$-projection
    \end{tabular}
  \end{center}
    \caption{Adaptively refined meshes for the target
      number 3 of Table \ref{tab:table1}. (a) - (c) once refined, (d) - (f) twice refined, (g) - (i) three times refined mesh.}
\label{fig:fig7}
\end{figure}


\begin{figure}[tbp]
\begin{center}
\begin{tabular}{cc}
{\includegraphics[scale=0.18, angle=-90, clip=]{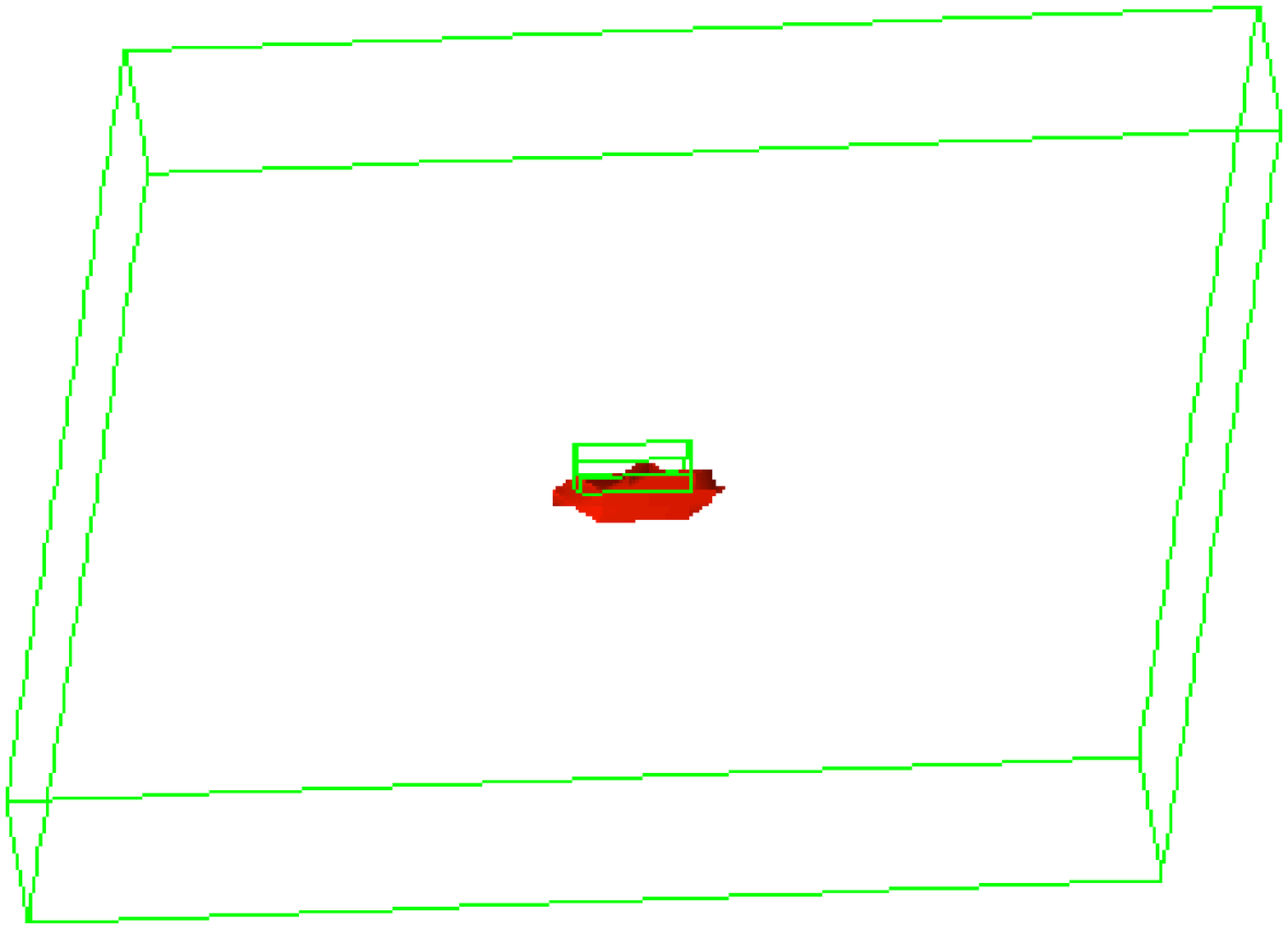}} & 
{\includegraphics[scale=0.18, angle=-90, clip=]{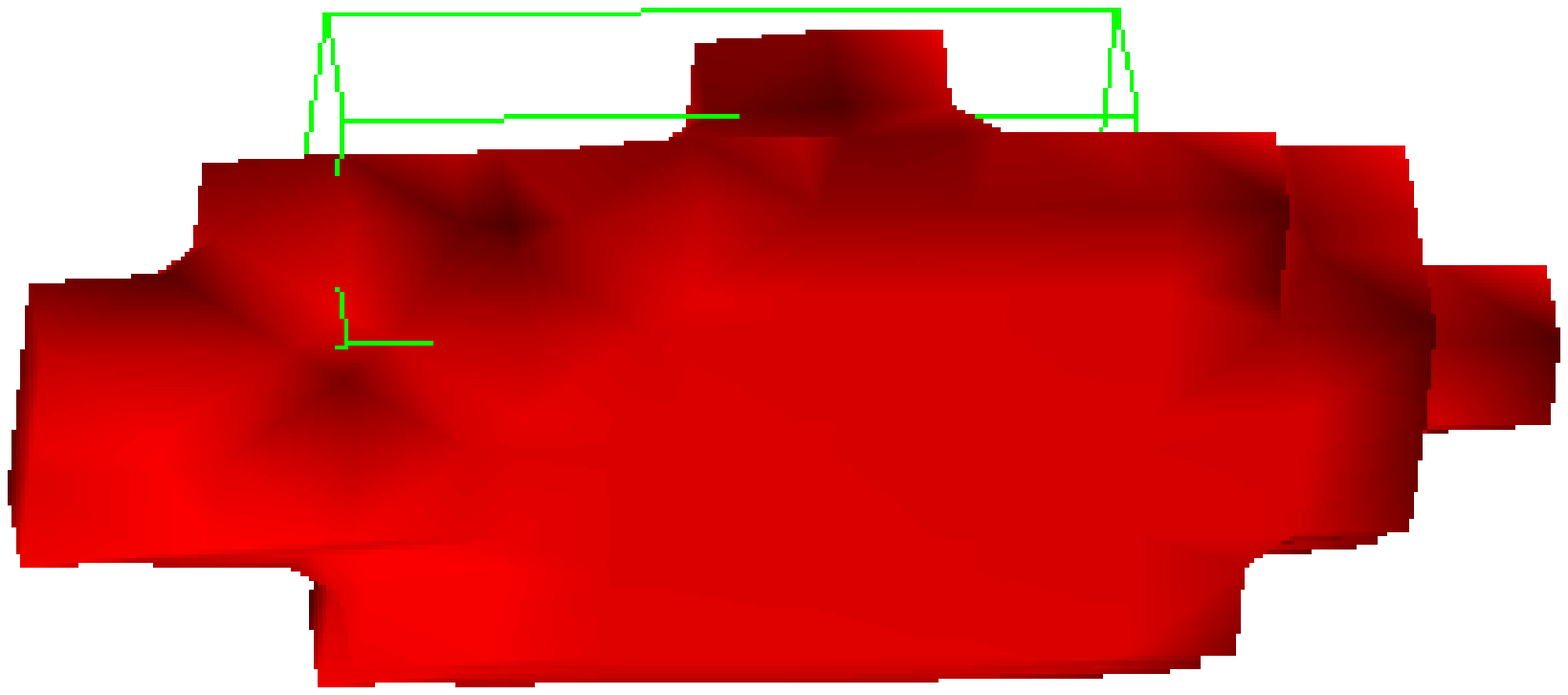}} \\
a)  coarse mesh  & b) zoom \\
{\includegraphics[scale=0.18, angle=-90, clip=]{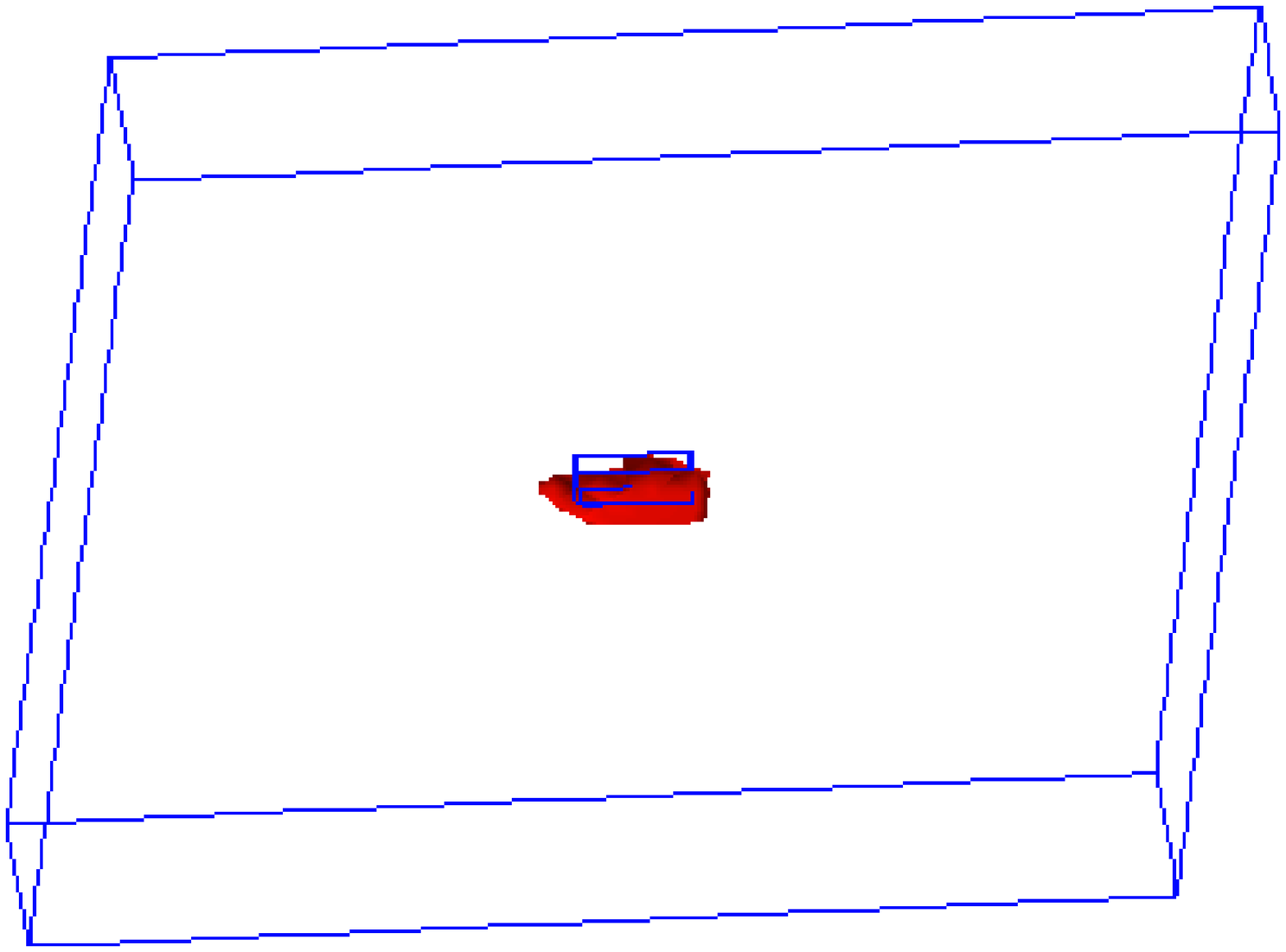}} & 
{\includegraphics[scale=0.18, angle=-90, clip=]{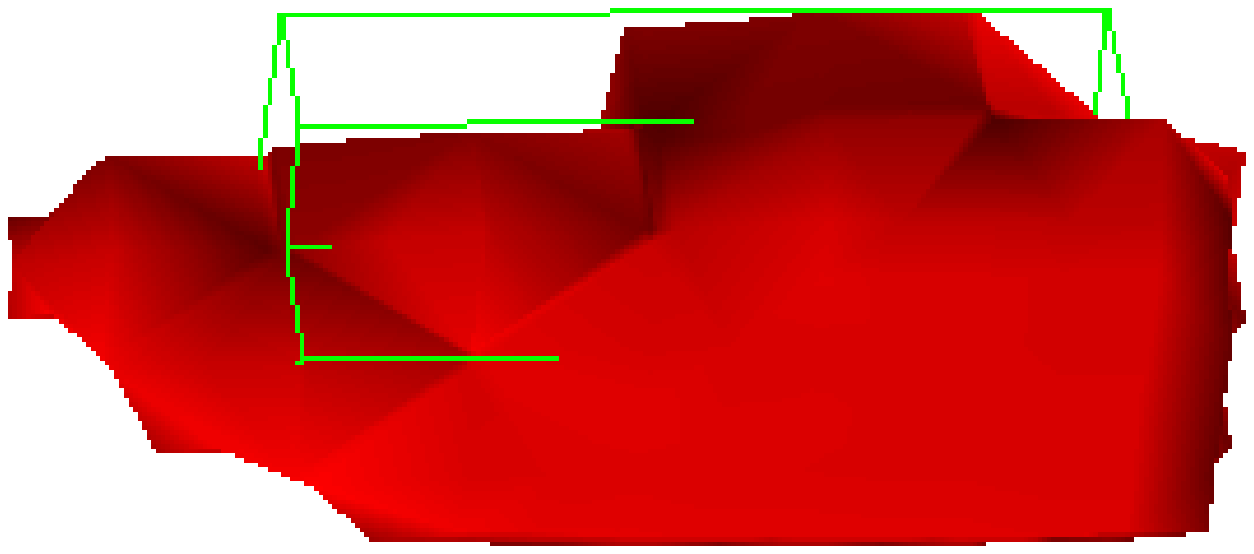}} \\
c)  once refined & d) zoom \\
{\includegraphics[scale=0.18, angle=-90, clip=]{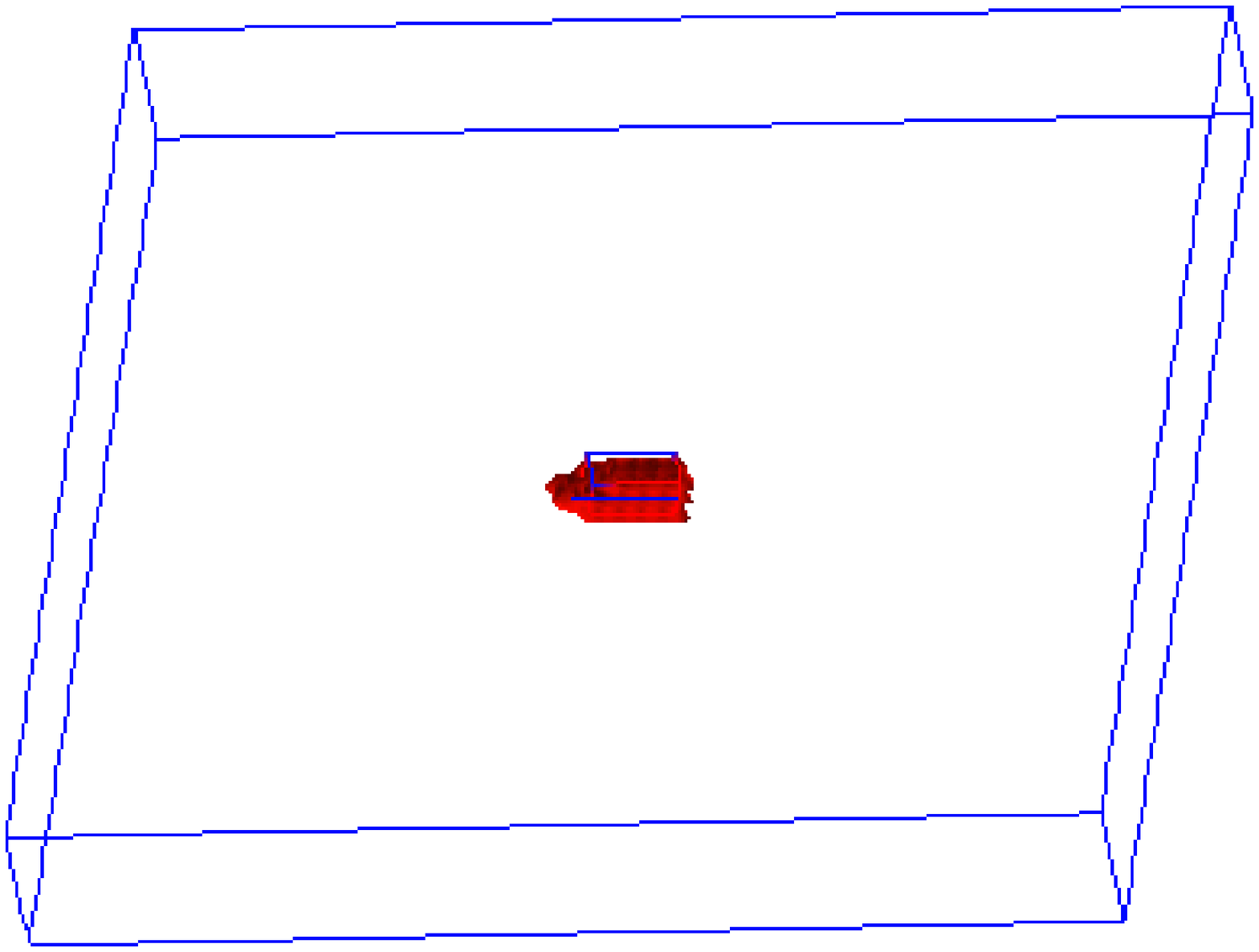}} & 
{\includegraphics[scale=0.18, angle=-90, clip=]{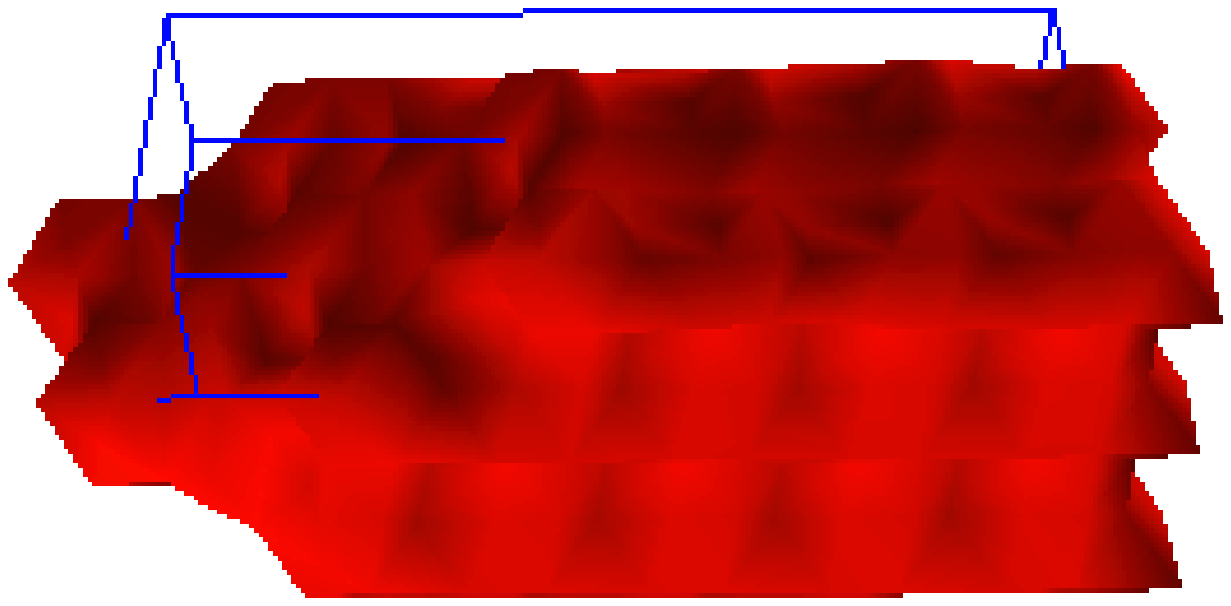}} \\
e) twice refined  & f) zoom \\
{\includegraphics[scale=0.18, angle=-90, clip=]{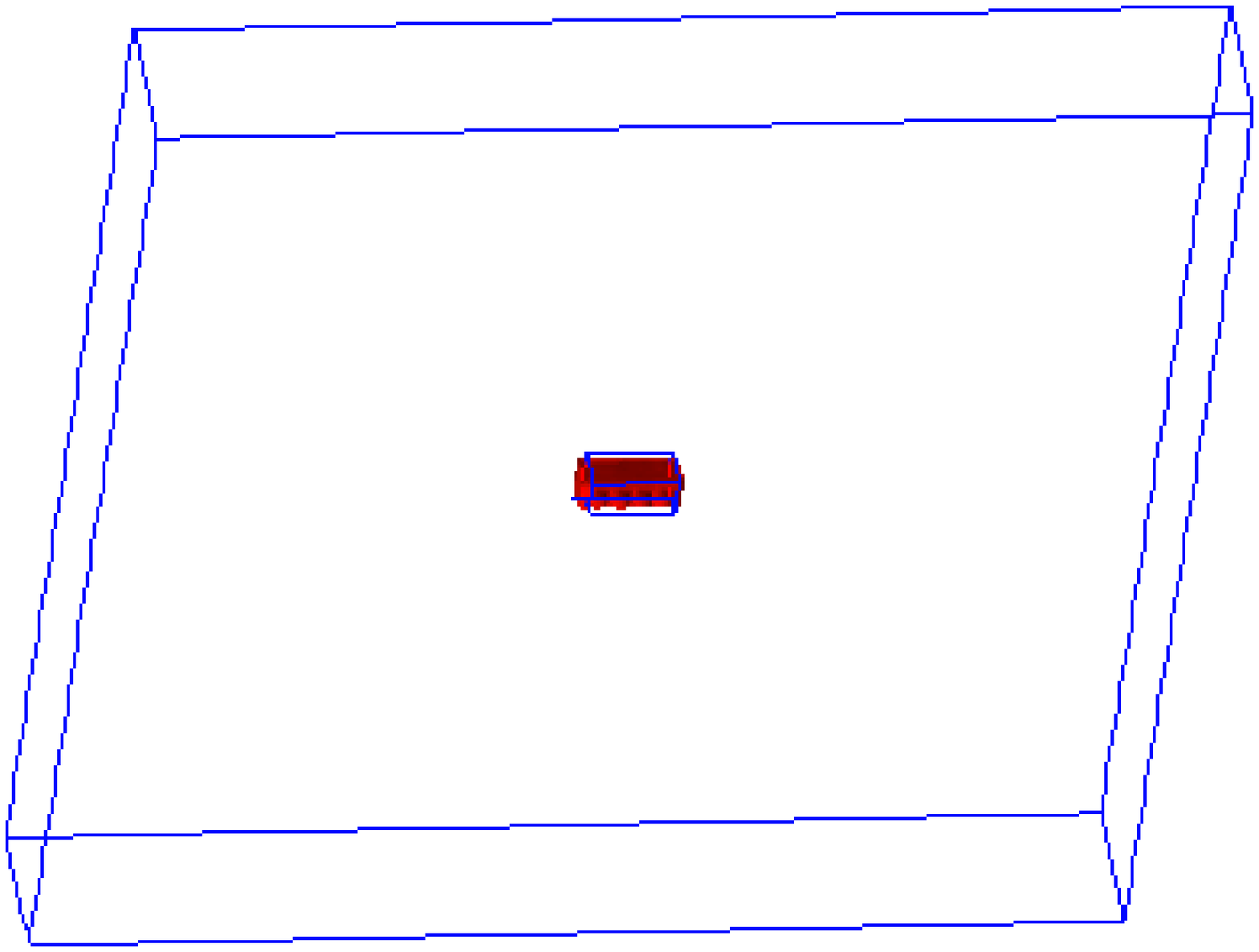}} & 
{\includegraphics[scale=0.18, angle=-90, clip=]{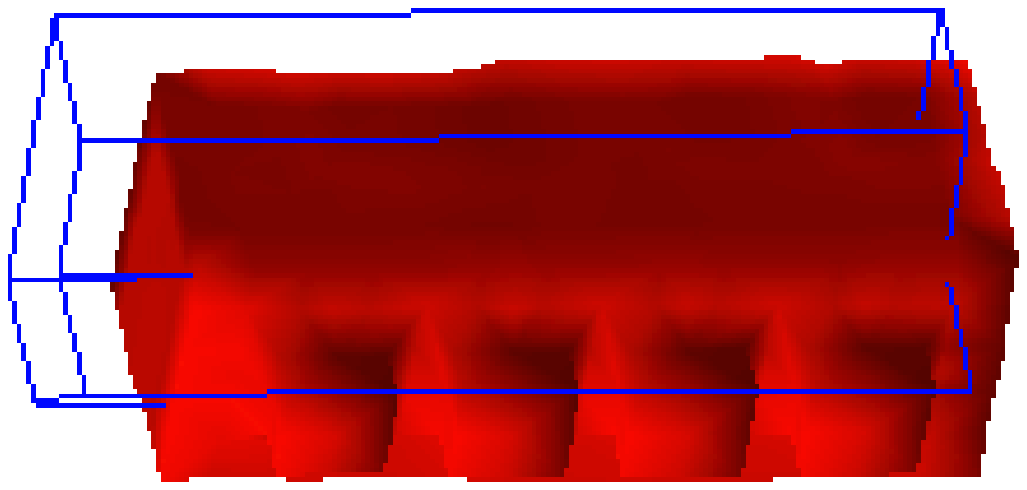}} \\
g)   three times refined  & h)  zoom  \\
{\includegraphics[scale=0.18, angle=-90, clip=]{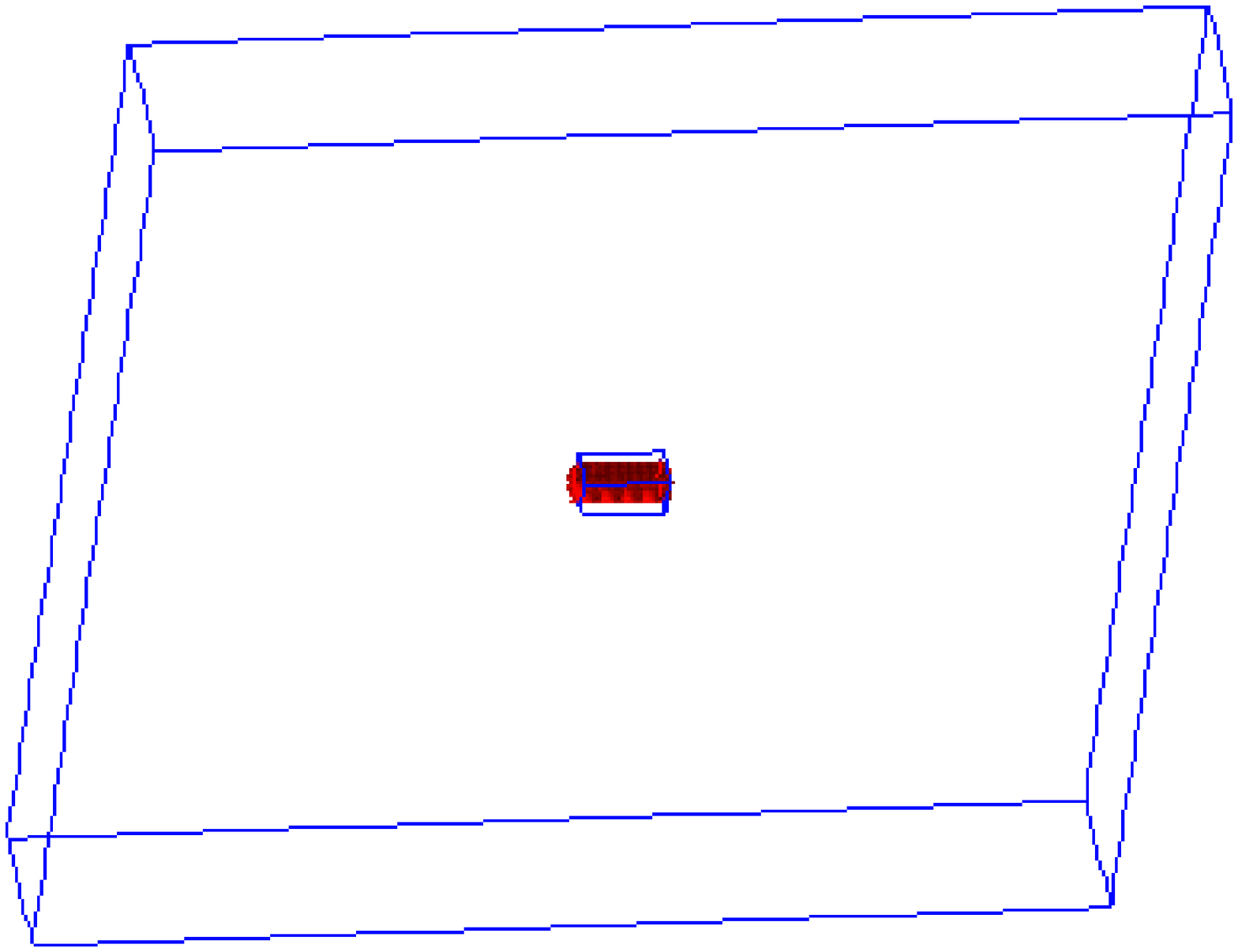}} & 
{\includegraphics[scale=0.18, angle=-90, clip=]{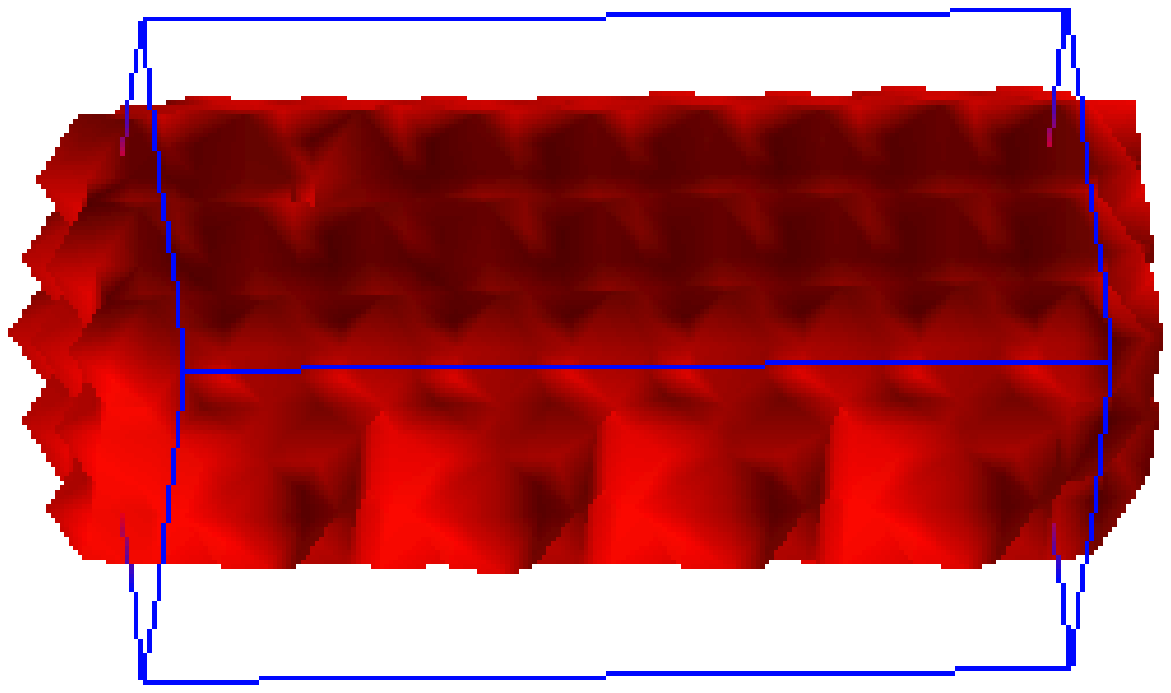}} \\
i)    four times refined & j) zoom 
\end{tabular}%
\end{center}
\caption{{\protect\small \emph{Computed images of target number 4
      of Table \ref{tab:table1} on four times adaptively refined
      meshes. Compare with Figure \ref{fig:globconv}-c).
 }}}
\label{fig:fig8}
\end{figure}

\clearpage
\begin{figure}[ht!]
  \begin{center}
    \begin{tabular}{cccc}
      \includegraphics[width=5cm, clip=true, trim = 11cm 4cm 11cm 4cm]{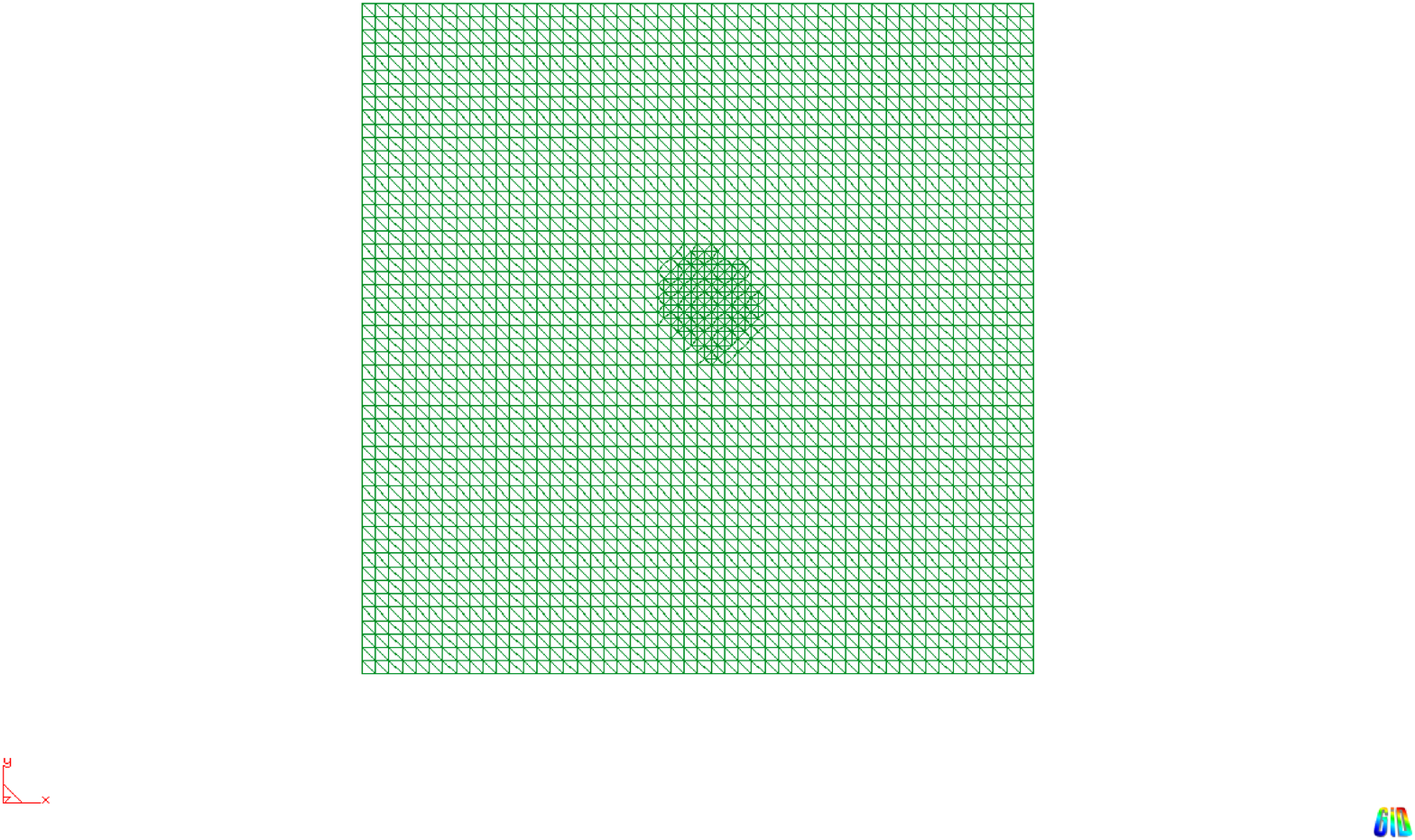}&
      \includegraphics[width=5cm, angle = 90, clip = true, trim = 11cm 12.5cm 11cm 12.5cm]{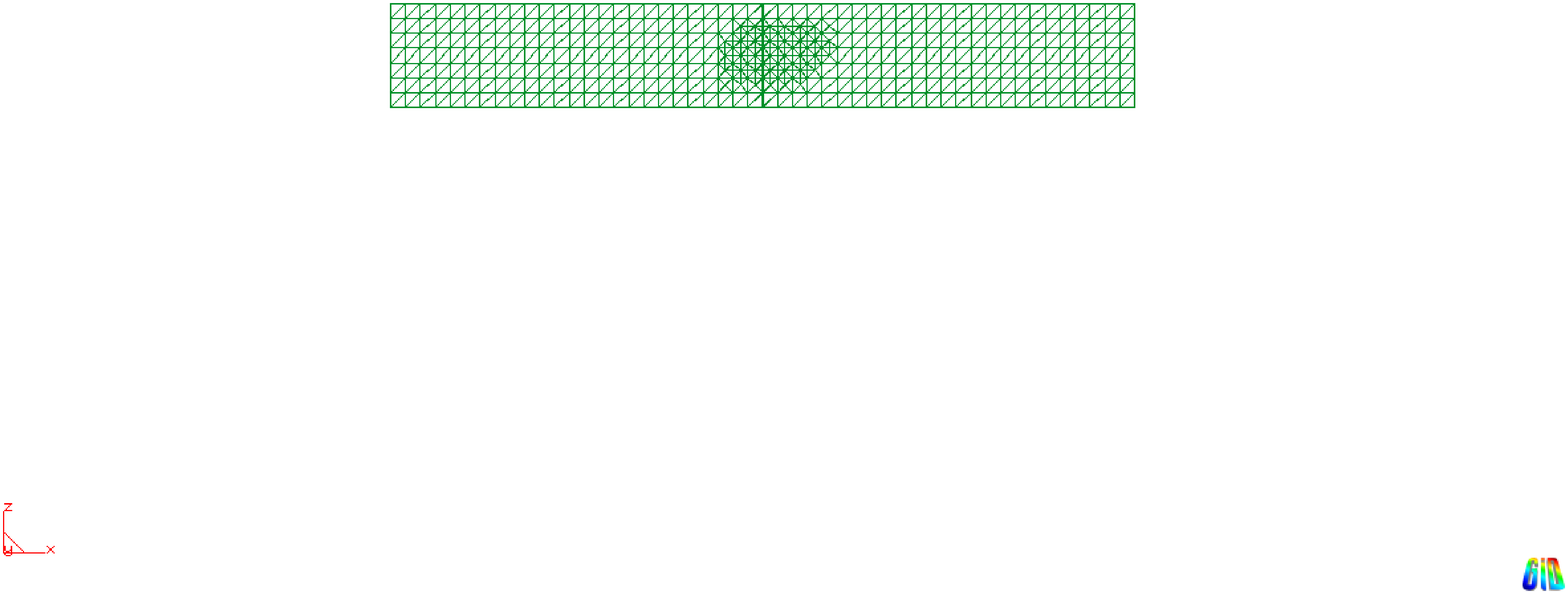}&
      \includegraphics[width=5cm, angle = 90, clip = true, trim = 11cm 12.5cm 11cm 12.5cm]{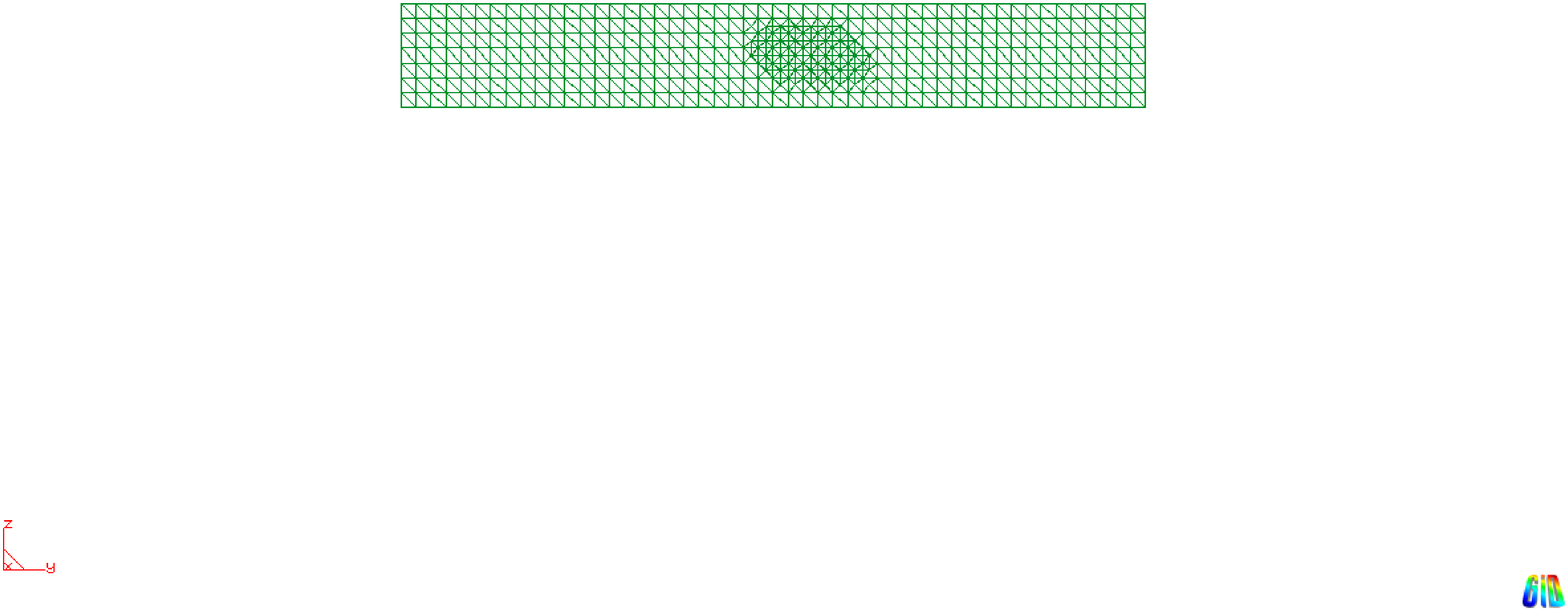}&
      \includegraphics[width=4.0cm, clip = true, trim = 1.6cm 0.0cm 4.0cm 0.0cm]{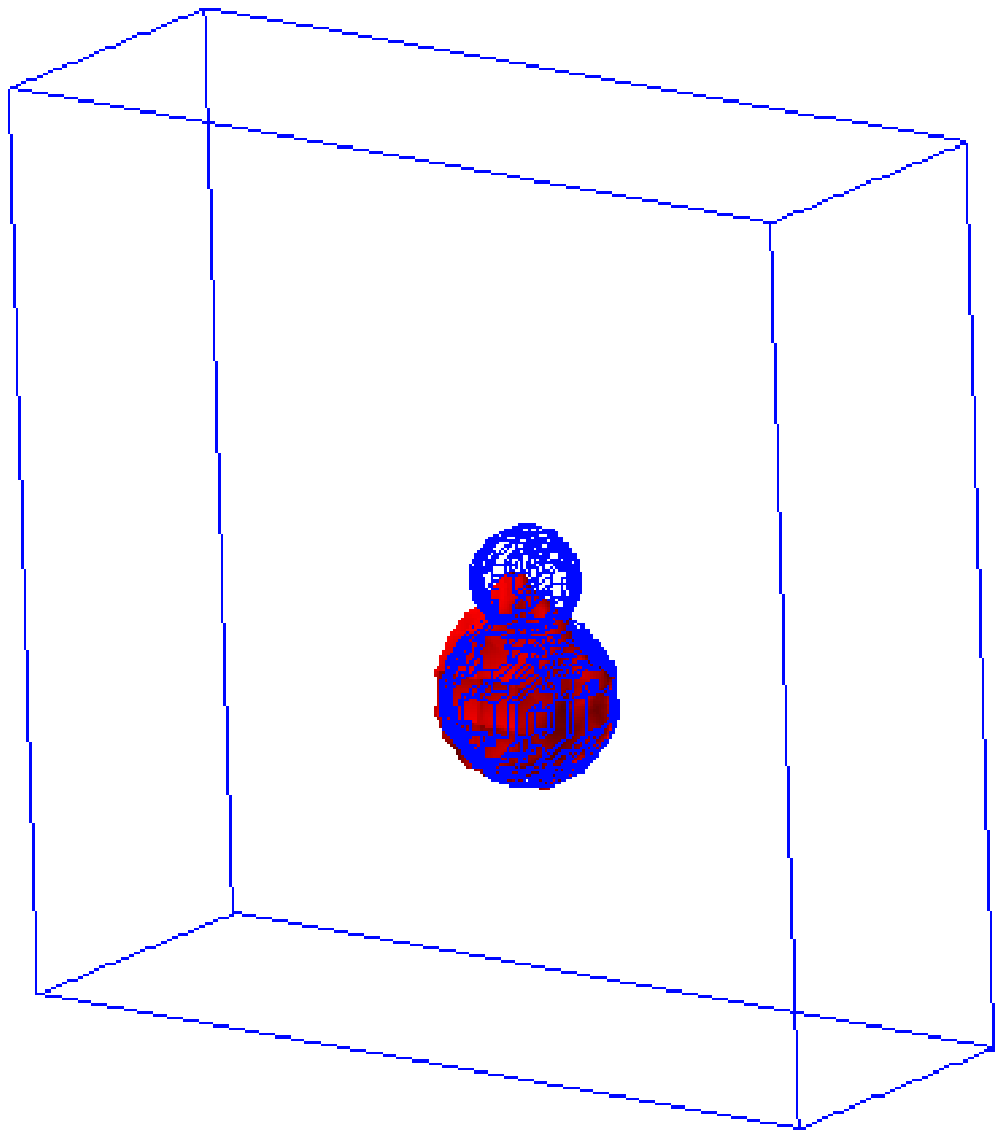}\\
      (a) &
      (b) &
      (c) &
      (d) 
    \end{tabular}
\end{center}
    \caption{(a) $xy$-projection, (b) $xz$-projection, and (c)
      $yz$-projection of the three times refined mesh; d) Computed
      image of target number 7 (doll, air inside) of Table \ref{tab:table1} on that
      mesh. Thin lines indicate correct shape.}
\label{fig:fig9}
\end{figure}

\begin{figure}[ht!]
  \begin{center}
    \begin{tabular}{cccc}
      \includegraphics[width=5cm, clip=true, trim = 11cm 4cm 11cm 4cm]{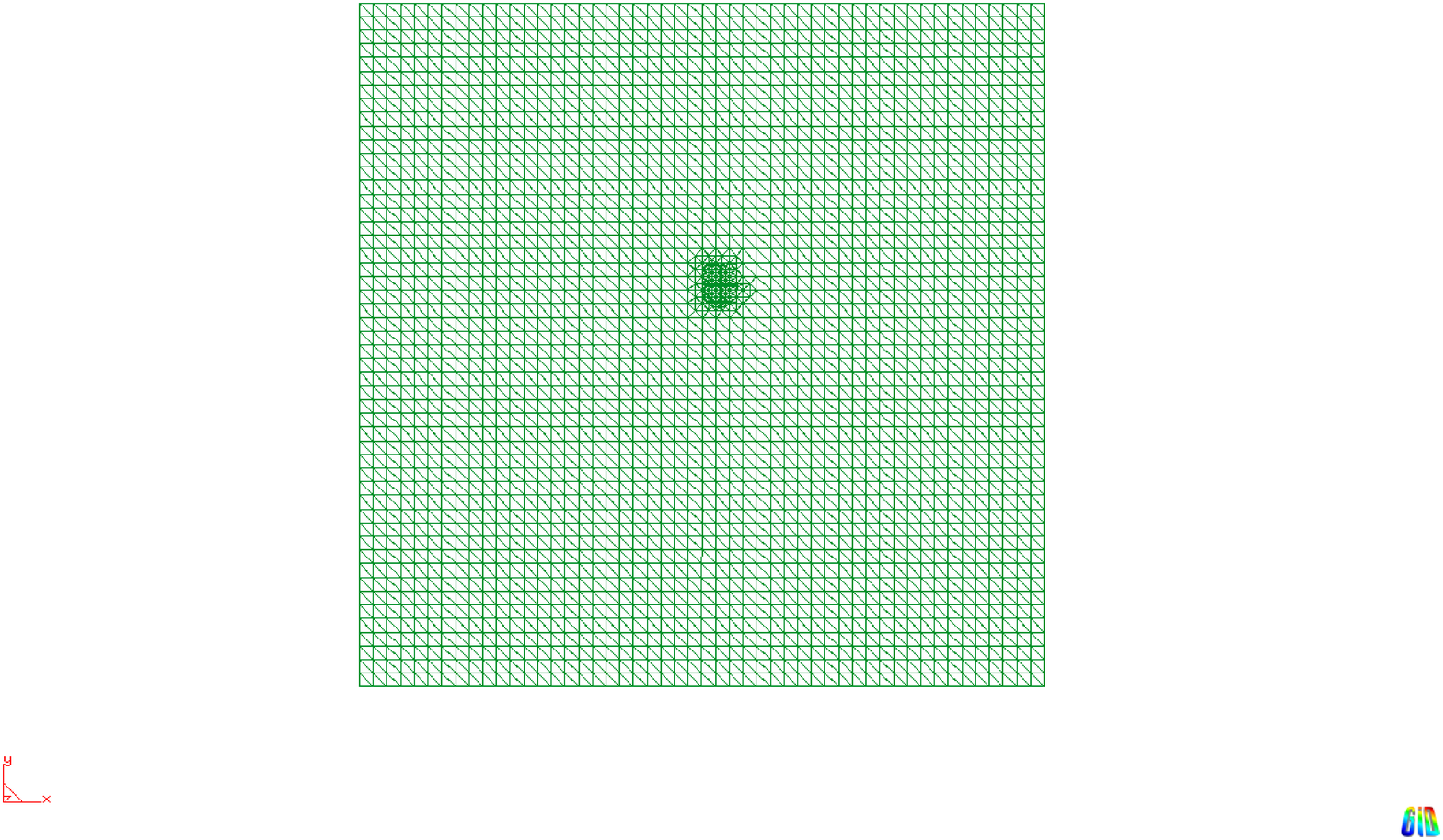}&
      \includegraphics[width=5cm, angle = 90, clip = true, trim = 11cm 12.5cm 11cm 12.5cm]{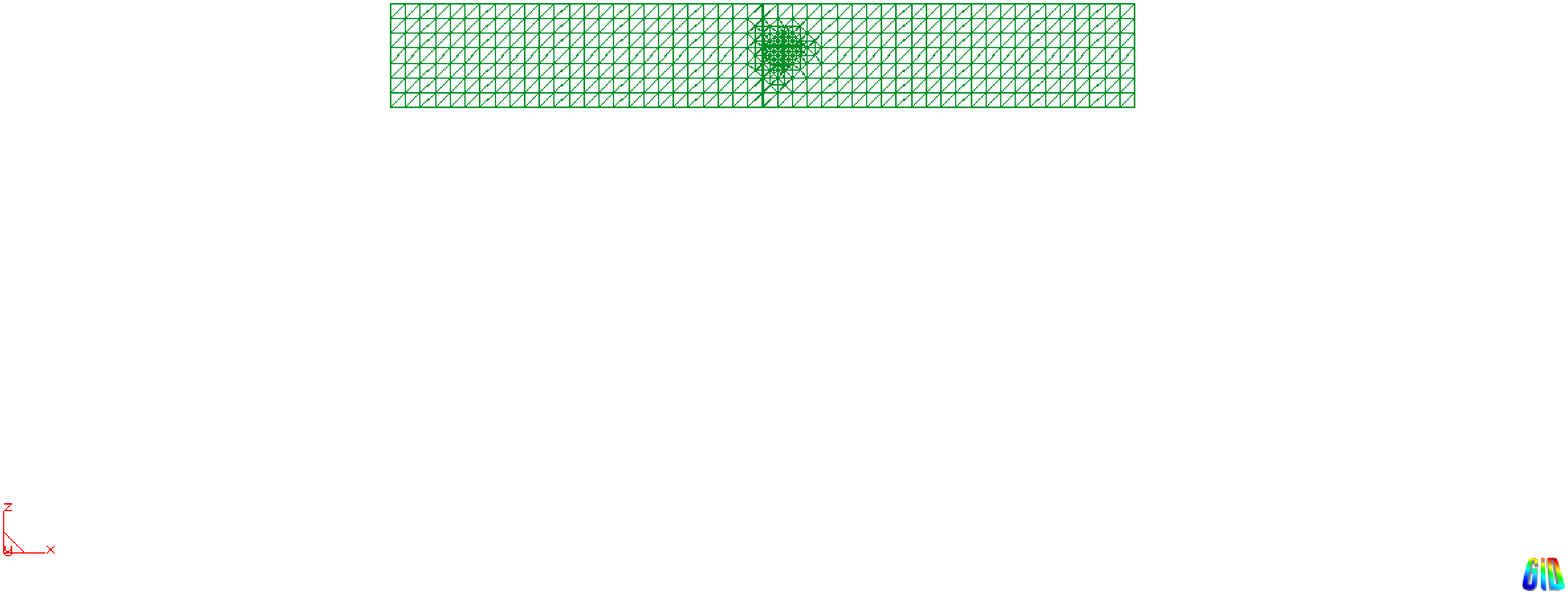}&
      \includegraphics[width=5cm, angle = 90, clip = true, trim = 11cm 12.5cm 11cm 12.5cm]{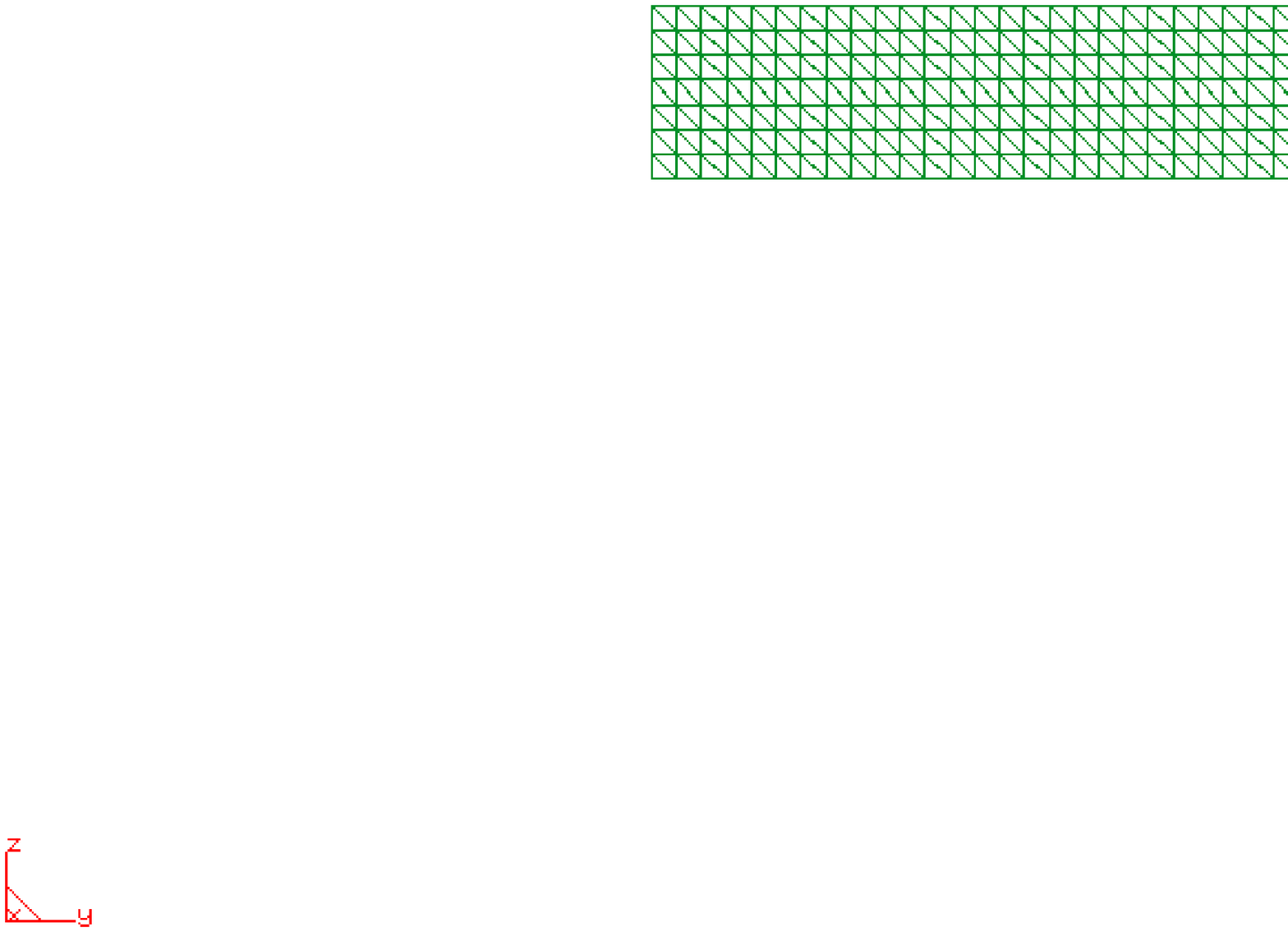}&
      \includegraphics[width=4.0cm, clip = true, trim = 1.8cm 0.0cm 4.0cm 0.0cm]{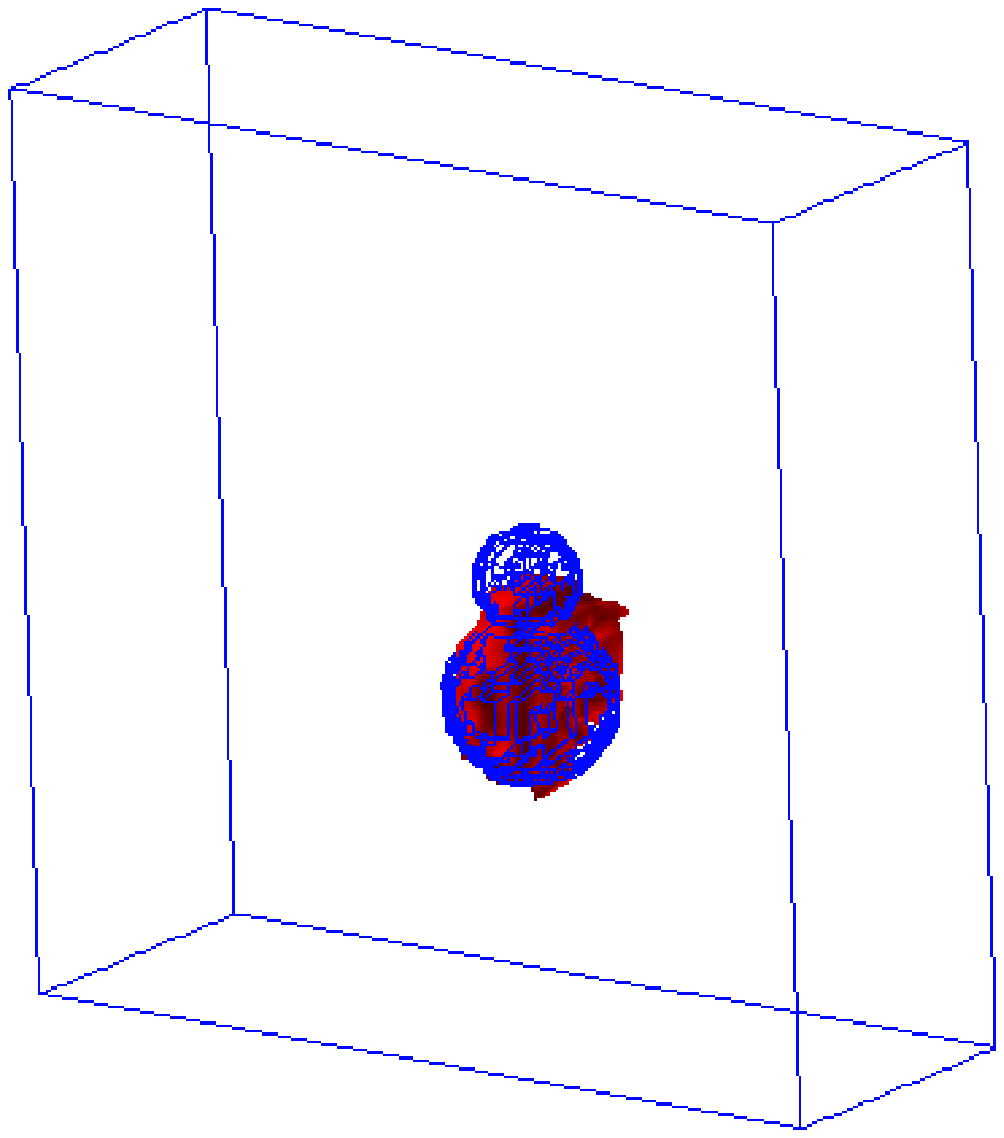}\\
      (a) &
      (b) &
      (c) &
      (d) 
    \end{tabular}
 \end{center}
    \caption{(a) $xy$-projection, (b) $xz$-projection, and (c) $yz$-projection of the three times refined mesh and the reconstruction (d) of target number  8  on that mesh.}
\label{fig:fig10}
\end{figure}

\begin{figure}
  \begin{center}
    \begin{tabular}{ccc}
      \includegraphics[width=4.0cm, clip = true, trim = 1.6cm 0.0cm 4.0cm 0.0cm]{16/eps_perspective.eps} &
      \includegraphics[width=4.0cm, clip = true, trim = 1.8cm 0.0cm 4.0cm 0.0cm]{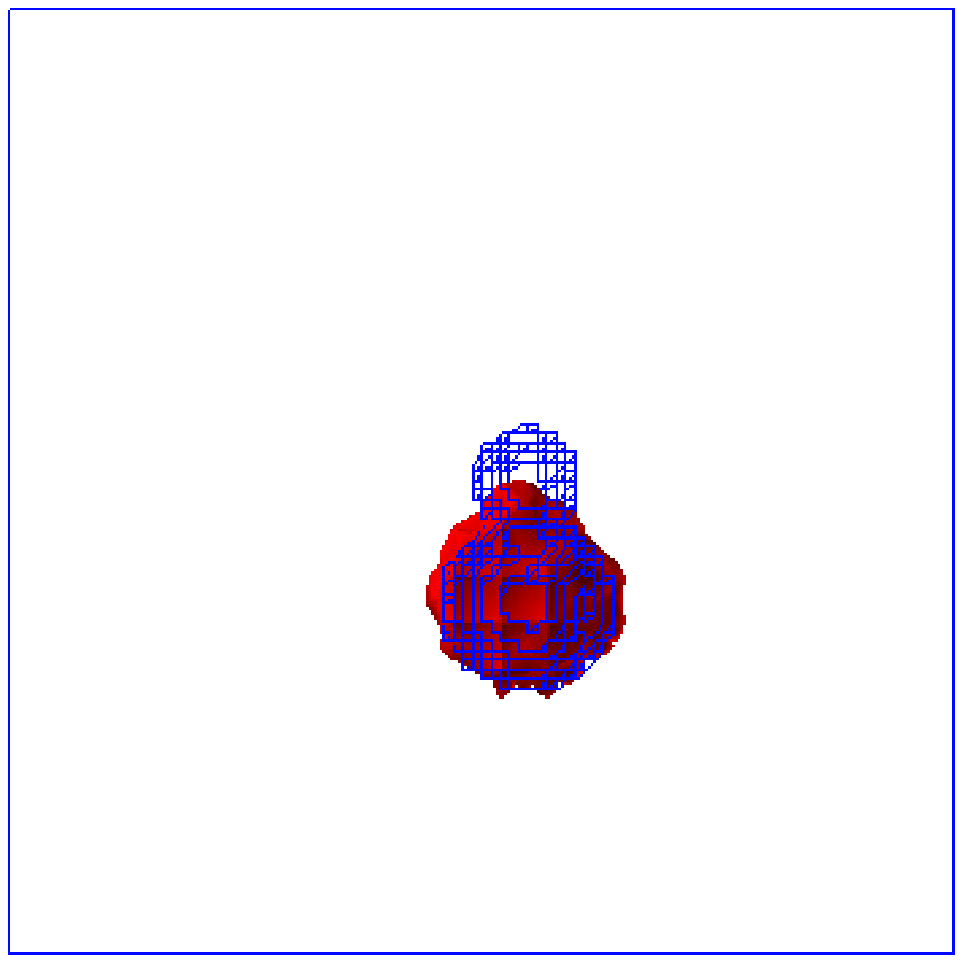} &
      \includegraphics[width=4.0cm, clip = true, trim = 1.8cm 0.0cm 4.0cm 0.0cm]{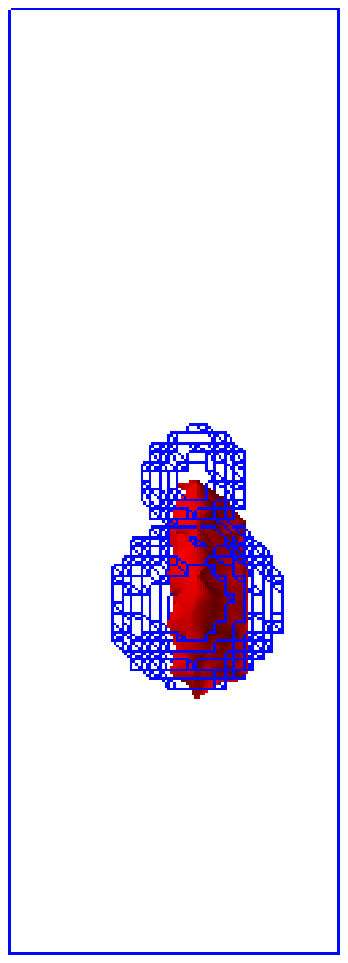}\\
      (a) Perspective view&
      (b) Front view&
      (c) Side view\\
      \includegraphics[width=4.0cm, clip = true, trim = 2.0cm 0.0cm 4.0cm 0.0cm]{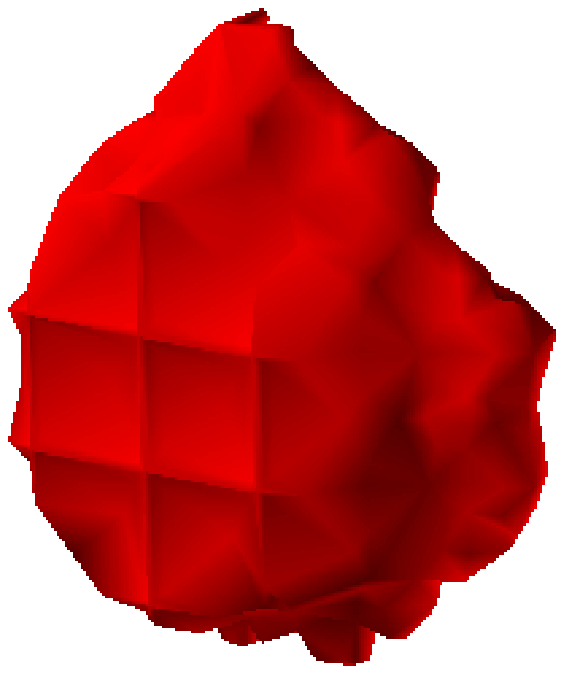} &
      \includegraphics[width=4.0cm, clip = true, trim = 2.0cm 0.0cm 4.0cm 0.0cm]{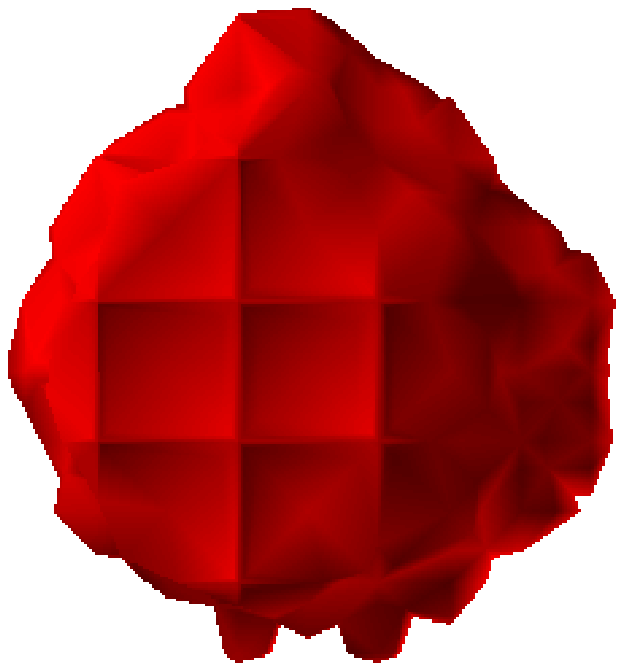} &
      \includegraphics[width=4.0cm, clip = true, trim = 2.0cm 0.0cm 4.0cm 0.0cm]{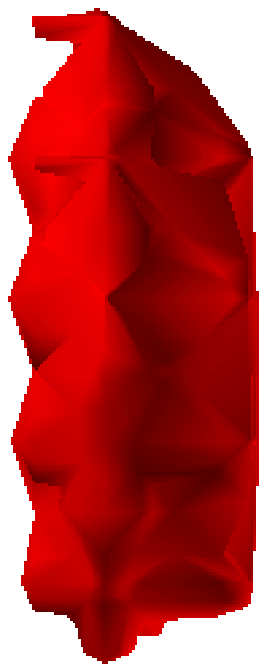}\\
      (d) Zoom, perspective&
      (e) Zoom, front&
      (f) Zoom, side \\
  \includegraphics[width=4.0cm, clip = true, trim = 1.8cm 0.0cm 4.0cm 0.0cm]{17/eps_perspective.eps} &
      \includegraphics[width=4.0cm, clip = true, trim = 1.8cm 0.0cm 4.0cm 0.0cm]{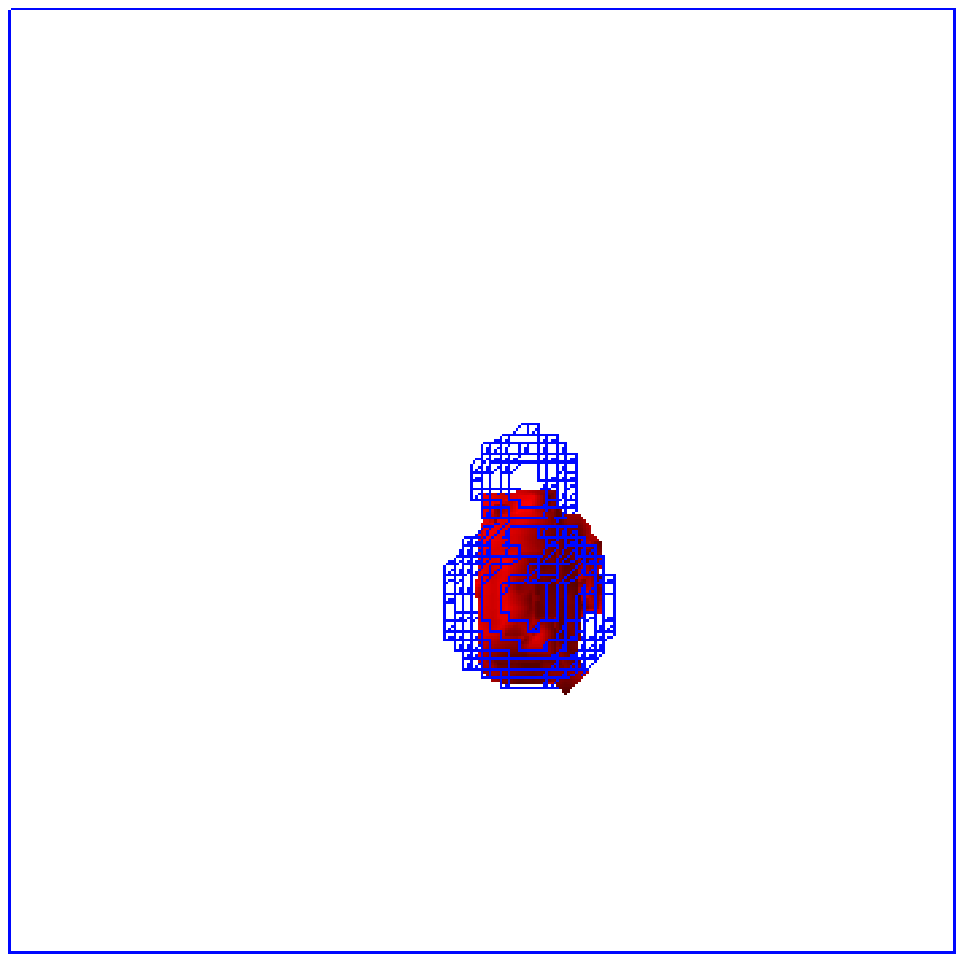} &
      \includegraphics[width=4.0cm, clip = true, trim = 1.8cm 0.0cm 4.0cm 0.0cm]{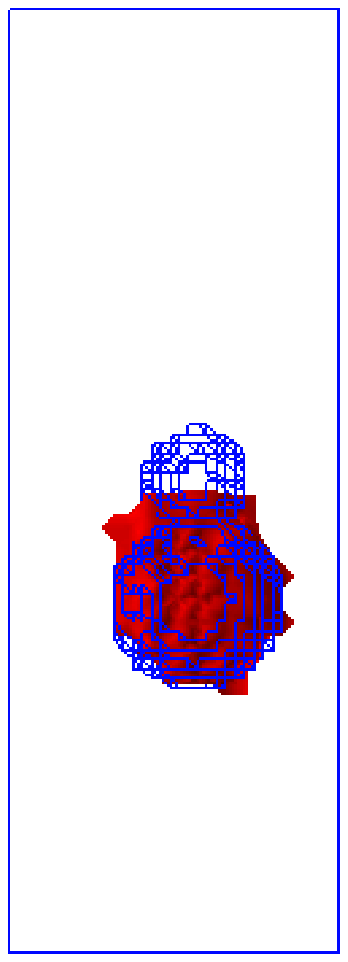}\\
      (g) Perspective view&
      (h) Front view&
      (i) Side view\\
      \includegraphics[width=4.0cm, clip = true, trim = 2.0cm 0.0cm 4.0cm 0.0cm]{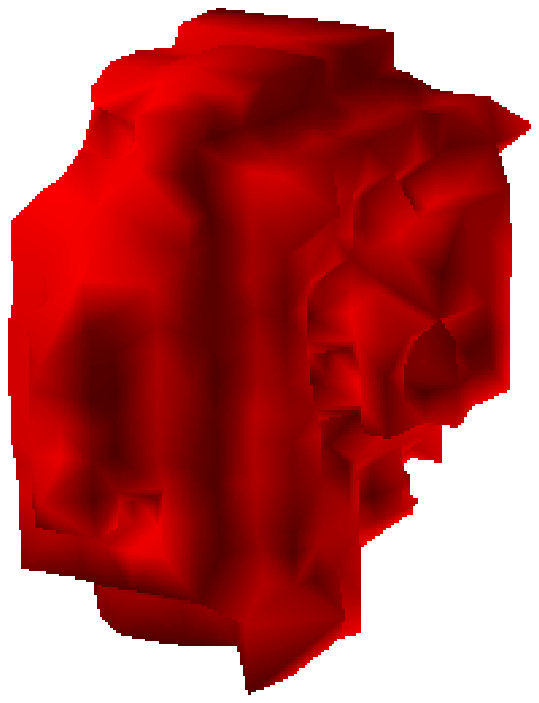} &
      \includegraphics[width=4.0cm, clip = true, trim = 2.0cm 0.0cm 4.0cm 0.0cm]{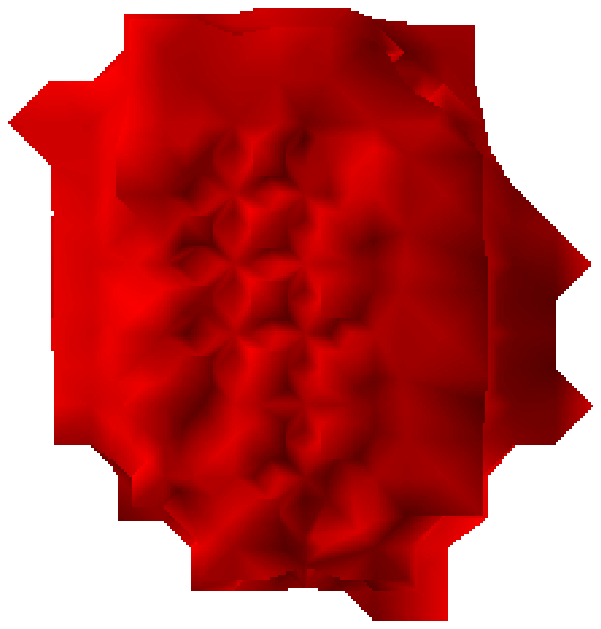} &
      \includegraphics[width=4.0cm, clip = true, trim = 2.0cm 0.0cm 4.0cm 0.0cm]{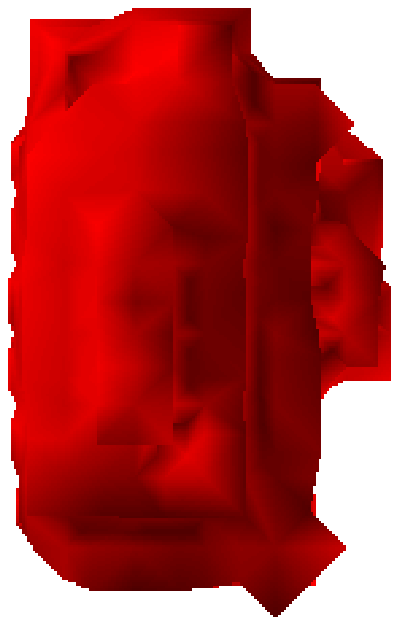}\\
      (j) Zoom, perspective&
      (k) Zoom, front&
      (l) Zoom, side
    \end{tabular}
  \end{center}
    \caption{Three views and zooms of targets number 7 (figures a)-f),
      doll, air inside) and number 8 (figures g)-l), doll, metal
      inside) of Table \ref{tab:table1} on three times refined
      mesh. Thin lines indicate correct shape. We observe that
      on d) and e) even the head of the doll is indicated, which was
      the most difficult part. As to images g)-l), it follows from
      (\ref{8.7}) that they display mostly the metallic part. Comparison of d), j)
      with Figures \ref{fig:globconv}-d), e) again shows a significant
      improvement of the image due to the adaptivity. }
\label{fig:fig11}
\end{figure}

\begin{figure}
  \begin{center} 
    \begin{tabular}{ccc}
      \includegraphics[width=5cm, clip=true, trim = 11cm 4cm 11cm 4cm]{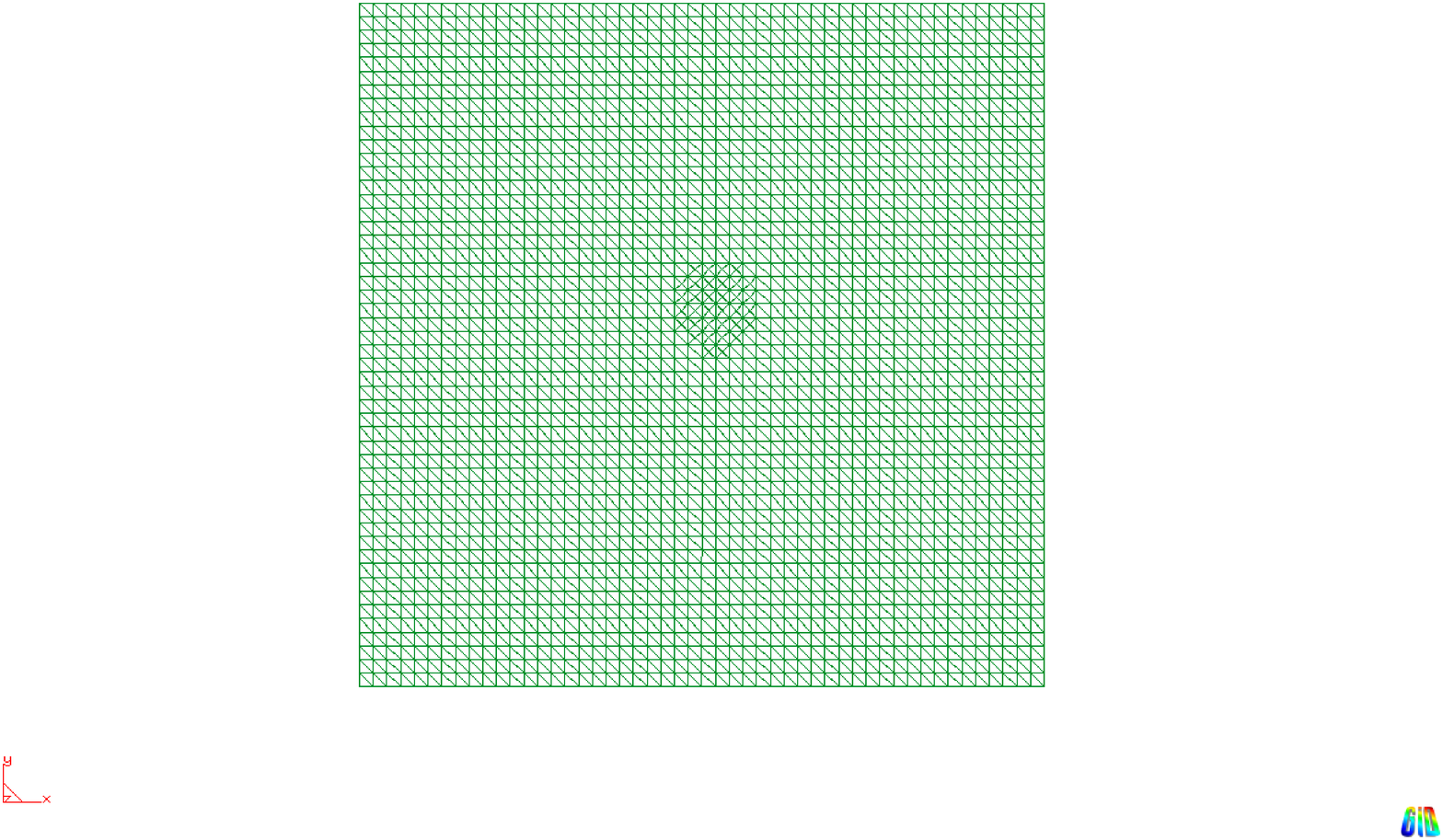}&
      \includegraphics[width=5cm, angle = 90, clip = true, trim = 11cm 12.5cm 11cm 12.5cm]{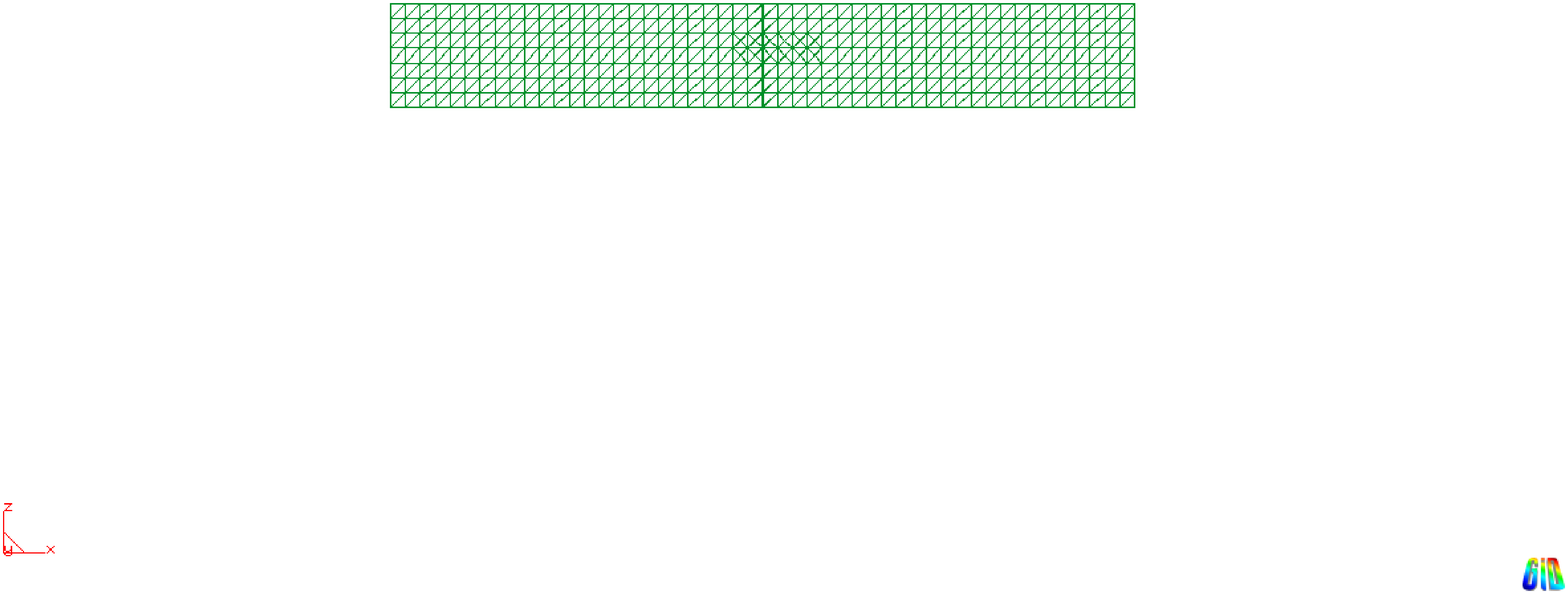}&
      \includegraphics[width=5cm, angle = 90, clip = true, trim = 11cm 12.5cm 11cm 12.5cm]{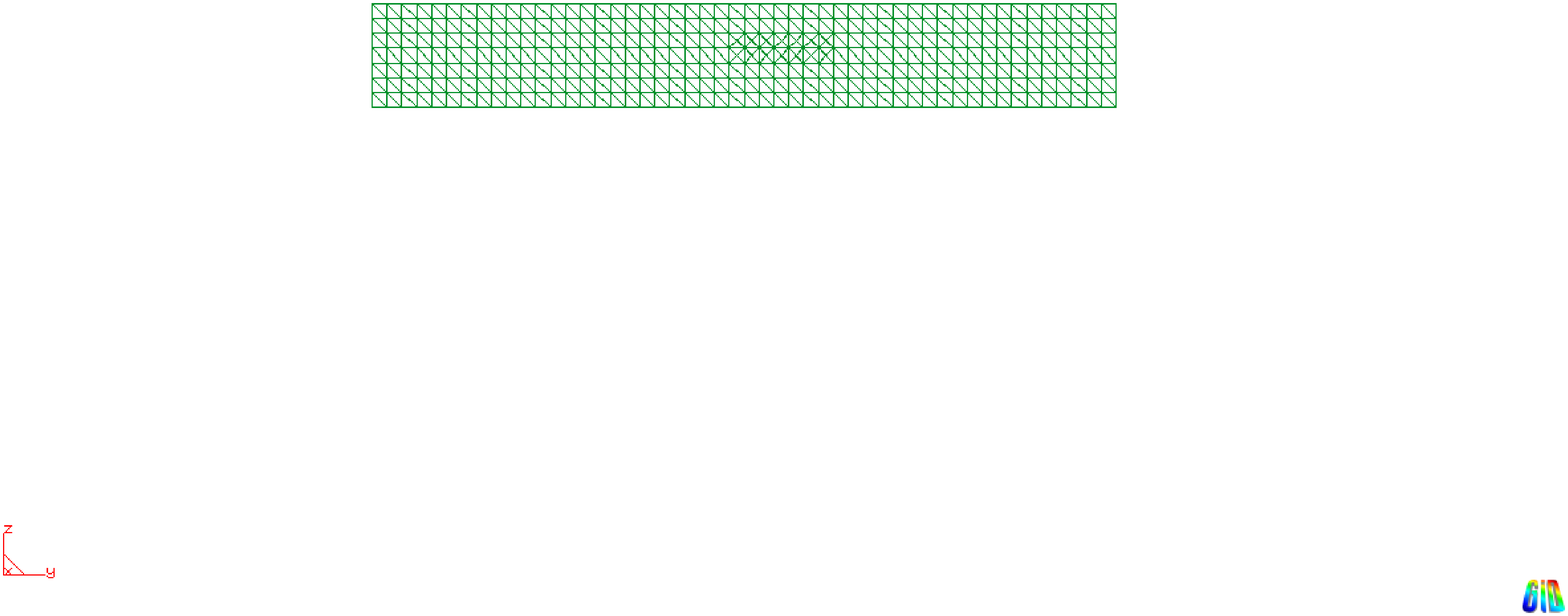}\\
      (a) $xy$-projection &
      (b) $xz$-projection &
      (c) $yz$-projection\\
      \includegraphics[width=5cm, clip=true, trim = 11cm 4cm 11cm 4cm]{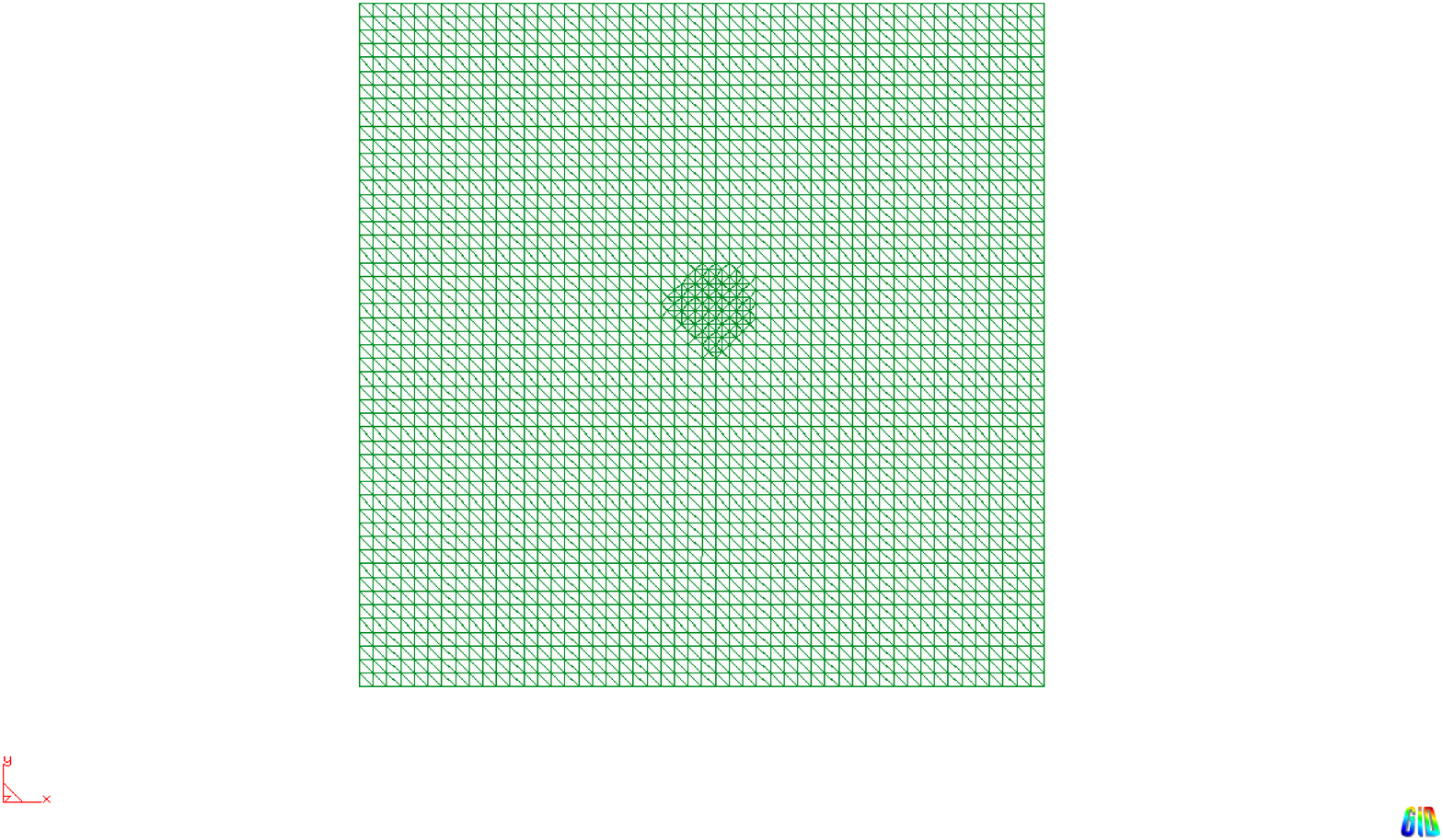}&
      \includegraphics[width=5cm, angle = 90, clip = true, trim = 11cm 12.5cm 11cm 12.5cm]{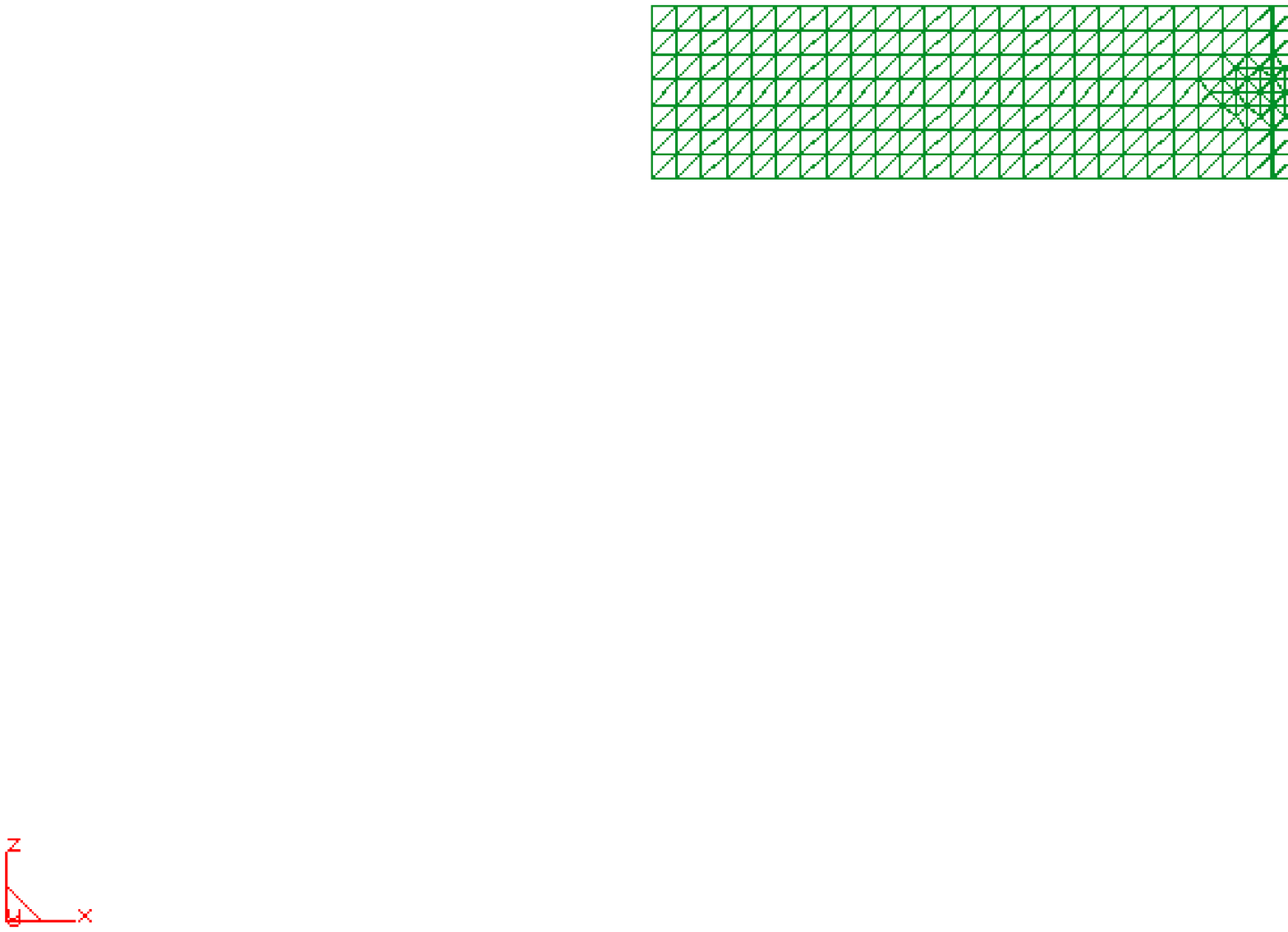}&
      \includegraphics[width=5cm, angle = 90, clip = true, trim = 11cm 12.5cm 11cm 12.5cm]{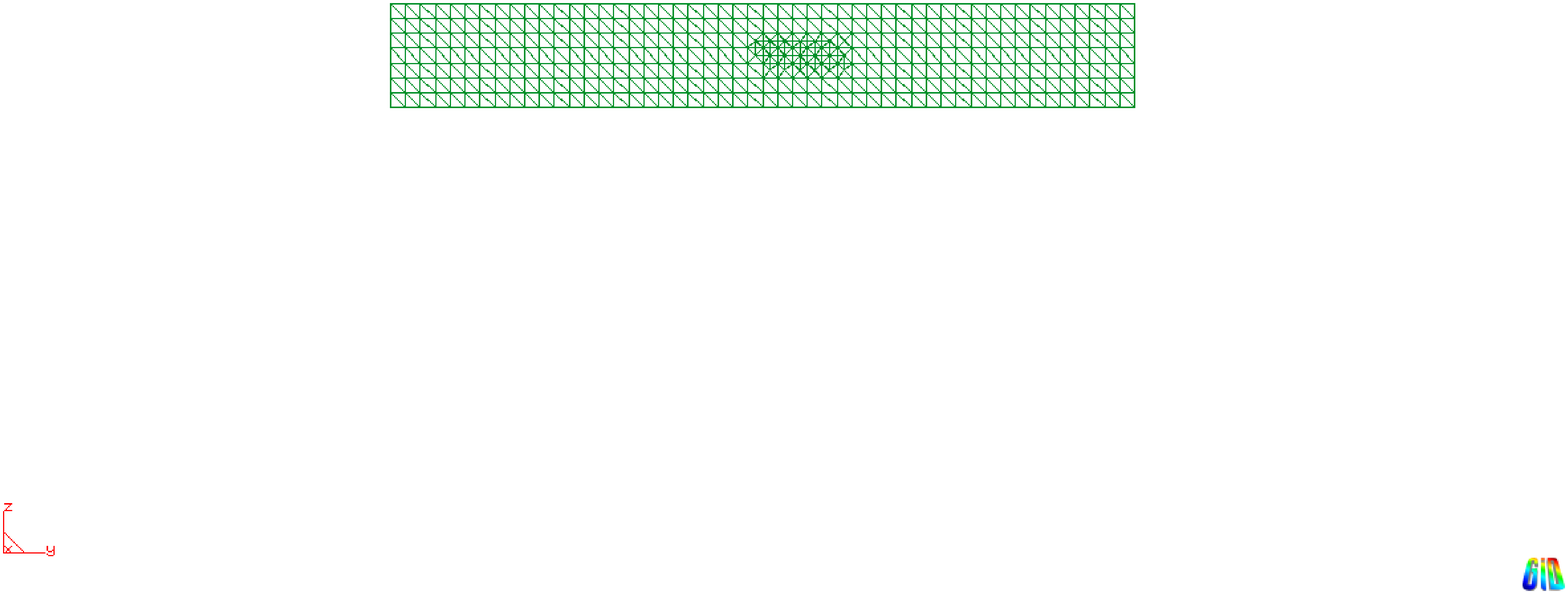}\\
      (d) $xy$-projection &
      (e) $xz$-projection &
      (f) $yz$-projection \\
      \includegraphics[width=5cm, clip=true, trim = 11cm 4cm 11cm 4cm]{16/ref3xy.eps}&
      \includegraphics[width=5cm, angle = 90, clip = true, trim = 11cm 12.5cm 11cm 12.5cm]{16/ref3xz.eps}&
      \includegraphics[width=5cm, angle = 90, clip = true, trim = 11cm 12.5cm 11cm 12.5cm]{16/ref3yz.eps}\\
      (g) $xy$-projection &
      (h) $xz$-projection &
      (i) $yz$-projection
    \end{tabular}
 \end{center}
    \caption{Adaptively refined meshes for target number 7 (doll, air
      inside) of Table \ref{tab:table1} used in our computations. (a)
      - (c) once refined, (d) - (f) twice refined, (g) - (i) three
      times refined.}
\label{fig:fig12}
\end{figure}

\begin{figure}
  \begin{center}
    \begin{tabular}{cccc}
      \includegraphics[width=5cm, clip=true, trim = 11cm 4cm 11cm 4cm]{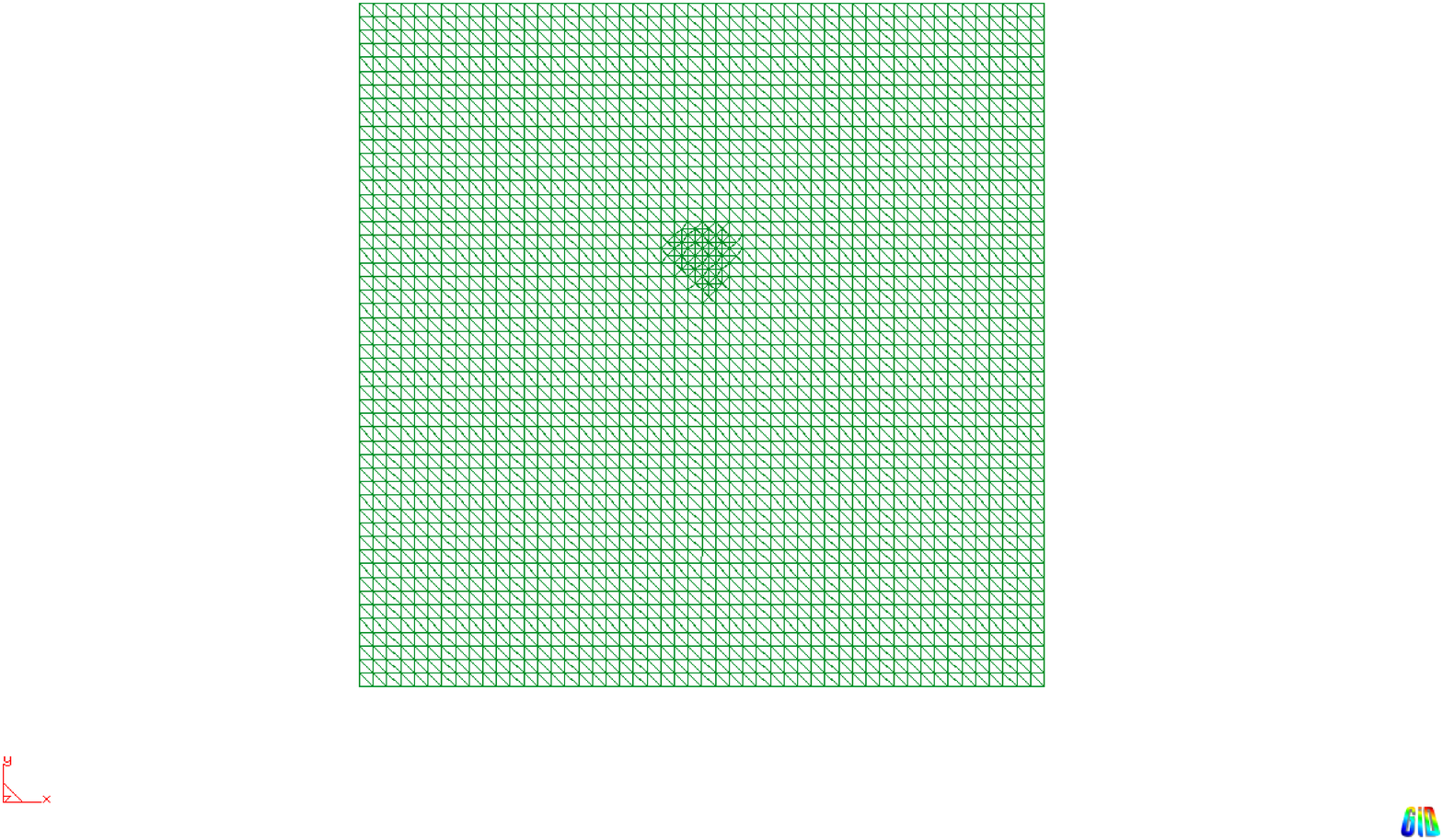}&
      \includegraphics[width=5cm, angle = 90, clip = true, trim = 11cm 12.5cm 11cm 12.5cm]{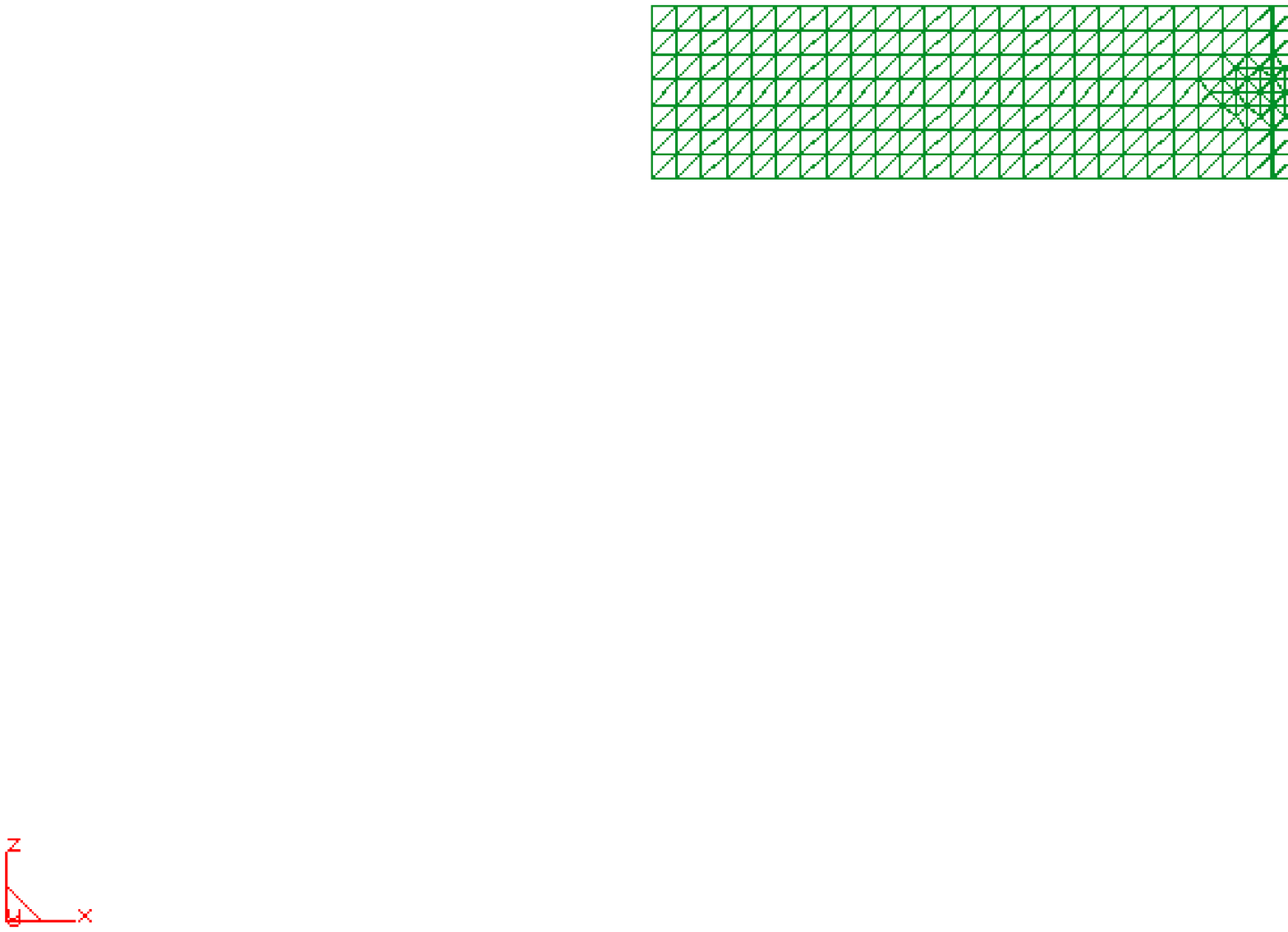}&
      \includegraphics[width=5cm, angle = 90, clip = true, trim = 11cm 12.5cm 11cm 12.5cm]{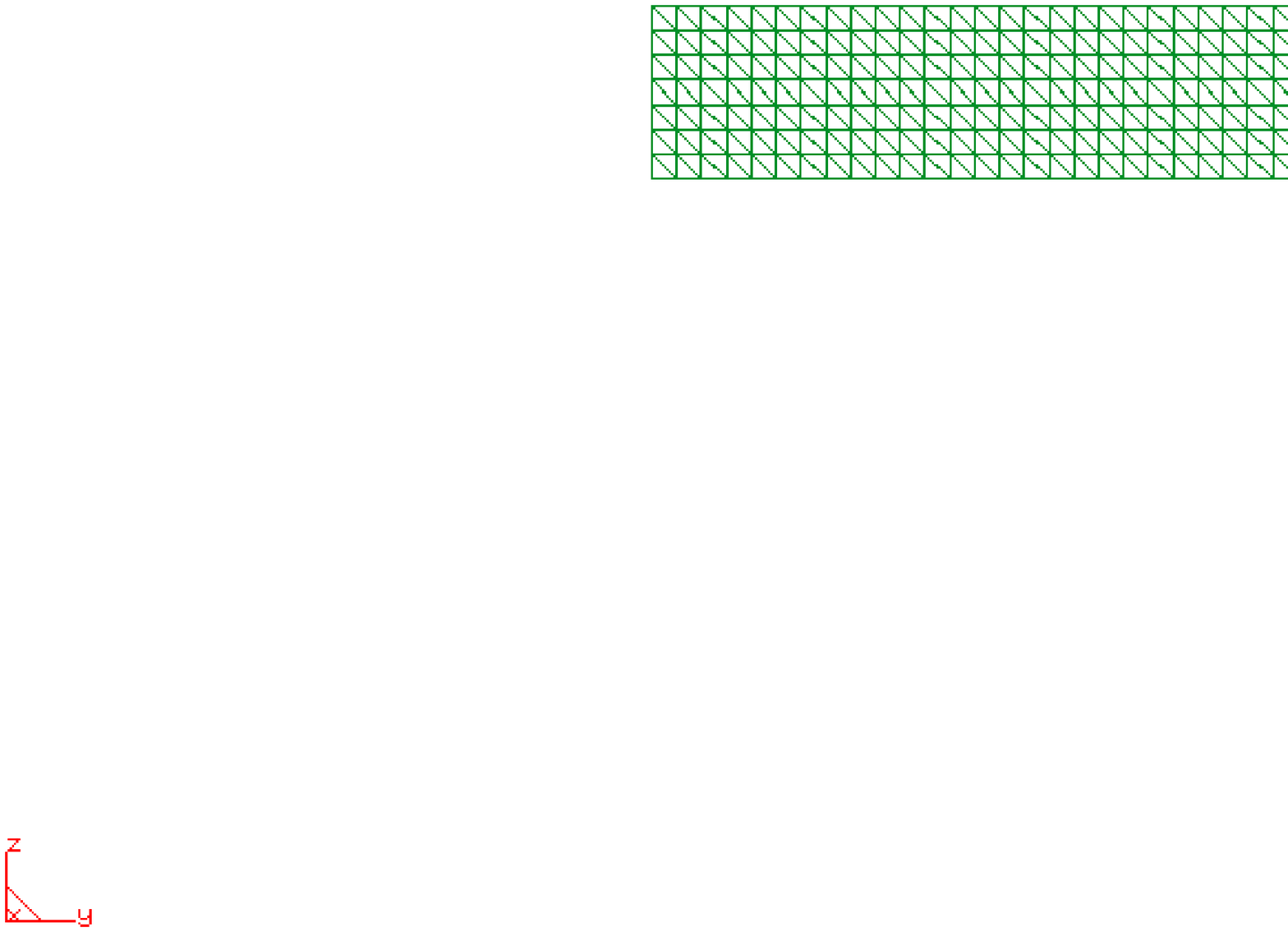}&
      \includegraphics[width=4.0cm, clip = true, trim = 1.6cm 0.0cm 4.0cm 0.0cm]{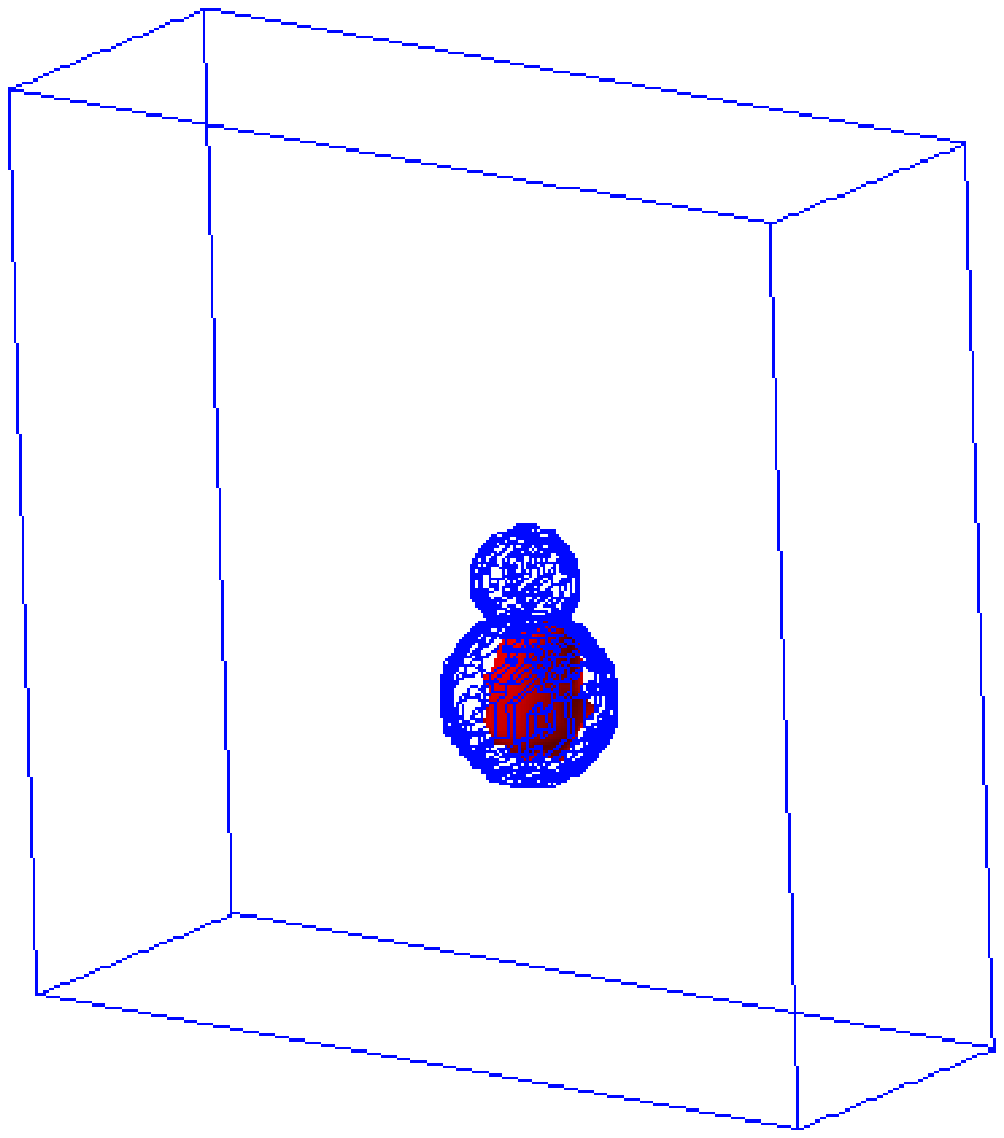}\\
      (a) &
      (b) &
      (c) &
      (d) 
    \end{tabular}
 \end{center}
    \caption{(a) $xy$-projection, (b) $xz$-projection, and (c) $yz$-projection of the twice  refined (optimal) mesh and the reconstruction (d) of target number  9 on that mesh.}
\label{fig:fig13}
\end{figure}

\begin{figure}
  \begin{center}
    \begin{tabular}{ccc}
      \includegraphics[width=4.0cm, clip = true, trim = 1.6cm 0.0cm 4.0cm 0.0cm]{18/eps_perspective.eps} &
      \includegraphics[width=4.0cm, clip = true, trim = 1.8cm 0.0cm 4.0cm 0.0cm]{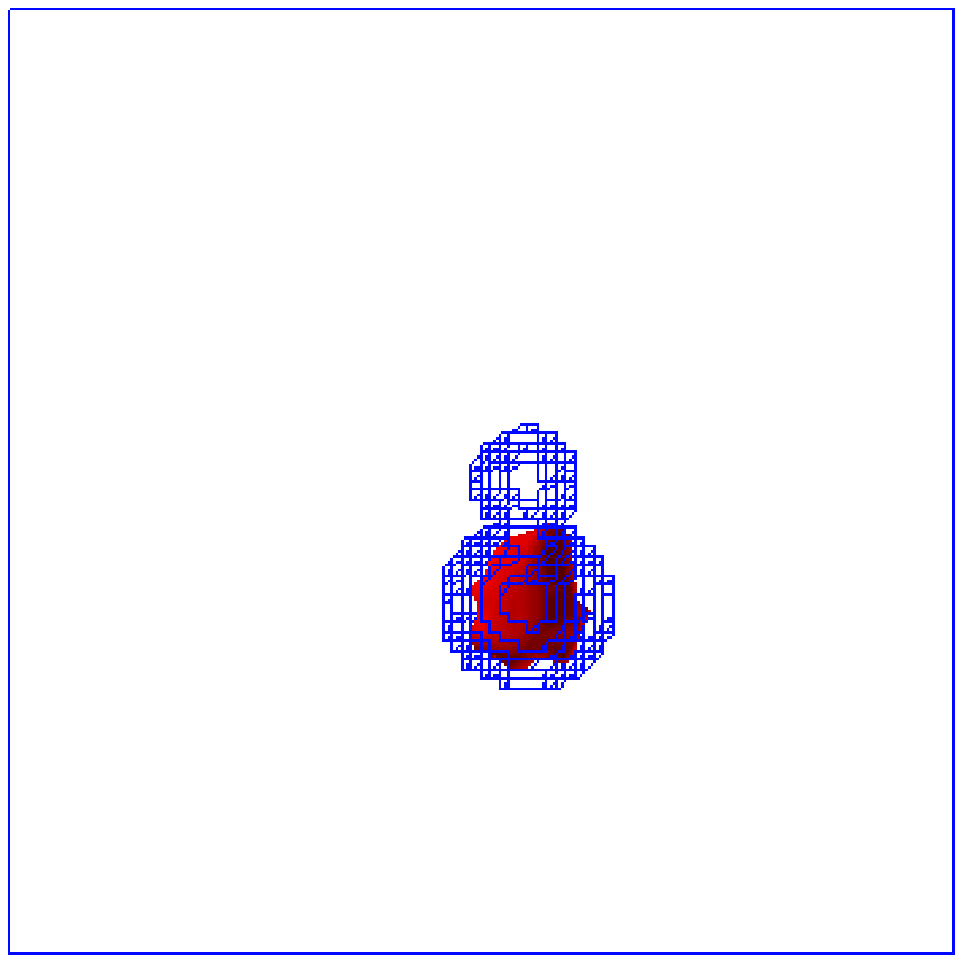} &
      \includegraphics[width=4.0cm, clip = true, trim = 1.8cm 0.0cm 4.0cm 0.0cm]{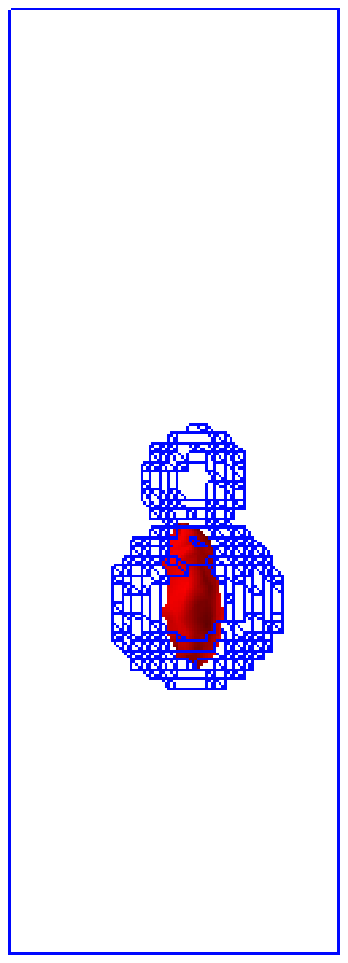}\\
      (a) Perspective view&
      (b) Front view&
      (c) Side view\\
      \includegraphics[width=4.0cm, clip = true, trim = 2.0cm 0.0cm 4.0cm 0.0cm]{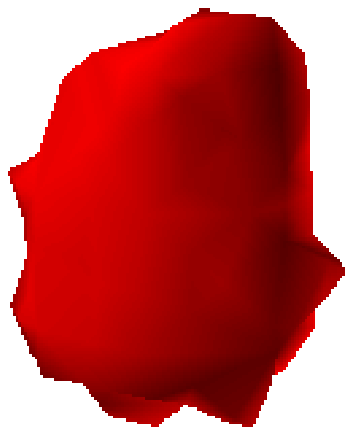} &
      \includegraphics[width=4.0cm, clip = true, trim = 2.0cm 0.0cm 4.0cm 0.0cm]{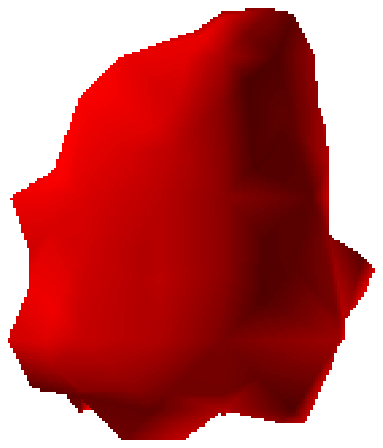} &
      \includegraphics[width=4.0cm, clip = true, trim = 2.0cm 0.0cm 4.0cm 0.0cm]{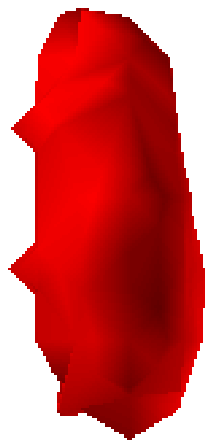}\\
      (d) Zoom, perspective&
      (e) Zoom, front&
      (f) Zoom, side
    \end{tabular}
 \end{center}
    \caption{Three views and zooms of target number 9 of Table \ref{tab:table1} on twice refined
      mesh. Thin lines indicate correct shape. Comparison of d) with Figure \ref{fig:globconv}-f) again shows a significant     improvement of the image due to the adaptivity. }
\label{fig:fig14}
\end{figure}

\begin{center}
\textbf{Acknowledgments}
\end{center}

This research was supported by US Army Research Laboratory and US Army
Research Office grant
 W911NF-11-1-0399, the
Swedish Research Council, the Swedish Foundation for Strategic
Research (SSF) through the Gothenburg Mathematical Modelling Centre
(GMMC). The computations were performed on
resources at Chalmers Centre for Computational Science and Engineering
(C3SE) provided by the Swedish National Infrastructure for Computing
(SNIC).


\begin{thebibliography}{99}
\bibitem{Ass} F.~Assous, P.~Degond, E.~Heintze and P.~Raviart, On a
finite-element method for solving the three-dimensional Maxwell equations, 
\emph{J. Comput.Physics}, 109, 222--237, 1993.

\bibitem{BakKok} A.B.\ Bakushinsky and M.Yu. Kokurin, \emph{Iterative
Methods for Approximate Solutions of Inverse Problems}, Springer, New York,
2004.

\bibitem{BKrelax} L. Beilina and M.V. Klibanov, Relaxation property
  for the adaptivity for ill-posed problems, \emph{Applicable
    Analysis}, 93, 223-253, 2014.

\bibitem{BK1} L. Beilina and M.V. Klibanov, A globally convergent
  numerical method for a coefficient inverse problem, \emph{SIAM
    J. Scientific Computing}, 31, 478-509, 2008.

\bibitem{BK2} L. Beilina and M.V. Klibanov, Synthesis of global convergence
and adaptivity \ for a hyperbolic coefficient inverse problem in 3D, \emph{%
J.\ Inverse and Ill-Posed Problems}, 18, 85-132, 2010.

\bibitem{BK} L. Beilina and M.V. Klibanov, \emph{Approximate Global
Convergence and Adaptivity for Coefficient Inverse Problems}, Springer, New
York, 2012.

\bibitem{BKJIIP12} L. Beilina and M.V. Klibanov, A new approximate
mathematical model for global convergence for a coefficient inverse problem
with backscattering data, \emph{J. Inverse and Ill-Posed Problems}, 20,
513-565, 2012.

\bibitem{BTKF} L. Beilina, Nguyen Trung Th\`{a}nh, M. V. Klibanov and M. A.
Fiddy, Reconstruction from blind experimental data for an inverse problem
for a hyperbolic equation, \emph{Inverse Problems} 30, 025002, 2014.

\bibitem{BM} L.\ Beilina, Energy estimates and numerical verification of the
stabilized domain decomposition finite element/finite difference approach
for the Maxwell's system in time domain, \emph{Central European Journal of
Mathematics}, 11, 702-733, 2013.

\bibitem{BMaxwell2} L. Beilina, Adaptive Finite Element Method for a
coefficient inverse problem for the Maxwell's system,\emph{\ Applicable
Analysis}, 90, 1461-1479, 2011.

\bibitem{BJ} L. Beilina and C. Johnson, A posteriori error estimation in
computational inverse scattering, \emph{Mathematical Models in Applied
Sciences}, 1, 23-35, 2005.

\bibitem{Brenner} S.~C.~Brenner and L.~R.~Scott, \emph{The Mathematical
theory of finite element methods}, Springer-Verlag, Berlin, 1994.

\bibitem{EM} B. Engquist and A. Majda, Absorbing boundary conditions for the
numerical simulation of waves, \emph{\ Math. Comp.,} 31, 629-651, 1977.

\bibitem{Isakov} V.\ Isakov, Inverse obstacle problems, \emph{Inverse
Problems}, 25, 123002, 2009.


\bibitem{joly} P.~Joly, \emph{Variational Methods for Time-Dependent Wave
Propagation Problems}, Lecture Notes in Computational Science and
Engineering, Springer, New York, 2003.

\bibitem{KBKSNF} A.V. Kuzhuget, L. Beilina, M.V. Klibanov, A. Sullivan,\ L.
Nguyen and M.A. Fiddy, Blind experimental data collected in the field and an
approximately globally convergent inverse algorithm, \emph{Inverse Problems}%
, 28, 095007, 2012.


%

\bibitem{Lad} O.A.\ Ladyzhenskaya, \emph{Boundary Value Problems of
Mathematical Physics}, Springer, New York, 1985.

\bibitem{Lak1} A. Lakhal, A decoupling based imaging method for inverse
medium scattering for Maxwell's equations, \emph{Inverse Problems}, 26,
015007, 2010.

\bibitem{Lak2} A. Lakhal, KAIRUAIN-algorithm applied on electromagnetic
imaging, \emph{Inverse Problems}, 29, 095001, 2013. 

\bibitem{Li} J.\ Li, J. Xie and J. Zou, An adaptive finite element
reconstruction of distributed fluxes, \emph{Inverse Problems}, 27, 075009,
2011.

\bibitem{Liu1} J. Li, H. Y. Liu, H. Sun and J. Zou, Reconstructing acoustic
obstacles by planar and cylindrical waves, \emph{J. Math. Phys}., 53,
103705, 2012.

\bibitem{Liu2} J. Li, H. Y. Liu, H. Sun and J. Zou, Imaging acoustic
obstacles by singular and hypersingular point sources, \emph{Inverse
Problems and Imaging}, 7, 545--563, 2013. 

\bibitem{Liu} Y.~Liu, J.~Su, Z.-J. Lin, S.~Teng, A.~Rhoden, N.~Pantong, and
H.~Liu. \newblock Reconstructions for continuous-wave diffuse optical
tomography by a globally convergent method. \newblock 2013. \newblock %
Preprint, available online at http://www.ma.utexas.edu/mp\_arc/, preprint
number 13 - 87.

\bibitem{Peron} O.Pironneau, \emph{Optimal shape design for elliptic systems}%
, Springer-Verlag, Berlin, 1984.

\bibitem{Stolt:1978} R.~Stolt, Migration by Fourier transform, \newblock
\emph{Geophysics}, 43, 23--48, 1978.

\bibitem{NBKF} Nguyen Trung Th\`{a}nh, L.~Beilina, M.~V. Klibanov and M.~A.
Fiddy, Reconstruction of the refractive index from experimental
backscattering data using a globally convergent inverse method, \emph{SIAM
J. Scientific Computing}, accepted for publication; preprint: \emph{arXiv}%
:1306.3150 [math-ph], 2013.

\bibitem{TBKF:2013-2}
N.~T. Th\`anh, L.~Beilina, M.~V. Klibanov, and M.~A. Fiddy.
\newblock Imaging of buried objects from experimental backscattering radar
  measurements using a globally convergent inverse method.
\newblock 2014.
\newblock Submitted.


\bibitem{T1} M. Sini and Nguyen Trung Th\`{a}nh, Inverse acoustic obstacle
scattering problems using multifrequency measurements, \emph{Inverse
Problems and Imaging}, 6, 749--773, 2012.

\bibitem{T2} M. Sini and Nguyen Trung Th\`{a}nh, Convergence rates of
recursive Newton-type methods for multifrequency scattering, \emph{arXiv}%
:1310.5156 [math.NA], 2013.

\bibitem{waves} WavES, the software package, http://www.waves24.com

\bibitem{Yilmaz:1987} O.~Yilmaz. \newblock {\em Seismic Data Imaging}. %
\newblock Society of Exploration Geophysicists, Tulsa Oklahoma, 1987.

\end{thebibliography}
\end{document}